\documentclass[aap,seceqn,MSNbibl,citesort,dvips]{arximspdf}

\doi{10.1214/10-AAP710}
\volume{21}
\issue{2}
\pubyear{2011}
\firstpage{397}
\lastpage{463}

\makeatletter

\newtheorem{theorem}{Theorem}[section]
\newproclaim{definition}[theorem]{Definition}
\newproclaim{example}[theorem]{Example}
\newtheorem{lemma}[theorem]{Lemma}
\newtheorem{proposition}[theorem]{Proposition}
\newproclaim{remark}[theorem]{Remark}

\newcommand{\re}{{\mathbb R}}
\newcommand{\E}{{\mathbb E}}
\newcommand{\tr}{\operatorname{Tr}}
\newcommand{\supp}{\operatorname{supp}}
\newcommand{\diag}{\operatorname{diag}}
\newcommand{\rk}{\operatorname{rk}}
\newcommand{\sgn}{\operatorname{sgn}}
\newcommand{\esssup}{\operatorname{essup}}
\newcommand{\Dya}{\operatorname{Dya}}
\newcommand{\vecnew}{\operatorname{vec}}

\newcommand{\im}{{i}}
\newcommand{\R}{{\mathbb R}}
\newcommand{\D}{{\mathbb D}}
\newcommand{\C}{{\mathbb C}}
\newcommand{\N}{{\mathbb N}}
\newcommand{\Pa}{{\mathbb P}}
\newcommand{\Acal}{{\mathcal A}}
\newcommand{\Bcal}{{\mathcal B}}
\newcommand{\Dcal}{{\mathcal D}}
\newcommand{\Fcal}{{\mathcal F}}
\newcommand{\Gcal}{{\mathcal G}}
\newcommand{\Hcal}{{\mathcal H}}
\newcommand{\Scal}{{\mathcal S}}

\makeatother

\begin{document}
\begin{frontmatter}

\title{Affine processes on positive semidefinite matrices}
\runtitle{Affine processes on positive semidefinite matrices}

\begin{aug}
\author[A]{\fnms{Christa} \snm{Cuchiero}\thanksref{a1}\ead
[label=e1]{christa.cuchiero@math.ethz.ch}},
\author[B]{\fnms{Damir} \snm{Filipovi\'c}\thanksref{a2}\ead
[label=e2]{damir.filipovic@epfl.ch}},
\author[C]{\fnms{Eberhard} \snm{Mayerhofer}\thanksref{a2}\ead
[label=e3]{eberhard.mayerhofer@vif.ac.at}}\\ and
\author[A]{\fnms{Josef} \snm{Teichmann}\corref{}\thanksref{a1}\ead
[label=e4]{josef.teichmann@math.ethz.ch}}
\runauthor{Cuchiero, Filipovi\'c, Mayerhofer and Teichmann}
\affiliation{ETH Z\"urich, \'Ecole Polytechnique F\'ed\'erale de
Lausanne, Vienna Institute of Finance and ETH Z\"urich}
\address[A]{C. Cuchiero\\
J. Teichmann\\
ETH Z\"urich\\
Departement Mathematik\\
R\"amistrasse 101, 8092 Z\"urich\\
Switzerland\\
\printead{e1}\\
\phantom{E-mail: }\printead*{e4}}
\address[B]{D. Filipovi\'c\\
\'Ecole Polytechnique F\'ed\'erale\\
de Lausanne\\
Swiss Finance Institute\\
Quartier UNIL-Dorigny\\
Extranef Building, 1015 Lausanne\\
Switzerland\\
\printead{e2}}
\address[C]{E. Mayerhofer\\
Vienna Institute of Finance\\
Heiligenst\"adter Strasse 46--48\\
1190 Vienna\\
Austria\\
\printead{e3}}
\end{aug}
\thankstext{a1}{Supported by FWF-Grant Y328 (START prize from the Austrian Science Fund).}
\thankstext{a2}{Supported by WWTF (Vienna Science and Technology Fund) and Swissquote.}

% HISTORY:
\received{\smonth{10} \syear{2009}}
\revised{\smonth{4} \syear{2010}}

% ABSTRACT
%
\begin{abstract}
This article provides the mathematical foundation for stochastically
continuous affine processes on the cone of positive semidefinite
symmetric matrices. This analysis has been motivated by a large and
growing use of matrix-valued affine processes in finance,
including multi-asset option pricing with stochastic volatility and
correlation structures, and fixed-income models with stochastically
correlated risk factors and default intensities.
\end{abstract}

% KEYWORDS
%

\end{frontmatter}

\tableofcontents
\setcounter{footnote}{2}

%s1 ###
\section{Introduction}

This paper provides the mathematical foundation for stochastically
continuous affine processes on the cone of positive semidefinite
symmetric $d\times d$-matrices $S_d^+$. These matrix-valued affine
processes have arisen from a large and growing range of useful
applications in finance, including multi-asset option pricing with
stochastic volatility and correlation structures, and fixed-income
models with stochastically correlated risk factors and default
intensities.

For illustration, let us consider a multi-variate stochastic
volatility model consisting of a $d$-dimensional logarithmic price
process with risk-neutral dynamics
%
%e1.1 ###
%
\begin{equation}\label{stockex}
dY_t=\bigl(r\mathbf{1}-\tfrac{1}{2}X^{\diag}_{t}\bigr)\,dt+\sqrt{X_t}\,dB_t,\qquad
Y_0=y,
\end{equation}
and stochastic covariation process $X=\langle Y,Y\rangle$, which is
a proxy for the instantaneous covariance of the price returns. Here
$B$ denotes a standard $d$-dimensional Brownian motion, $r$ the
constant interest rate, $\mathbf{1}$ the vector whose entries are
all equal to one and $X^{\diag}$ the vector containing the diagonal
entries of $X$.

The necessity to specify $X$ as a process in $S_d^+$ such that it
qualifies as
covariation process is one of the mathematically interesting and
demanding aspects of such models. Beyond that, the modeling of $X$
must allow for enough flexibility in order to reflect the stylized
facts of financial data and to adequately capture the dependence
structure of the different assets. If these requirements are met,
the model can be used as a basis for financial decision-making in
the area of portfolio optimization, pricing of multi-asset options
and hedging of correlation risk.

The tractability of such a model crucially depends on the dynamics
of $X$. A~large part of the literature in the area of multivariate
stochastic volatility modeling has proposed the following affine
dynamics for $X$:
%
%e1.2 ###
%
\begin{eqnarray}
\label{eq: cov process}\quad
dX_t&=&(b+H X_t+X_tH^{\top})\,dt+\sqrt{X_t}\,dW_t\,\Sigma+\Sigma^\top
\,dW_t^\top\,\sqrt{X_t}+ dJ_t,\nonumber\\[-8pt]\\[-8pt]
X_0&=&x\in S_d^+,
\nonumber
\end{eqnarray}
where $b$ is some suitably chosen matrix in $S_d^+$, $H, \Sigma$
some invertible matrices, $W$ a standard $d\times d$-matrix of
Brownian motions possibly correlated with $B$, and $J$ a pure jump
process whose compensator is an affine function of
$X$.\footnote{This affine multi-variate stochastic volatility model
generalizes the well-known one-dimensional models of
Heston \cite{heston}, for the diffusion case, or the
Barndorff-Nielsen Shepard model \cite{bns}, for the pure jump case.}

The main reason for the analytic tractability of this model is that,
under some technical conditions, the following affine transform
formula holds:
\[
\E_{x,y}\bigl[ e^{-\tr(zX_t)+v^\top Y_t}\bigr] =
e^{\Phi(t,z,v)+\tr(\Psi(t,z,v) x)+ v^\top y}
\]
for appropriate
arguments $z\in S_d\times\im S_d$ and $v\in\C^d$. The functions
$\Phi$ and $\Psi$ solve a system of nonlinear ordinary differential
equations (ODEs), which are determined by the model parameters.
Setting $v=0$, $\phi(t,z)=-\Phi(t,z,0)$ and $\psi(t,z)=-\Psi(t,z,0)$
and taking $z=u \in S_d^+$, we arrive at
%
%e1.3 ###
%
\begin{equation}\label{introaffine}
\E_x\bigl[ e^{-\tr(uX_t) }\bigr] =
e^{-\phi(t,u)-\tr(\psi(t,u) x)},\qquad u\in S_d^+.
\end{equation}

In this paper, we characterize the class of all stochastically
continuous time-homogeneous Markov processes with the key
property (\ref{introaffine})---henceforth called
affine processes---on $S_d^+$. Our main result shows that an affine
process is necessarily a Feller process whose generator has affine
coefficients in the state variables. The parameters of the generator
satisfy some well-determined admissibility conditions, and are in a
one-to-one relation with those of the corresponding ODEs for $\phi$
and $\psi$. Conversely, and more importantly for applications, we
show that for any admissible parameter set there exists a unique
well-behaved affine process on $S_d^+$. Furthermore, we prove that
any stochastically continuous infinitely decomposable Markov process
on $S_d^+$ is affine with zero diffusion, and vice versa.

On the one hand, our findings extend the model class (\ref{eq: cov
process}), since a more general drift and jumps are possible. Indeed,
we allow for full generality in $b$, as long as
$b-(d-1)\Sigma^T\Sigma\in S_d^+$, for a general linear drift part
$B(x)=\sum_{ij} x_{ij}\beta^{ij}$ and for an inclusion of (infinite
activity) jumps. This of course enables more flexibility in
financial modeling. For example, due to the general linear drift
part, the volatility of one asset can generally depend on the other
ones, which is not possible for $B(x)=Hx+xH^{\top}$. On the other
hand, we now know the exact assumptions under which affine processes
on $S_d^+$ actually exist. Our characterization of affine processes
on $S_d^+$ is thus exhaustive. Beyond that, the equivalence of
infinitely decomposable Markov processes with state space $S_d^+$
and affine processes without diffusion is interesting in its own
right.

This paper complements Duffie, Filipovi{\'c} and Schachermayer \cite{dfs}, who analyzed
time-homogeneous affine processes on the state space
$\R^m_+\times\R^n$.\footnote{For the diffusion case see also
\cite{ADPTA} or \cite{fil09}, Chapter 10. Time-inhomogeneous affine
processes on $\R^m_+\times\R^n$ have been explored in\vspace*{1pt}
\cite{filaff05}.} Matrix-valued affine processes seem to have been
studied systematically for the first time in the literature by
Bru \mbox{\cite{bru89,bru}}, who introduced the so called Wishart
processes. These are generalizations of squares of matrix
Ornstein--Uhlenbeck processes, that is, of the form (\ref{eq: cov
process}) for $J=0$ and $b=k\Sigma^\top\Sigma$, for some real
parameter $k>d-1$. Note that $k>d-1$ is a stronger assumption than
what we require on $b$ and $\Sigma^\top\Sigma$. Bru \cite{bru} then
establishes existence and uniqueness of a local\footnote{Up to the
first collision time of the eigenvalues.} $S_d^+$-valued solution to
(\ref{eq: cov process}) under the additional assumptions that $X_0$
has distinct eigenvalues, $-H\in S_d^+$, and that $H$ and $\Sigma$
commute (see \cite{bru}, Theorem 2$''$). In the more special case
where $H=0$ and $k > d-1$, Bru \cite{bru} shows global existence and
uniqueness for (\ref{eq: cov process}) for any $X_0$ with distinct
eigenvalues (see \cite{bru}, Theorem 2 and last part of Section
3).\footnote{Actually, Bru \cite{bru} establishes existence
and uniqueness of solutions also for $k=1,\ldots,d-1$. But these are
degenerate solutions, as they are only defined on lower-dimensional
subsets of the boundary of $S_d^+$ (see \cite{bru}, Corollary 1).}
Bru's results concerning strong solutions
have recently been extended to the case of matrix valued
jump-diffusions; see \cite{pfaffel}.

Wishart processes have subsequently been introduced in the financial
literature by Gourieroux and Sufana
\cite{gourierouxsufana,gousuf04} and Gourieroux et al.
\cite{gourieerouxmonfortsufana}. Financial applications thereof have
then been taken up and carried further by various authors, including
Da Fonseca et al.
\cite
{fonsecagrasselliielpo1,fonsecagrasselliielpo2,fonsecaetal2,fonsecaetal1}
and Buraschi, Cieslak and Trojani \cite{burcietro07,buraschiporchiatrojani}. Grasselli and
Tebaldi \cite{grasselli} give some general results on the solvability
of the corresponding Riccati ODEs. Barndorff-Nielsen
and Stelzer \cite{barndorffstelzer} provide a theory for a certain
class of matrix-valued L\'evy driven Ornstein--Uhlenbeck processes of
finite variation. Leippold and Trojani \cite{leippoldtrojani}
introduce $S_d^+$-valued affine jump diffusions and provide
financial examples, including multi-variate option pricing,
fixed-income models and dynamic portfolio choice. All of these
models are contained in our framework.

We want to point out that the full characterization of positive
semidefinite matrix-valued affine processes needs a multitude of
methods. In order to prove the fundamental property of regularity of
affine processes another adaption of the famous analysis of
Montgommery and Zippin is necessary, which has been worked out
in \cite{keller} and \cite{kst} for the state space $\R^m_+\times
\R^n$. For the necessary conditions on drift, diffusion and jump
parameters we need the theory of infinitely divisible distributions
on $S_d^+$. Most interestingly, the constant drift part $b$ must
satisfy a condition depending on the magnitude of the diffusion
component (see Proposition \ref{th: constant drift}), which is in
accordance with the choice of the drift in Bru's work \cite{bru} on
Wishart processes, as explained above. This enigmatic additional
condition on the drift $b$ is derived by studying the process with
respect to well chosen test functions, including in our case the
determinant of the process. It is worth noting, as already visible in
dimension one, that a naive application of classical geometric
invariance conditions does not bring the correct necessary result on
the drift but a stronger one. Indeed, take a one-dimensional affine
diffusion process $X$ solving
\[
dX_t = b \,dt + \sqrt{X_t}\,dW_t.
\]
Then a back-of-the-envelope calculation would yield the Stratonovich
drift at the boundary point $ x = 0 $ of value $ b - \frac{1}{4} $,
leading to the necessary parameter restriction $ b \geq\frac{1}{4}
$, which is indeed too strong. It is well known that the correct parameter
restriction is $ b \geq0 $. We see two reasons why geometric
conditions on the drift cannot be applied: first, precisely at the
boundary of our state spaces the diffusion coefficients are not
Lipschitz continuous anymore, and, second, the boundary of the cone of
positive semi-definite matrices is not a smooth submanifold but a more
complicated object.

For the sufficient direction refined methods from stochastic invariance theory
are applied. Having established viability of a particular class of
jump-diffusions on $S_d^+$,
existence of affine processes on $S_d^+$---under the necessary
parameter conditions---is shown through the solution of a
martingale problem. Uniqueness follows by semigroup methods which
need the theory of multi-dimensional Riccati equations.

Summing up, we face two major problems in the analysis of positive
matrix valued affine processes. First, the candidate stochastic
differential equations necessarily lead to volatility terms which
are not Lipschitz continuous at the boundary of the state space.
This makes every existence, uniqueness and invariance question
delicate. Second, the jump behavior
transversal to the boundary is of finite
total variation.

%s1.1 ###
\subsection{Program of the article}

For affine processes on $S_d^+$, results and proofs deviate in
essential points from the theory on state spaces of the form
$\re^m_+\times\re^n$ given in \cite{dfs,keller}, which is a
consequence of the more involved geometry of this nonpolyhedral
cone. The program of the paper as outlined below therefore includes
a comparison with the approach in \cite{dfs}.

Section \ref{section: Definition of affine property} contains the
main definition and a summary of the results of this article. In
Section \ref{subsec: regular and feller}, we then derive two main
properties, namely the regularity of the process and the Feller
property of the associated semigroup. The Feller property, in turn,
is a simple consequence of an important positivity result of the
characteristic exponents $\phi,\psi$, which is proved in
Lemma \ref{prop: feller prop}. This lemma is further employed as a
tool for the treatment of the generalized Riccati differential
equations in Section~\ref{section: Riccati} (see proof of
Proposition \ref{prop_ricc_sol}). The global existence and
uniqueness of these equations is then used to show uniqueness of the
martingale problem for affine processes (see proof of
Proposition \ref{th: existence Markov}).

In Section \ref{section: Nec paramter restrictions}, we define a set
of admissible parameters specifying the infinitesimal generator of
affine semigroups and prove the necessity of the parameter
restrictions (see Proposition \ref{th: necessary admissibility}).

The sufficient direction is then treated in Section \ref{sec:
existence}. It is known that, for $d\ge2$, there exist continuous
affine processes on $S_d^+$ which are---in contrast to those on the
state space $\R^m_+\times\R^n$---\textit{not} infinitely divisible (see
Example \ref{counterex: Bru}). The analysis of this paper reveals
the failure of infinite divisibility as a consequence of the drift
condition (see proof of Theorem \ref{th: char infdec infdiv}). This
has substantial influence on the approach chosen here to prove
existence of affine processes associated with a given parameter set:
Being in general hindered to recognize the solutions of the
generalized Riccati differential equations as
cumulant generating functions of sub-stochastic measures, as done
in \cite{dfs}, Section 7, we solve the martingale problem for the
associated L\'evy type generator on $S_d^+$, as exposed in
Section \ref{sec: existence} and Appendix \ref{app: viability}. In
Section \ref{subsec: alternative existence}, however, we deliver
a variant of the existence proof of \cite{dfs} for pure jump
processes, which is possible in this case due to the absence of a
diffusion component.

Finally, Section \ref{secproofs} contains the proofs of the main
results which build on the propositions of the previous sections.

%s1.2 ###
\subsection{Notation}\label{subsec: notation}
For the stochastic background and notation, we refer to standard text
books such as \cite{jacod} and \cite{revuzyor}. We write
$\re_+=[0,\infty)$ and $\re_{++}=(0,\infty)$. Moreover:
\begin{itemize}
\item$S_d$ denotes the space of symmetric $d\times d$-matrices
equipped with the scalar product $\langle x, y\rangle= \tr(xy)$. Note
that $S_d$ is isomorphic, but not isometric, to the standard
Euclidean space $\re^{d(d+1)\slash2}$. We denote by $\{c^{ij},i \leq
j\}$ the standard basis of $S_d$, that is, the $(kl)$th
component of $c^{ij}$ is given by
$c^{ij}_{kl}=\delta_{ik}\delta_{jl}+\delta_{jk}\delta_{il}(1-\delta_{ij})$,
where $\delta_{ij}$ denotes the Kronecker delta. Additionally, we
sometimes consider the following basis elements $\{e^{ij}, i
\leq j\}$ which are positive semidefinite and form a basis of
$S_d$:
\[
e^{ij}=\cases{c^{ii}, &\quad if $i=j$,\cr
c^{ii}+c^{ij}+c^{jj}, &\quad if $i \neq j$.}
\]
\item$S_d^+$ stands for the cone of symmetric $d\times d$-positive
semidefinite matrices, $S_d^{++}$~for its interior in $S_d$, the cone
of strictly
positive definite matrices. The boundary is denoted by $\partial
S_d^+=S_d^+\setminus S_d^{++}$, the complement is denoted by
$(S_d^+)^c$, and $S_d^+ \cup\{\Delta\}$ denotes the one-point
compactification. Recall that $S_d^+$ is self-dual [w.r.t. the
scalar product $\langle x,y\rangle=\tr(xy)$], that is,
\[
S_d^+ =\{ x\in S_d\mid\langle x,y\rangle\ge0 ,\forall y\in
S_d^+\}.
\]
Both cones, $S_d^+$ and $S_d^{++}$, induce a partial and strict
order relation on $S_d$, respectively: we write $x \preceq y$ if
$y-x\in S_d^+$, and $x\prec y$ if $y-x\in S_d^{++}$.

\item$M_d$ is the space of $d\times d$-matrices and $O(d)$ the
orthogonal group of dimension $d$ over $\re$.
\item$I_d$ denotes the $d\times d$-identity matrix.
\end{itemize}
Throughout this paper, a function $f\dvtx S_d\rightarrow\mathbb R$ is
understood as the restriction $f=g|_{S_d}$ of a function $g\dvtx
M_d\rightarrow\mathbb R$ which satisfies $g(x)=g(x^\top)$ for all
$x\in M_d$. Without loss of generality $g(x)=f((x+x^\top)/2)$. We
avoid using the $vech$ operator, that is, to identify $x\in S_d$
with a vector in $\re^{d(d+1)\slash2}$ by stringing the columns of
$x$ together, while only taking the entries $x_{ij}$ with $i\leq j$.

Throughout this article, we shall consider the following function
spaces for measurable $U \subseteq S_d$. We write $\Bcal(U)$ for the
Borel $\sigma$-algebra on $U$. $bU$ corresponds to the Banach space
of bounded real-valued Borel measurable functions $f$ on $U$ with
norm $\|f\|_{\infty}=\sup_{x\in U}|f(x)|$. We write $C(U)$ for the
space of real-valued continuous functions $f$ on $U$, $C_b(U)$ for
$C(U) \cap bU$, $C_c(U)$ for the space of functions $f \in C(U)$
with compact support and $C_0(U)$ for the Banach space of functions
$f \in C(U)$ with $\lim_{x \rightarrow\Delta}f(x)=0$ and norm
$\|f\|_{\infty}={\sup_{x\in U}}|f(x)|$. Furthermore,
$C^k(U)$ is the space of $k$ times differentiable functions $f$ on
$U^{\circ}$, the interior of $U$, such that all partial derivatives
of $f$ up to order $k$ belong to $C(U)$. As usual, we set
$C^\infty(U)=\bigcap_{k\ge1} C^k(U)$, and we write $C_c^k(U)=C_c(U)
\cap C^k(U)$ and $C^k_b(U)=C_b(U)\cap C^k(U)$, for $k\le\infty$.

%s2 ###
\section{Definition and characterization of affine processes}\label
{section: Definition of affine property}

We consider a time-homogeneous Markov process $X$ with state space
$S_d^+$ and semigroup $(P_t)_{t \geq0}$ acting on functions $f \in
bS_d^+$,
\[
P_tf(x)=\int_{S_d^+} f(\xi)p_t(x,d\xi),\qquad x \in S_d^+.
\]
We note that $X$ may not be conservative. Then there is a standard
extension of the transition probabilities to\vadjust{\goodbreak} the one-point
compactification $S_d^+\cup\{\Delta\}$ of $S_d^+$ by defining
\[
p_t(x,\{\Delta\})=1-p_t(x,S_d^+),\qquad
p_t(\Delta,\{\Delta\})=1
\]
for all $t$ and $x\in S_d^+$, with the
convention that $f(\Delta)=0$ for any function $f$ on $S_d^+$. Thus
$X$ becomes conservative on $S_d^+\cup\{\Delta\}$.

%Here is our main definition below.
%
\begin{definition}\label{def: affine process Sd+}
The Markov process $X$ is called \textit{affine} if:
\begin{longlist}[(ii)]
\item[(i)] it is stochastically continuous, that is, $\lim_{s\to t}
p_s(x,\cdot)=p_t(x,\cdot)$ weakly on $S_d^+$ for every $t$ and $x\in
S_d^+$, and

\item[(ii)] its Laplace transform has
exponential-affine
dependence on the initial state
%
%e2.1 ###
%
\begin{equation}\label{def: affine process}
P_te^{-\langle u, x\rangle}=\int_{S_d^+}e^{-\langle u, \xi
\rangle}p_t(x,d\xi)=e^{-\phi(t,u)-\langle\psi(t,u),x\rangle},
\end{equation}
for all $t$ and $u,x\in S_d^+$, for some functions $\phi\dvtx \re_+
\times S_d^+ \rightarrow\re_+$ and $\psi\dvtx \re_+ \times S_d^+
\rightarrow S_d^+$.
\end{longlist}
\end{definition}

Note that stochastic continuity of $X$ implies that $\phi(t,u)$ and
$\psi(t,u)$ are jointly continuous in $(t,u)$; see Lemma \ref{lem:
order preserving affine}(iii) below.
Moreover, due to the Markov property, this also means that
$p_t(x, \{\Delta\}) <1$ for all $x \in S_d^+$ and $t \geq0$. In contrast
to~\cite{dfs}, we take stochastic continuity as part of the
definition of affine processes, and consider the Laplace transform
instead of the characteristic function. The latter is justified by
the nonnegativity of $X$, the former is by convenience since, as we
will see in Proposition \ref{th: regularity} below, it automatically
implies regularity in the following sense.
\begin{definition}
The affine process $X$ is called \textit{regular} if the derivatives
%
%e2.2 ###
%
\begin{equation}\label{eq: F,R}
F(u) = \frac{\partial\phi(t,u)}{\partial t} \bigg|_{t=0+},\qquad
R(u) = \frac{\partial\psi(t,u)}{\partial t} \bigg|_{t=0+}
\end{equation}
exist and are continuous at $u=0$.
\end{definition}

We remark that there are simple examples of Markov processes which
satisfy Definition \ref{def: affine process Sd+}(ii) but are not stochastically continuous; see
\cite{dfs}, Remark 2.11. However, such processes are of limited
interest for applications and will not be considered.

In the following, we shall provide an equivalent characterization of
the affine property in terms of the generator of $X$. As we shall see
in (\ref{eq: generator}), the diffusion, drift, jump and killing
characteristics of $X$ depend in an affine way on the underlying
state. We denote by $\chi\dvtx S_d\rightarrow S_d$ some bounded continuous
truncation function with $\chi(\xi)=\xi$ in a neighborhood of $0$.
Then the involved parameters are admissible in the following
sense.
\begin{definition}\label{def: necessary admissibility}
An \textit{admissible parameter set} $(\alpha, b,\beta^{ij}, c,
\gamma, m, \mu)$ associated with $\chi$ consists of:\vadjust{\goodbreak}
\begin{itemize}
\item a linear diffusion coefficient
%
%e2.3 ###
%
\begin{equation}\label{eq: alpha}
\alpha\in S_d^+,
\end{equation}

\item a constant drift term
%
%e2.4 ###
%
\begin{equation}\label{eq: b}
b \succeq(d-1)\alpha,
\end{equation}
\item a constant killing rate term
%
%e2.5 ###
%
\begin{equation}\label{eq: c}
c \in\re^+,
\end{equation}

\item a linear killing rate coefficient
%
%e2.6 ###
%
\begin{equation}\label{eq: gamma}
\gamma\in S_d^+,
\end{equation}
\item a constant jump term: a
Borel measure $m$ on $S_d^+\setminus\{0\}$ satisfying
%
%e2.7 ###
%
\begin{equation}\label{eq: m}
\int_{S_d^+\setminus\{0\}}(\|\xi\|\wedge1)m(d\xi) <
\infty,
\end{equation}
\item a linear jump coefficient:
a $d\times d$-matrix $\mu=(\mu_{ij})$ of finite signed measures on
$S_d^+\setminus\{0\}$ such that $\mu(E)\in S_d^+$ for all $E\in
\mathcal{B}(S_d^+\setminus\{0\})$ and the kernel
%
%e2.8 ###
%
\begin{equation}\label{eq: mudef}
M(x,d\xi):=\frac{\langle x, \mu(d\xi)
\rangle}{\|\xi\|^2 \wedge1}
\end{equation}
satisfies
%
%e2.9 ###
%
\begin{equation}\label{eq: mu}
\int_{S_d^+\setminus\{0\}} \langle\chi(\xi), u\rangle M(x,d\xi) <
\infty\qquad\mbox{for all $x,u\in S_d^+$ with $\langle x,
u\rangle=0$,}
\end{equation}
\item a linear drift coefficient: a family $\beta^{ij}=\beta^{ji}\in
S_d$ such
that the linear map $B\dvtx S_d\to S_d$ of the form
%
%e2.10 ###
%
\begin{equation}\label{eq: betaijdef}
B(x)=\sum_{i,j} \beta^{ij}x_{ij}
\end{equation}
satisfies
%
%e2.11 ###
%
\begin{eqnarray}\label{eq: betaij}
&&\langle B(x),u\rangle-
\int_{S_d^+\setminus\{0\}}\langle\chi(\xi),u\rangle M(x,
d \xi) \geq0 \nonumber\\[-8pt]\\[-8pt]
&&\eqntext{\mbox{for all $x,u\in S_d^+$ with $\langle x,
u\rangle=0$.}}
\end{eqnarray}
\end{itemize}
\end{definition}

We shall comment more on the admissibility conditions in
Section \ref{subsecdis} below. The following three theorems contain
the main results of this article. Their proofs are given in
Section \ref{secproofs}. First, we provide a characterization of
affine processes on $S_d^+$ in terms of the admissible parameter set
introduced in Definition \ref{def: necessary admissibility}. As for
the domain of the generator, we consider the space $\Scal_+$ of
rapidly decreasing $C^\infty$-functions on $S_d^+$, defined in
(\ref{defS+}) below. It is shown in Appendix \ref{secstoneweier} that
$e^{-\langle u,\cdot\rangle}\in\Scal_+$, for $u\in S_d^{++}$, as
well as $C^\infty_c(S_d^+)\subset\Scal_+$.
\begin{theorem}\label{th: main theorem}
Suppose $X$ is an affine process on $S_d^+$. Then $X$ is regular and
has the Feller property. Let $\mathcal{A}$ be its infinitesimal
generator on $C_0(S_d^+)$. Then $\Scal_+ \subset D(\mathcal{A})$ and
there exists\vadjust{\goodbreak} an admissible parameter set $(\alpha, b, \beta^{ij}, c,
\gamma, m, \mu)$ such that, for $f \in\Scal_+$,
%
%e2.12 ###
%
\begin{eqnarray}
\label{eq: generator}\qquad
\mathcal{A}f(x)&=&\frac{1}{2}\sum_{i,j,k,l}A_{ijkl}(x)\,
\frac{\partial^2
f(x)} {\partial x_{ij}\,\partial x_{kl}}\nonumber\\
&&{}+\sum
_{i,j}\bigl(b_{ij}+B_{ij}(x)\bigr)\,\frac{\partial f(x)}{\partial
x_{ij}}
-(c+\langle\gamma,x\rangle)f(x)\nonumber\\[-8pt]\\[-8pt]
&&{} +\int_{S_d^+\setminus\{0\}}
\bigl(f(x+\xi)-f(x)\bigr)m(d\xi)\nonumber\\
&&{} +\int_{S_d^+\setminus\{0\}}\bigl(f(x+\xi)-f(x)-\langle
\chi(\xi), \nabla f(x)\rangle\bigr)M(x,d\xi),\nonumber
\end{eqnarray}
where $B(x)$ is defined by (\ref{eq: betaijdef}), $M(x,d\xi)$
by (\ref{eq: mudef}) and
%
%e2.13 ###
%
\begin{equation}\label{eq: cijkl}
A_{ijkl}(x)=x_{ik}\alpha_{jl}+x_{il}\alpha_{jk}+x_{jk}\alpha
_{il}+x_{jl}\alpha_{ik}.
\end{equation}
Moreover, $\phi(t,u)$ and $\psi(t,u)$ in (\ref{def: affine process})
solve the generalized Riccati differential equations, for $u \in
S_d^+$,
%
%e2.15 ###
%e2.14 ###
%
\begin{eqnarray}
\label{eq: F-Riccati}
\frac{\partial\phi(t,u)}{\partial t}&=&F(\psi(t,u)),\qquad \phi(0,u)=0,
\\
\label{eq: R-Riccati}
\frac{\partial\psi(t,u)}{\partial t}&=&R(\psi(t,u)),\qquad
\psi(0,u)=u,
\end{eqnarray}
with
%
%e2.17 ###
%e2.16 ###
%
\begin{eqnarray}\qquad
\label{eq: F}
F(u)&=&\langle b, u\rangle+ c - \int_{S_d^+\setminus\{0\}
}\bigl(e^{-\langle u, \xi\rangle}-1\bigr)m(d\xi),\\
\label{eq: R}
R(u)&=&-2u\alpha u+B^\top(u)+\gamma\nonumber\\[-8pt]\\[-8pt]
&&{}-\int_{S_d^+\setminus\{0\}}
\biggl(\frac{e^{-\langle u, \xi\rangle}-1+\langle
\chi(\xi), u \rangle}{\|\xi\|^2 \wedge
1}\biggr)\mu(d \xi),\nonumber
\end{eqnarray}
where $B^\top_{ij}(u)=\langle\beta^{ij},u\rangle$.

Conversely, let $(\alpha, b,\beta^{ij}, c, \gamma, m, \mu)$ be an
admissible parameter set. Then there exists a unique affine process
on $S_d^+$ with infinitesimal generator (\ref{eq: generator})
and (\ref{def: affine process}) holds for all $(t,u) \in\re_+\times
S_d^+$, where $\phi(t,u)$ and $\psi(t,u)$ are given by (\ref{eq:
F-Riccati}) and (\ref{eq: R-Riccati}).
\end{theorem}
\begin{remark}\label{remcons}
It can be proved as in \cite{maysmi09} that $X$ is conservative if
and only if $c=0$ and $\psi(t,0)\equiv0$ is the only $S_d^+$-valued
local solution of (\ref{eq: R-Riccati}) for $u=0$. The latter
condition clearly requires that $\gamma=0$.

Hence, a sufficient condition for $X$ to be conservative is $c=0$
and $\gamma=0$ and
\[
\int_{{S_d^+}\cap\{\|\xi\|\geq
1\}}\|\xi\|\bigl(\mu^+_{ij}(d\xi)+\mu^-_{ij}(d\xi)\bigr)< \infty\qquad
\mbox{for all } 1\leq i \leq j\leq d,
\]
where $\mu_{ij}=\mu_{ij}^+-\mu_{ij}^-$ denotes the Jordan
decomposition of $\mu_{ij}$. Indeed, it can be shown similarly as in
\cite{dfs}, Section 9, that the latter property implies Lipschitz
continuity of $R(u)$ on $S_d$.
\end{remark}

Due to the Feller property, as established in Theorem \ref{th: main
theorem}, any affine process $X$ on $S_d^+$ admits a c\`adl\`ag
modification, still denoted by $X$ (see, e.g.,
\cite{revuzyor}, Chapter~III.2). It can and will thus be realized on
the space
$\Omega=\D(S_d^+\cup\{\Delta\})$ of c\`adl\`ag paths $\omega\dvtx\R
_+\to
S_d^+\cup\{\Delta\}$ with $\omega(s)=\Delta$ for $s>t$ whenever
$\omega(t-)=\Delta$ or $\omega(t)=\Delta$. For every $x\in S_d^+$,
we denote by $\Pa_x$ the law of $X$ given $X_0=x$ and by
$(\mathcal{F}^X_t)$ the natural filtration generated by $X_t$. We
also consider the usual augmentation
%
%e2.18 ###
%
\begin{equation} \label{eq: filt}
\mathcal{\widetilde{F}}_t:=\bigcap_{x \in S_d^+} \mathcal{F}^{(x)}_t
\end{equation}
of $(\mathcal{F}^X_t)$, where $(\mathcal{F}^{(x)}_t)$ is the
augmentation of $(\mathcal{F}^X_t)$ with respect\vspace*{2pt} to
$\mathbb{P}_x$. Then $(\mathcal{\widetilde{F}}_t)$ is right continuous
and $X$ is still a Markov process under $(\mathcal{\widetilde{F}}_t)$.
We shall now relate conservative affine processes to semimartingales,
where semimartingales are understood with respect to the stochastic
basis $(\Omega,\mathcal{\widetilde{F}},(\mathcal{\widetilde{F}})_t,
\mathbb{P}_x)$ for every $x$.
\begin{theorem}\label{thmsemim}
Let $X$ be a conservative affine process on $S_d^+$ and let
$(\alpha, b, \beta^{ij},c=0, \gamma=0, m, \mu)$ be the related
admissible parameter set associated with the truncation function
$\chi$. Then $X$ is a semimartingale whose characteristics $(B, A,
\nu)$ with respect to $\chi$ are given by
%
%e2.21 ###
%e2.20 ###
%e2.19 ###
%
\begin{eqnarray}
\label{eq: char c}
A_{t,ijkl}&=&\int_0^t A_{ijkl}(X_{s})\,ds,\\
\label{eq: char b}
B_t&=&\int_0^t \biggl(b+\int_{S_d^+\setminus\{0\}}\chi(\xi)m(d\xi
)+B(X_{s})\biggr)\,ds,\\
\label{eq: char nu}
\nu(dt,d\xi)&=&\bigl(m(d\xi)+M(X_{t},d\xi)\bigr)\,dt,
\end{eqnarray}
where $B(x)$ is given by (\ref{eq: betaijdef}), $A_{ijkl}(x)$
by (\ref{eq: cijkl}) and $M(x,d\xi)$ by (\ref{eq: mudef}).
Furthermore, there exists, possibly on an enlargement of the
probability space, a $d \times d$-matrix of standard Brownian
motions $W$ such that $X$ admits the following representation:
%
%e2.22 ###
%
\begin{eqnarray}\label{eq: BM rep}
X_t&=&x+B_t+\int_0^t\bigl(\sqrt{X_{s}}\,dW_s\,\Sigma+\Sigma^{\top}\,dW_s\,\sqrt
{X_{s}}\bigr)\nonumber\\
&&{} +\int_0^t\int_{{S_d^+}\setminus\{0\}}\chi(\xi)\bigl(\mu^X(ds,d\xi
)-\nu(ds,d\xi)\bigr)\\
&&{} +\int_0^t\int_{{S_d^+}\setminus\{0\}}\bigl(\xi-\chi(\xi)\bigr)\mu
^X(ds,d\xi),
\nonumber
\end{eqnarray}
where $\Sigma\in M_d$ satisfies $\Sigma^{\top}\Sigma=\alpha$ and
$\mu^X$ denotes the random measure associated with the jumps of $X$.

Hence, $X$ is continuous if and only if $m$ and $\mu$ vanish.
\end{theorem}

Let $\mathcal P$ be the set of all families of probability measures
$(\mathbb P_x)_{x\in S_d^+}$ on the canonical probability space
$(\Omega,\mathcal F^X)$ such that $(X,(\mathbb P_x)_{x\in S_d^+})$
is a stochastically continuous Markov processes on $S_d^+$ with
$\mathbb P_x[X_0=x]=1$, for all $x\in S_d^+$. Note that in
contrast to \cite{dfs}, there is no need to impose regularity of
$X$. For two probability measures $\mathbb P,\mathbb Q$ on
$(\Omega,\mathcal F^X)$, the convolution $\mathbb P*\mathbb Q$ is
defined as the push-forward of $\mathbb P\times\mathbb Q$ under the
map $(\omega,\omega')\mapsto
\omega+\omega'\dvtx (\Omega\times\Omega, \mathcal F^X \otimes\mathcal
F^X) \to(\Omega, \mathcal F^X )$.
\begin{definition}\label{def: infinite divisibility}
An element $ (\mathbb P_x)_{x\in S_d^+}\in\mathcal P$ is called:
\begin{longlist}[(ii)]
\item[(i)] infinitely decomposable, if for each $k\geq1$, there exists
$(\mathbb P^{(k)}_x)_{x\in S_d^+}\in\mathcal P$ such that
\[
\mathbb P_{x^{(1)}+\cdots+x^{(k)}}=\mathbb
P^{(k)}_{x^{(1)}}*\cdots*\mathbb P^{(k)}_{x^{(k)}};
\]
\item[(ii)] infinitely divisible, if the
one-dimensional marginal distributions $\mathbb P_x\circ X_t^{-1}$
are infinitely divisible, for all $(t,x)\in\re_+\times S_d^+$.
\end{longlist}
\end{definition}

In \cite{dfs} it was shown that regular affine processes on
$\R^m_+\times\R^n$ are infinitely decomposable Markov processes, and
vice versa. In fact, this property was the core for the existence
proof of affine processes in \cite{dfs}. On $S_d^+$ the situation is
different. The following counterexample reveals that not all
affine processes on $S_d^+$ are infinitely divisible.
\begin{example}\label{counterex: Bru}
The affine process $X$ on $S_d^+$ corresponding to the parameter set
$(\alpha=I_d,b=\delta I_d,0,0,0,0,0)$, where $\delta\in[d-1,\infty)$,
is the diffusion process initially studied by Bru \cite{bru}.
By \cite{bru}, Theorem 3, the Laplace-transforms
\[
\mathbb E_x\bigl[e^{-\langle
X_t,u\rangle}\bigr]=\bigl(\det(I+2tu)\bigr)^{-\delta/2}e^{-\langle
(I+2tu)^{-1}u,x\rangle}
\]
are those of noncentral Wishart distributions $\operatorname{WIS}(\delta,d,x)$. By
a well-known result due to Paul L\'evy, these Wishart distributions are
not infinitely divisible if $d \geq2$ (see \cite{donatimartin},
Section 2.C).
\end{example}

Here, is our main result on infinite divisibility of affine processes
on $S_d^+$.
\begin{theorem}\label{th: char infdec infdiv}
Let $d \geq 2$ and $(\mathbb P_x)_{x\in S_d^+}\in\mathcal P$. The following
assertions are
equivalent:
\begin{longlist}[(iii)]
\item[(i)] $(\mathbb P_x)_{x\in S_d^+}$ is infinitely
decomposable.
\item[(ii)] $(X,(\mathbb P_x)_{x\in S_d^+})$ is affine
with vanishing diffusion
parameter $\alpha=0$.
\item[(iii)] $(X,(\mathbb P_x)_{x\in S_d^+})$ is affine
and infinitely divisible.
\end{longlist}
\end{theorem}

%s2.1 ###
\subsection{Discussion of the parameters}\label{subsecdis}

We discuss and highlight some properties of the admissible
parameter set $(\alpha, b,\beta^{ij}, c, \gamma, m, \mu)$ of an
affine process.

Let us therefore define the normal cone
%
%e2.23 ###
%
\begin{equation}\label{eq: normal cone Sd+}
N_{S_d^+}(x)=\{u \in S_d^+ |\langle u, x \rangle=0\},
\end{equation}
containing the inward pointing normal vectors, to $S_d^+$ at $x\in
S_d^+$.\footnote{Indeed, we obtain (\ref{eq: normal cone Sd+}) from
the general definition in (\ref{defnorcone}) below by choosing $y=0$
and $y=2x$, and using the self-duality of $S_d^+\dvtx\langle
u,y\rangle\ge0$ for all $y,u\in S_d^+$.} It will be shown in
Lemma \ref{lem: zero divisor} below that $N_{S_d^+}(x)\neq\{0\}$
only for boundary elements $x\in\partial S_d^+$.

%s2.1.1 ###
\subsubsection{Diffusion} The diffusion term does not admit a constant
part, and its linear
part is of the very specific form
\[
\langle u, A(x) u\rangle= 4 \langle x, u\alpha u\rangle.
\]
This property of $A(x)$ has also been stated in the setting of
symmetric cones in~\cite{grasselli}. We could thus write the second order
differential operator in (\ref{eq: generator}) as
\[
\frac{1}{2}\sum_{i,j,k,l}A_{ijkl}(x)
\,\frac{\partial^2 f(x)} {\partial x_{ij}\,\partial x_{kl}} = 2\langle
x,\nabla\alpha\nabla f(x)\rangle.
\]
The reason why we introduce and
use the symmetrization (\ref{eq: cijkl}) of $A(x)$ is that it
corresponds to the quadratic characteristic (\ref{eq: char c}) of
the semimartingale $X$.

%s2.1.2 ###
\subsubsection{Drift} The remarkable drift condition (\ref{eq: b})
has been assumed in many previous papers. Here is the first time
where necessity and sufficiency of (\ref{eq: b}) are proved in the
full generality in the presence of jumps.
Note that in dimension $d=1$, the drift condition simply reduces to
nonnegativity $b\ge0$. But for dimension $d\ge2$, the boundary of
the state space $S_d^+$ becomes curved and kinked, implying a
nontrivial trade-off between diffusion $\alpha$ and $b$.

Concerning the form of $B$, let us
note the following: condition (\ref{eq: betaij}) implies in
particular
%
%e2.24 ###
%
\begin{equation}\label{eq: betaii}
\beta^{ii}_{\setminus\{i\}} -\int_{S_d^+\setminus\{0\}}
\frac{\chi(\xi)_{\setminus\{i\}}}{\|\xi\|^2 \wedge1}
\mu_{ii}(d\xi)\in S_{d-1}^+ \qquad\mbox{for all } 1\le i\le d,
\end{equation}
where for any matrix $u \in S_d$, $u_{\setminus\{i\}}$ denotes the
matrix where the $i$th row and column are deleted. Indeed,
inserting $x=c^{ii}$ in condition (\ref{eq: betaij}) yields
\[
\langle B(c^{ii}),u\rangle- \int_{S_d^+\setminus\{0\}}
\frac{\langle\chi(\xi),u\rangle}{\|\xi\|^2 \wedge1}
\mu_{ii}(d\xi)\geq0
\]
for all $u \in S_d^+$ with $\langle c^{ii},u\rangle=0$.
Since the $i$th column and row of such an element $u \in S_d^+$
is zero, it follows that
%
%e2.25 ###
%
\begin{eqnarray}\label{eq: Bchi}
&&\langle B(c^{ii}),u\rangle- \int_{S_d^+\setminus\{0\}}
\frac{\langle\chi(\xi),u\rangle}{\|\xi\|^2 \wedge1}
\mu_{ii}(d\xi)\nonumber\\[-8pt]\\[-8pt]
&&\qquad=\bigl\langle\beta^{ii}_{\setminus\{i\}},u_{\setminus\{i\}}
\bigr\rangle-\int_{S_d^+\setminus\{0\}}
\frac{\langle\chi(\xi)_{\setminus\{i\}},u_{\setminus\{i\}}\rangle
}{\|\xi\|^2
\wedge1} \mu_{ii}(d\xi)\geq0.\nonumber
\end{eqnarray}
By choosing appropriate elements $u_{\setminus\{i\}} \in S_{d-1}^+$,
we can further derive the integrability of $\chi(\xi)_{kl}$ for all
$k\neq i, l \neq i$, which implies
%
%e2.26 ###
%
\begin{equation}\label{eq: intchi}\qquad
\int_{S_d^+\setminus\{0\}}
\frac{\langle\chi(\xi)_{\setminus\{i\}},u_{\setminus\{i\}}\rangle
}{\|\xi\|^2
\wedge1} \mu_{ii}(d\xi)=\biggl\langle
\int_{S_d^+\setminus\{0\}}\frac{\chi(\xi)_{\setminus\{i\}}}{\|\xi
\|^2
\wedge1} \mu_{ii}(d\xi),u_{\setminus\{i\}}\biggr\rangle.
\end{equation}
As (\ref{eq: Bchi}) and (\ref{eq: intchi}) must hold true for all
$u_{\setminus\{i\}} \in S_{d-1}^+$, assertion (\ref{eq: betaii}) is
proved.

Note that the $(ij)$th component of the adjoint operator $B^\top$
is given by
%
%e2.27 ###
%
\begin{equation}\label{eq: B, B transpose}
B^\top_{ij}(u)=\langle\beta^{ij},u\rangle,
\end{equation}
since $\langle B(x),u\rangle=\langle
\sum_{i,j}\beta^{ij}x_{ij},u\rangle=\sum_{i,j}\langle
\beta^{ij},u\rangle x_{ij}=\langle B^\top(u),x\rangle$.

In most previous papers, $B(x)$ is of the form
%
%e2.28 ###
%
\begin{equation}\label{eq: wishart B}
B(x)= Hx+xH^{\top}.
\end{equation}
In this case,
%
%e2.29 ###
%
\begin{equation}\label{eq: wishart drift}
\langle B(x),u\rangle=\langle H x+x H^{\top},u\rangle=
0 \qquad\mbox{for all $x, u \in S_d^+$ with $\langle x, u\rangle
=0$,}\hspace*{-32pt}
\end{equation}
and hence (\ref{eq: betaij}) is equivalent to
\[
\int_{S_d^+\setminus\{0\}} \langle\chi(\xi),u\rangle M(x,d\xi)=0,
\]
for all $x, u \in S_d^+$ with $\langle x, u\rangle=0$.

If $B(x)$ is of the form
%
%e2.30 ###
%
\begin{equation}\label{eq: form B}
B(x)= Hx+xH^{\top}+\Gamma(x),
\end{equation}
where $H \in M_d$ and $\Gamma\dvtx S_d\to S_d$ linear satisfying
$\Gamma(S_d^+) \subseteq S_d^+$, then, in view of (\ref{eq: wishart
drift}), condition (\ref{eq: betaij}) holds true as long as
\[
\langle\Gamma(x), u\rangle-\int_{S_d^+\setminus\{0\}} \langle
\chi(\xi),u\rangle M(x,d\xi)\geq0
\]
for all $x, u \in S_d^+$ with $\langle x, u\rangle=0$. As a bold
conjecture, we claim that any $B(x)$ satisfying (\ref{eq: betaij})
is of form (\ref{eq: form B}).

Here is a simple example where $B(x)$ is of the form (\ref{eq: form
B}) but not of the usual form (\ref{eq: wishart B}): let $d=2$ and
\[
B(x)=\pmatrix{x_{22} & x_{12}\cr x_{12} &
x_{11}}.
\]
It can be easily checked that (\ref{eq: betaij}) is satisfied, while
$B(x)$ cannot be brought into the
form (\ref{eq: wishart B}). If $x_{ii}$ models the (squared)
volatility of the $i$th stock price, as in~(\ref{stockex}), then
this drift specification admits level impacts of the volatility of
stock~1 on the volatility of stock 2, and vice versa.

%s2.1.3 ###
\subsubsection{Killing} See Remark \ref{remcons}.

%s2.1.4 ###
\subsubsection{Jumps}\label{sec2.1.4}
 Condition (\ref{eq: m}) means that jumps
described by $m$, which can
for instance appear at $x=0$, should be of finite variation entering
the cone $S_d^+$, since infinite variation transversal to the
boundary would let the process leave the state space. Similarly,
condition (\ref{eq: mu}) asserts finite variation
for the inward
pointing directions, while we could a priori have a general jump
behavior (supported by
$S_d^+$ due to the affine structure) parallel to the boundary. Note that
in the case
$d=1$, which corresponds to $\mathbb{R}_+$, the linear jump part can
have infinite total variation
(see \cite{dfs}, equation (2.11)). However, due to the geometry of the
cone $S_d^+$, we conjecture that
in higher dimensions $d \geq 2$ such a behavior is no longer possible
and that all jumps are in fact of finite total variation.
In any case, for $d \geq 2$, affine positive matrix valued diffusion
processes cannot be approximated (in law) by pure jump processes, since
this
would yield a contradiction to condition (\ref{eq: b}). See also
Remark~\ref{remapproxinfdiv}
below.
\section{Affine processes are regular and Feller}\label{subsec:
regular and feller}

Suppose $X$ is an affine process on~$S_d^+$. The main result of this
section is that $X$ is regular in the sense of Definition~\ref{def:
affine process Sd+}. In addition, we shall prove that $P_t$ is a
Feller semigroup on $C_0(S_d^+)$. In order to show both properties,
we shall mainly rely on Lemma \ref{prop: feller prop} below. The
Feller property is then a simple consequence of this statement and
regularity is obtained by arguing as in Keller-Ressel, Schachermayer
and Teichmann \cite{kst}, who obtained the corresponding statements
for affine processes on the state space $\R^m_+\times\R^n$;
see~\cite{kst}, Theorem 4.3, and also the Ph.D. thesis of
Keller-Ressel \cite{keller}. We observe that most arguments
of \cite{kst} translate to our setting without major changes. It is
only required to tailor some technicalities to the cone $S_d^+$. We
start with the following elementary observations.
\begin{lemma}\label{lem: boundary}
If $u\in\partial S_d^+$ and $S_d^+\ni v \preceq u$, then $v\in
\partial S_d^+$.
\end{lemma}
\begin{pf}
Let $x\in S_d^+\setminus\{0\}$ such that $\langle x, u\rangle=0$.
Then, $S_d^+\ni v \preceq u$ implies $0 \leq\langle v,x\rangle
\leq\langle u,x\rangle=0$. Hence, $v \in
\partial S_d^+$.
\end{pf}

We now derive some first properties of the functions $\phi$ and
$\psi$ in (\ref{def: affine process}).
\begin{lemma}\label{lem: order preserving affine}
Let $X$ be an affine process on $S_d^+$. Then, we have:
\begin{longlist}[(iii)]
\item[(i)] The functions $\phi$ and $\psi$ satisfy
%
%e3.2 ###
%e3.1 ###
%
\begin{eqnarray}
\label{eq: flow phi}
\phi(t+s,u)&=&\phi(t,u)+\phi(s,\psi(t,u)),\\
\label{eq: flow psi}
\psi(t+s,u)&=&\psi(s,\psi(t,u))
\end{eqnarray}
for all $t,s \in\re_+$.
\item[(ii)] For all $u, v\in S_d^+$ with $v\preceq u$ and
for all $t \geq0$, the order relations
%
%e3.3 ###
%
\begin{equation}\label{eq: order pres}
\phi(t,v)\leq\phi(t,u) \quad\mbox{and}\quad \psi(t,v) \preceq
\psi(t,u)
\end{equation}
hold true.
\item[(iii)] The functions $\phi$ and $\psi$ are
jointly continuous in $\re_+\times S_d^+$. Furthermore, $u \mapsto
\phi(t,u)$ and $u\mapsto\psi(t,u)$
are analytic on $S_d^{++}$.
\end{longlist}
\end{lemma}
\begin{pf}
Assertion (i) follows directly from the
Chapman--Kolmogorov equation,
\begin{eqnarray*}
e^{-\phi(t+s,u)-\langle\psi(t+s,u),x\rangle}&=&\int_{S_d^+}
p_s(x,d\xi)\int_{S_d^+}e^{-\langle u,\widetilde{\xi}\rangle
}p_{t}(\xi,d\widetilde{\xi})\\
&=&e^{-\phi(t,u)}\int_{S_d^+}e^{-\langle\psi(t,u),\xi\rangle
}p_s(x,d\xi)\\
&=&e^{-\phi(t,u)-\phi(s,\psi(t,u))-\langle\psi(s,\psi
(t,u)),x\rangle}.
\end{eqnarray*}
For the proof of (ii), note that $v \preceq u$ is
equivalent to $\langle v,x\rangle\leq\langle u, x\rangle$ for all
$x \in S_d^+$. By the monotonicity of the exponential function, we
have for all $x\in S_d^+$ and for all $t\geq0$,
\begin{eqnarray*}
e^{-\phi(t,v)-\langle\psi(t,v),x\rangle}&=&\int_{S_d^+}e^{-\langle
v,\xi\rangle}p_t(x,{d}\xi) \geq\int_{S_d^+}e^{-\langle
u,\xi\rangle}p_t(x,{d}\xi)\\
&=&e^{-\phi(t,u)-\langle
\psi(t,u),x\rangle},
\end{eqnarray*}
and the assertion follows by taking logarithms.

Concerning statement (iii), note that stochastic
continuity of $X$ implies joint continuity of $P_t e^{-\langle
u,x\rangle}$ in $(t,u)\in\R_+\times S_d^+$ (this follows, e.g., from
\cite{bau96}, Lem\-ma~23.7), for all $x\in S_d^+$. This in turn
yields continuity of the functions $(t,u) \mapsto\phi(t,u)$ and
$(t,u) \mapsto\psi(t,u)$. The second assertion follows from
analyticity properties of the Laplace transform.
\end{pf}

The following property of $\psi$ is crucial.
\begin{lemma}\label{prop: feller prop}
Let $\psi\dvtx\mathbb R_+\times S_d^+\rightarrow S_d^+$ be any map
satisfying $\psi(0,u)=u$ and the properties \textup{(i)--(iii)} of
Lemma \ref{lem: order preserving affine} (regarding the function $\psi
$). Then $\psi(t,u) \in
S_d^{++}$ for all $(t,u) \in\R_+ \times S_d^{++}$.
\end{lemma}
\begin{pf}
We adapt the proof of \cite{keller}, Proposition 1.10, to our
setting. Assume by contradiction that there exists some $(t,u) \in
\re_+\times S_d^{++}$ such that $\psi(t,u) \in
\partial S_d^+$. Let us consider the interval $(0, \lambda
_{\min}(u)]\neq\varnothing$, where $\lambda_{\min}(u)>0$
denotes the smallest eigenvalue of $u$. Then for all $v \in(0,
\lambda_{\min}(u)]$ we have $vI_d \preceq u$. Since $\psi(t,u)$
admits property (ii) of Lemma \ref{lem: order preserving affine},
we obtain
\[
S_d^+ \ni\psi(t,vI_d) \preceq\psi(t,u) \in\partial S_d^+
\]
for all $v \in(0, \lambda_{\min}(u)]$. Consequently,
Lemma \ref{lem: boundary} yields that $\psi(t,vI_d)\in\partial
S_d^+$. Hence,
\[
\det(\psi(t,vI_d))=0
\]
for all $v \in(0, \lambda_{\min}(u)]$. The analyticity of $u \mapsto
\psi(t,u)$ on $S_d^{++}$ carries over to $u \mapsto\det(\psi(t,u))$
and implies that $\det(\psi(t,vI_d))=0$ for all $v \in\re_{++}$.
Indeed, the set of zeros of $\det(\psi(t,vI_d))$ has an accumulation
point in $\re_{++}$, which implies that $\det(\psi(t,vI_d))$
vanishes entirely on $\re_{++}$. The same statement holds true for
$t$ replaced by $\frac{t}{2}$. Indeed, if $\psi(\frac{t}{2},u)\in
\partial S_d^+$, then the assertion is shown by the same arguments
as above. Otherwise, if $\psi(\frac{t}{2},u)\in S_d^{++}$, we have
for all $v \in\re_{++}$ with $vI_d \preceq\psi(\frac{t}{2},u)$,
that is, for all $v \in(0, \lambda_{\min}(\psi(\frac{t}{2},u)]$
\[
S_d^+ \ni\psi\biggl(\frac{t}{2},vI_d\biggr) \preceq
\psi\biggl(\frac{t}{2},\psi\biggl(\frac{t}{2},u\biggr)\biggr)=\psi(t,u)
\in\partial S_d^+,
\]
which yields again $\psi(\frac{t}{2},vI_d)\in\partial S_d^+$ and
$\det(\psi(\frac{t}{2},vI_d))=0$ for all
$v \in(0$, $\lambda_{\min}(\psi(\frac{t}{2},u)]$. The same reasoning as
before then leads to $\det(\psi(\frac{t}{2},vI_d))=0$ for all $v \in
\re_{++}$. By reapplying this argument, we finally get for every $n
\in\mathbb{N}$ and for all $v \in\re_{++}$
\[
\det\biggl(\psi\biggl(\frac{t}{2^n},vI_d\biggr)\biggr)=0.
\]
From the continuity of the function $t \mapsto\psi(t,u)$ and of the
determinant, we deduce that for any $v \in\re_{++}$,
\[
0=\lim_{n\rightarrow
\infty}\det\biggl(\psi\biggl(\frac{t}{2^n},vI_d\biggr)\biggr)=\det(\psi(0,vI_d))=\det
(vI_d)=v^d>0,
\]
a contradiction, and the assertion is proved.
\end{pf}

We may now formulate the main result of this section.
\begin{proposition}\label{th: regularity}
Let $X$ be an affine process with state space $S_d^+$. Then, we
have:
\begin{longlist}[(ii)]
\item[(i)] $X$ is a Feller process.
\item[(ii)] $X$ is regular.
\end{longlist}
\end{proposition}
\begin{pf}
In order to prove (i), it suffices to show that for
all $f \in C_0(S_d^+)$
%
%e3.5 ###
%e3.4 ###
%
\begin{eqnarray}
\label{eq: fellerprop1}
\lim_{t\rightarrow0+}P_tf(x) &=&f(x)\qquad\mbox{for all $x \in
S_d^+$,}\\
\label{eq: fellerprop2}
P_tf &\in& C_0(S_d^+)\qquad\mbox{for all $t \in\re_+$},
\end{eqnarray}
(see, e.g., \cite{revuzyor}, Propostion III.2.4).
Property (\ref{eq: fellerprop1}) is a consequence of stochastic
continuity, which implies for all $f \in C_0(S_d^+)$ and $x\in S_d^+$
\[
\lim_{t\rightarrow0^+}P_tf(x)= f(x).
\]
Concerning (\ref{eq: fellerprop2}), it suffices to verify this property
for a dense subset of $C_0(S_d^+)$. By a locally compact version of
Stone--Weierstrass' theorem (see, e.g., \cite{semadeni}), the linear
span of the set $\{e^{-\langle u,x\rangle} \mid u\in S_d^{++}\}$ is dense
in $C_0(S_d^+)$. Indeed, it is a subalgebra of $C_0(S_d^+)$, separates
points and vanishes nowhere, as all elements are strictly positive
functions on $S_d^+$. From Lemma \ref{prop: feller prop}, we can deduce
that $P_te^{-\langle u,x\rangle} \in C_0(S_d^+)$ if $u \in S_d^{++}$,
since $\psi(t,u) \in S_d^{++}$ and $\langle\psi (t,u),x\rangle> 0$ for
$x \neq0$ implying that
\[
P_te^{-\langle u,x\rangle}=e^{-\phi(t,u)-\langle\psi(t,u),x\rangle}
\]
goes to $0$ as $x \rightarrow\Delta$. Hence, statement (i) is proved.

The proof of (ii) follows precisely the lines
of \cite{kst}, proof of Theorem 4.3. Using Lemma \ref{prop:
feller prop}, one may mimic the proof of \cite{kst}, Theorem 4.3, to
obtain that differentiability of $\psi(t,u)$ in $u\in S_d^{++}$,
which follows from Lemma \ref{lem: order preserving
affine}(iii), implies differentiability of
$\psi(t,u)$ in $t$ for $t=0$ and for all $u \in S_d^{+}$.
\end{pf}

By the regularity of $X$, we are now allowed to differentiate the
equations (\ref{eq: flow phi}) and (\ref{eq: flow psi}) with respect
to $t$ and evaluate them at $t=0$. As a consequence, $\phi$ and
$\psi$ satisfy the system of differential equations
\begin{eqnarray*}
\frac{\partial\phi(t,u)}{\partial t}&=&F(\psi(t,u)),\qquad \phi(0,u)=0, \\
\frac{\partial\psi(t,u)}{\partial t}&=&R(\psi(t,u)),\qquad
\psi(0,u)=u \in S_d^+,
\end{eqnarray*}
where $F$ and $R$ are defined as in (\ref{eq: F,R}). The analysis of
these (generalized Riccati) differential equations is subject of
Section \ref{section: Riccati}, whereas the specific form of $F$ and
$R$ is elaborated in the following.

%s4 ###
\section{Necessary parameter restrictions}\label{section: Nec
paramter restrictions}

In this section, we derive \textit{necessary} parametric restrictions,
that is, given an affine process on $S_d^+$, we determine necessary
implications on a set of parameters which only ensue from
Definition \ref{def: affine process Sd+}. These conditions are
precisely the conditions on the admissible parameter set as of
Definition \ref{def: necessary admissibility}. The form of the
functions $F$ and $R$ as defined by (\ref{eq: F,R}) is then
characterized by means of this parameter set, which is stated in
Proposition \ref{th: necessary admissibility} below. For its proof,
we first provide a number of technical prerequisites.
\begin{lemma}\label{lem: zero divisor}
Let $x, u \in S_d^+$ and
%
%e4.1 ###
%
\begin{equation}\label{eq: form x}
x=O\Lambda O^{\top}=O\diag(\lambda_1 >0,\ldots,
\lambda_{d-r}>0,0,\ldots, 0)O^\top
\end{equation}
be the diagonalization of $x$ with $r\geq0$ and $O \in O(d)$. Then
the following assertions are equivalent:
\begin{longlist}[(iii)]
\item[(i)] $ux=xu=0$,
\item[(ii)] $\langle x,u\rangle=0$,
\item[(iii)] $u$ is of form
%
%e4.2 ###
%
\begin{equation}\label{eq: form zerodivisor}
u=O\pmatrix{
0 & 0\cr
0 & w}
O^\top
\end{equation}
with $w \in S^+_{r}$.
\end{longlist}
\end{lemma}
\begin{pf}
The direction (i)${}\Rightarrow{}$(ii) is obvious. In
order to prove the implication (ii)${}\Rightarrow{}$(iii),
define $v$ as $v=O^{\top}uO$. Then we have
\[
0=\langle x, u\rangle=\langle\Lambda, O^{\top}u
O\rangle=\sum_{i\leq d-r}\lambda_i v_{ii},
\]
which implies $v_{ii}=0$ for all $i \leq d-r$ and by the positive
definiteness of $v$ it must then be of form
\[
v=\pmatrix{
0 & 0\cr
0 & w}
\]
with $w \in S^+_{r}$. Thus, $u$ is given by (\ref{eq: form
zerodivisor}). This then implies that $ux=xu=0$, which proves the
direction (iii)${}\Rightarrow{}$(i).
\end{pf}
\begin{lemma}\label{lemprojector}
Let $p$ be an orthogonal projector, that is, $p\in S_d^+$ and
$p^2=p$ (see, e.g., Kato \cite{kat95}, Section \textup{I.6.7}), and define
$q=I_d-p$. Then $q$ is an orthogonal projector and the orthogonal
complement of $p$ in $S_d^+$ equals
\[
\{ v\in S_d^+\mid\langle p,v\rangle=0\}=\{quq\mid u\in
S_d^+\}.
\]
\end{lemma}
\begin{pf}
That $q$ is an orthogonal projector\vspace*{1pt} follows by inspection. The
diagonalization of $p$ is of the form $p=O\Lambda O^{\top}$ with
$\Lambda=\diag(1,\ldots,1,0,\ldots,0)$, and thus $q=O(I_d-\Lambda)
O^{\top}$. In view of Lemma \ref{lem: zero divisor}, we conclude that
$v\in S_d^+$ is orthogonal to $p$ if and only if $v=qvq$. This
proves the assertion.
\end{pf}
\begin{lemma}\label{lem: linear map}
Let $u$ be in $S_d$ and $x \in\partial S_d^+$ such that $ux=xu=0$.
Then, the linear map $T_u$ defined by
\[
T_u\dvtx S_d \rightarrow S_d,\qquad v \mapsto T_uv:=uvu
\]
has the following properties:
\begin{longlist}[(ii)]
\item[(i)] $T_u$ is self-adjoint and
$T_u(S_d^+-\re_+x)\subseteq S_d^+$.
\item[(ii)] There exists an element $v \in S_d$
such that $T_uv=u$.
\end{longlist}
\end{lemma}
\begin{pf}
The assertion (i) is obvious, since for every
$k \in\re_+$, $T_ukx=kuxu=0$ and $T_uv=uvu \in S_d^+$ if $v \in
S_d^+$. For proving part (ii), we use the fact
that $x$ is of form (\ref{eq: form x}) and that all zero divisors
$u$ in $S_d$ of $x$ can be represented by (\ref{eq: form zerodivisor})
with $w \in S_{r}$. Thus, setting
\[
v=O\pmatrix{
0 & 0\cr
0 & w^{+}}
O^\top,
\]
where $w^{+}$ satisfies $ww^{+}w=w$, yields $T_uv=u$.
\end{pf}
\begin{lemma}\label{lem: linear map cone}
Let $V$ denote a vector space\footnote{In the proof of
Proposition \ref{th: necessary admissibility} below, $V$ corresponds
to $S_d$, the vector space of linear maps $S_d\to S_d$, or the
vector space of finite signed measures on $S_d$.} over $\re$. Let
$L\dvtx S_d^+ \rightarrow V$ be an additive (resp., homogeneous additive)
map, that is, for all $x,y \in S_d^+$ and $\lambda=1$ (resp., for all
$\lambda\in\mathbb R_+$) we have
\[
L(x +\lambda y)=L(x)+\lambda L(y).
\]
Then $L(x)$ is the restriction of an additive (resp., $\mathbb
R$-linear) map on $S_d$.
\end{lemma}
\begin{pf}
We define the map $\widetilde L\dvtx S_d\rightarrow V$ as
\[
\widetilde L(x-y):=L(x)-L(y),\qquad x,y\in S_d^+.
\]
$\widetilde L$ is well defined, as for $u,v,x,y\in S_d^+$ such that
$u-v=x-y$ we have
\[
L(u)-L(v)=\widetilde L(u-v)=\widetilde L(x-y)=L(x)-L(y).
\]
Since $S_d^+-S_d^+=S_d$, the domain of $\widetilde L$ is all of
$S_d$. Also, $L(0)=0$ by the additivity of $L$. Hence, $ L$
is the restriction of $\widetilde L$ to $S_d^+$. Homogeneity of
$\widetilde L$
holds, as for $\lambda>0, z=x-y\in S_d$ we have by definition
\[
\widetilde L(\lambda z)=L(\lambda x)-L(\lambda y)=\lambda
L(x)-\lambda L(y)=\lambda{\widetilde L}(z).
\]
Finally, we show additivity of $\widetilde L$. Choose $w,z\in S_d$
such that $z=x-y, w=u-v$, hence $w+z=(x+u)-(y+v)$. By the
definition of $\widetilde{L}$, we have
\[
\widetilde L(z)=L(x)-L(y),\qquad \widetilde L(w)=L(u)-L(v),
\]
and by the additivity of $L$ we obtain
\begin{eqnarray*}
\widetilde L(w+z)&=&L(x+u)-L(y+v)=L(x)+L(u)-L(y)-L(v)\\
&=&\widetilde
L(z)+\widetilde L(w).
\end{eqnarray*}
\upqed\end{pf}

We now provide a convergence result for Laplace transforms (in fact
Laplace--Fourier transforms), which is most relevant for the analysis
of affine processes.
\begin{lemma}\label{lem: FourierLaplace}
Let $\nu_n$ be a sequence of measures on $S_d$ with
\[
L_n(u)=\int_{S_d}
e^{-\langle
u,\xi\rangle} \nu_n(d\xi)<\infty\quad\mbox{and}\quad \lim_{n \rightarrow
\infty} L_n(u) =
L(u),\qquad u\in S_d^+,
\]
pointwise, for some finite function $L$ on
$S_d^+$ continuous at $ u = 0 $. Then $\nu_n$ converges weakly to
some finite measure $\nu$ on $S_d$ and the Fourier--Laplace
transforms converge for $ u \in S_d^{++} \cup\{0\} $ and $ v \in
S_d $ to the Fourier--Laplace transforms of $ \nu$, that is,
\[
\lim_{n \rightarrow\infty} \int_{S_d} e^{-\langle u + \im v,\xi
\rangle} \nu_n(d\xi) =
\int_{S_d} e^{-\langle u + \im v,\xi\rangle} \nu(d\xi).
\]
In particular, $\nu(S_d)=\lim_{n \rightarrow\infty}\nu_n(S_d)$ and
\[
L(u) = \int_{S_d} e^{-\langle u ,\xi\rangle} \nu(d\xi),
\]
for all $ u \in S_d^{++} \cup\{ 0 \} $.
\end{lemma}
\begin{remark}
Instead of $ u = 0 $ we could take any set $K$ of points at the
boundary $ K \subset\partial S_d^{+} $: if we assume continuity of
$L$ at points in $K$, then we obtain the equality of $ L $ with the
Laplace transform of $ \nu$ for all points in $K$. Additionally,
continuity is too strong an assumption, since we only need right
continuity of $L$ along the segment $ u + \varepsilon
I_d $ for $ \varepsilon= 0 $ at the points from the
boundary under consideration.
\end{remark}
\begin{pf*}{Proof of Lemma \ref{lem: FourierLaplace}}
Since $\nu_n(S_d)=L_n(0)$ is bounded, we know by
general theory that $ \nu_n $ has a vague accumulation point $\nu$,
which is a finite measure on~$S_d$.

Since $L_n(u)<\infty$ on $S_d^+$, it follows by well-known
regularity properties of Laplace transforms (see, e.g.,
\cite{fil09}, Lemma 10.8) that the functions $L_n$ admit an
analytic extension on the strip $ S_d^{++} + \im S_d $, still denoted
by $L_n$:
\[
(u+ \im v ) \mapsto L_n(u+\im v)=\int_{S_d} e^{-\langle u + \im
v,\xi\rangle} \nu_n(d\xi).
\]
Moreover, pointwise convergence of the
finite convex functions $L_n$ to $L$ on $S_d^+$ implies that this
convergence is in fact uniform on compact subsets of $ S_d^{++} $
(see, e.g., Rockafellar \cite{roc97}, Theorem 10.8). Hence, the
functions $L_n$ are uniformly bounded on compact subsets of $S_d^{++}$ and
since $|L_n(u+ \im v)|\le L_n(u)$, also on compact subsets of $S_d^{++}
+\im S_d$.
Therefore, and since $ S_d^{++} $ is a set of uniqueness in $ S_d^{++}
+ \im S_d $,
it follows by Vitali's theorem (\cite{Narasimhan71}, Chapter 1, Proposition
7) that the analytic functions $L_n$ converge
uniformly on compact subsets of $ S_d^{++} + \im S_d $ to an analytic
limit thereon. By L\'evy's continuity theorem, we therefore know
that for any $ u \in S_d^{++} $ the
finite measures $ \exp(- \langle u, \xi\rangle) \nu_n (d\xi) $
converge weakly to a limit, which by uniqueness of the weak limit
has to equal $ \exp(- \langle u, \xi\rangle) \nu(d\xi) $. Whence
the only vague accumulation point of $ \nu_n $ is $ \nu$. Vague
convergence implies weak convergence if mass is conserved.
Continuity of $ L(u) $ at $ u = 0 $ implies this mass
conservation: indeed, by weak convergence of $ e^{-\langle\varepsilon
I_d,\xi\rangle} \nu_n $ we arrive at
\begin{eqnarray*}
L(\varepsilon I_d) &=& \lim_{n\rightarrow\infty} \int_{S_d}
e^{-\langle\varepsilon I_d,\xi\rangle} \nu_n (d\xi) = \int_{S_d}
e^{-\langle\varepsilon I_d,\xi\rangle} \nu(d\xi) \\
&=& \int_{S_d} e^{-\langle\varepsilon I_d, \xi\rangle} 1_{\{\langle
I_d, \xi\rangle\leq0\}} \nu(d\xi) + \int_{S_d} e^{-\langle
\varepsilon I_d, \xi\rangle} 1_{\{\langle I_d ,\xi\rangle> 0\}
} \nu(d\xi)
\end{eqnarray*}
and therefore---by dominated convergence---we obtain that the limit
$ \varepsilon\to0 $ yields
\[
L(0) = \int_{S_d} \nu(d\xi),
\]
which is the desired mass conservation, hence weak convergence,
which means in turn convergence of the Fourier--Laplace transform at
$ u = 0 $.
\end{pf*}

Finally, let us state a general comparison result for ODEs and
hereto introduce the notion of \textit{quasi-monotonicity}, which we
shall need several times throughout this article, in particular in
the proofs of Propositions \ref{th: necessary admissibility} and
\ref{prop_ricc_sol} below.
\begin{definition}\label{def: quasimono}
Let $U\subset S_d$ be an open set. A function $f\dvtx U \rightarrow S_d$
is called \textit{quasi-monotone increasing} if for all elements $x,
y \in U$, $u\in S_d^+$ which satisfy $x \preceq y$ and $\langle
x,u\rangle=\langle y,u\rangle$,
\[
\langle f(x),u\rangle\leq\langle f(y),u\rangle
\]
holds true. Accordingly, we call $f$ \textit{quasi-constant} if both
$f$ and $-f$ are quasi-monotone increasing.
\end{definition}

The following comparison result can be deduced from a more general
theorem proved by Volkmann \cite{Volkmann1973}.
\begin{theorem}\label{th: Volkmann}
Let $U\subset S_d$ be an open set. Let $f\dvtx[0,T)\times U\to
S_d$ be a continuous locally Lipschitz map such that $f(t,\cdot)$ is
quasi-monotone increasing on $U$ for all $t\in[0,T)$. Let
$0<t_0\leq T$ and $x,y\dvtx[0,t_0)\to U$ be differentiable maps
such that $x(0)\preceq y(0)$ and
\[
\dot x(t)-f(t,x(t))\preceq\dot y(t)-f(t,y(t)),\qquad 0\leq
t<t_0.
\]
Then we have $x(t)\preceq y(t)$ for all $t\in[0,t_0)$.
\end{theorem}

%s4.1 ###
\subsection{The functions F and R}

The main result of this section characterizes the form of the
functions $F$ and $R$ as defined by (\ref{eq: F,R}).
\begin{proposition}\label{th: necessary admissibility}
Let $X$ be an affine process with state space $S_d^+$. Then there
exist parameters $(\alpha, b, \beta^{ij}, c, \gamma, m, \mu)$, where
$\alpha, \beta^{ij}, c, \gamma, m, \mu$ satisfy the admissibility
conditions of Definition \ref{def: necessary admissibility} and $b
\in S_d^+$, such that the functions $F$ and $R$ are of the
form (\ref{eq: F}) and (\ref{eq: R}).
\end{proposition}
\begin{remark}
Note that for the moment we only obtain $b \in S_d^+$, and
not (\ref{eq: b}).
\end{remark}
\begin{pf*}{Proof of Proposition \ref{th: necessary admissibility}}
As the proof of Proposition \ref{th: necessary admissibility} is
rather long, we divide it into several steps:

\textit{Step} 1. \textit{Necessary admissibility conditions for $b, c,
\gamma, m$.} In order to derive the particular form of $F$ and $R$
with the above parameter restrictions, we follow the approach of
Keller-Ressel \cite{keller}, Theorem 2.6. Note that the
$t$-derivative of $P_te^{-\langle u,x\rangle}$ at $t=0$ exists for
all $x,u \in S_d^+$, since
%
%e4.3 ###
%
\begin{eqnarray}\label{Asharpeux}
\lim_{t \rightarrow0^+}\frac{P_te^{-\langle u,x\rangle}-e^{-\langle
u,x\rangle}}{t}&=&\lim_{t \rightarrow0^+}\frac{e^{-\phi
(t,u)-\langle\psi(t,u),x\rangle}-e^{-\langle
u,x\rangle}}{t}\nonumber\\[-8pt]\\[-8pt]
&=&\bigl(-F(u)-\langle R(u),x\rangle\bigr)e^{-\langle u, x\rangle}\nonumber
\end{eqnarray}
is well defined by Proposition \ref{th: regularity}. Moreover, we
can also write
\begin{eqnarray*}
&&
-F(u)-\langle R(u),x\rangle\\
&&\qquad=\lim_{t \rightarrow
0^+}\frac{P_te^{-\langle u,x\rangle}-e^{-\langle u,x\rangle}}{t
e^{-\langle u,x\rangle}} \\
&&\qquad=\lim_{t\rightarrow
0^+}\frac{1}{t}\biggl(\int_{S_d^+\setminus\{0\}}e^{-\langle u,
\xi-x\rangle}p_t(x,d\xi)-1\biggr)\\
&&\qquad=\lim_{t\rightarrow
0^+}\biggl(\frac{1}{t}\int_{S_d^+-x}\bigl(e^{-\langle
u,\xi\rangle}-1\bigr)p_t(x,d\xi+x) +\frac{p_t(x,S_d^+)-1}{t}\biggr).
\end{eqnarray*}
By the above equalities and the fact that $p_t(x,S_d^+)\leq1$, we
then obtain for $u=0$
\[
0 \geq\lim_{t\rightarrow0^+}\frac{p_t(x,S_d^+)-1}{t}=-F(0)-\langle
R(0),x\rangle.
\]
Setting $F(0)=c$ and $R(0)=\gamma$ yields $c \in\re^+$ as
in (\ref{eq: c}) and $\gamma\in S_d^+$ as in (\ref{eq: gamma}). We
thus obtain
%
%e4.4 ###
%
\begin{eqnarray}\label{eq: FRinfdiv}
&&-\bigl(F(u)-c\bigr) - \langle R(u)-\gamma,x\rangle\nonumber\\[-8pt]\\[-8pt]
&&\qquad=\lim_{t\rightarrow0^+}
\frac{1}{t}\int_{S_d^+-x}\bigl(e^{-\langle
u,\xi\rangle}-1\bigr)p_t(x,d\xi+x).\nonumber
\end{eqnarray}
For every fixed $t>0$, the right-hand side of (\ref{eq: FRinfdiv})
is the logarithm of the Laplace transform of a compound Poisson
distribution supported on $S_d^+-\re_+x$ with \mbox{intensity}
$p_t(x,S_d^+)/t$ and compounding distribution
$p_t(x,d\xi+x)/p_t(x,S_d^+)$. Concerning the support, note that the
compounding distribution is\vadjust{\goodbreak} concentrated on $S_d^+-x$, which implies
that the compound Poisson distribution has support on the convex
cone $S_d^+-\re_+ x$. By Lemma \ref{lem: FourierLaplace}, the
pointwise convergence of (\ref{eq: FRinfdiv}) for $t\rightarrow0$
to some function being continuous at $0$, implies weak convergence
of the compound Poisson distributions to some infinitely divisible
probability distribution $K(x,dy)$ supported on $S_d^+-\re_+x$.
Indeed, this follows from the fact that any compound Poisson
distribution is infinitely divisible and the class of infinitely
divisible distributions is closed under weak convergence
(\cite{sato}, Lemma 7.8). Again, by Lemma \ref{lem: FourierLaplace}
the Laplace transform of $K(x,dy)$ is then given as exponential of
the left-hand side of (\ref{eq: FRinfdiv}).

In particular, for $x=0$, $K(0,dy)$ is an infinitely divisible
distribution with support on the cone $S_d^+$. By the
L\'evy--Khintchine formula on proper cones (see \cite{skorohod},
Theorem 3.21), its Laplace transform is therefore of the form
\[
\exp\biggl(-\langle b,u\rangle+\int_{S_d^+\setminus\{0\}}\bigl(e^{-\langle
u, \xi\rangle}-1\bigr)m(d\xi)\biggr),
\]
where $b \in S_d^+$ and $m$ is a Borel measure supported on $S_d^+$
such that
\[
\int_{S_d^+\setminus\{0\}}(\|\xi\|\wedge1)m(d\xi) <
\infty,
\]
yielding (\ref{eq: m}). Therefore,
\[
F(u)=\langle b, u\rangle+ c -
\int_{S_d^+\setminus\{0\}}\bigl(e^{-\langle u, \xi\rangle}-1\bigr)m(d\xi).
\]

\textit{Step} 2. \textit{Necessary admissibility conditions for $\beta^{ij},
\mu$.} We next obtain the particular form of $R$. Observe that for
each $x \in S_d^+$ and $k\in\N$,
\[
\exp\bigl( -\bigl(F(u)-c\bigr)/k - \langle R(u)-\gamma,x\rangle\bigr)
\]
is the Laplace transform of the infinitely divisible distribution
$K(kx,dy)^{\ast{1/k}}$, where $\ast\frac{1}{k}$ denotes the
$\frac{1}{k}$ convolution power. For $k\to\infty$, these Laplace
transforms obviously converge to $\exp(-\langle
R(u)-\gamma,x\rangle)$ pointwise in $u$. Using again the same
arguments as before [an application of Lemma \ref{lem:
FourierLaplace} as below (\ref{eq: FRinfdiv})], we can deduce that
$K(kx,dy)^{\ast{1/k}}$ converges weakly to some infinitely
divisible distribution $ L(x,dy)$ on $S_d^+-\re_+x$ with Laplace
transform $\exp(-\langle R(u)-\gamma,x\rangle)$ for $u\in S_d^+$.

By the L\'evy--Khintchine formula on $S_d$ (\cite{sato},
Theorem 8.1,
indeed on $\re^{(d(d+1)/2)}$ by modifying the scalar product
appropriately), the characteristic function of
$L(x,dy)$ has the form
\begin{eqnarray*}
\widehat{L}(x,u)&=&\exp\biggl(\frac{1}{2}\langle u, A(x) u \rangle+\langle
B(x),u\rangle\\
&&\hspace*{21.2pt}{} +\int_{S_d\setminus\{0\}} \bigl(e^{-\langle u,
\xi\rangle}-1-\langle\chi(\xi), u
\rangle\bigr)M(x, d \xi)\biggr),
\end{eqnarray*}
for $u\in\im S_d$, where $A(x)$ is a symmetric positive semidefinite
linear operator on
$S_d$, $B(x) \in S_d$, $M(x,\cdot)$ a measure on $S_d\setminus\{0\}$
satisfying
\[
\int_{S_d\setminus\{0\}}(\|\xi\|^2 \wedge1) M(x,d\xi) < \infty,
\]
and $\chi$ some appropriate truncation function. Furthermore,
by \cite{sato}, Theorem 8.7,
%
%e4.5 ###
%
\begin{eqnarray}\label{eq: Levymeasure limit}
&&\int_{S_d\setminus\{0\}}f(\xi)\frac{1}{t}p_t(x,d\xi+x)\nonumber\\[-9pt]\\[-9pt]
&&\qquad\stackrel{t
\rightarrow0}{\longrightarrow}
\int_{S_d\setminus\{0\}}f(\xi)m(d\xi)+\int_{S_d\setminus\{0\}
}f(\xi)M(x,d\xi)\nonumber
\end{eqnarray}
holds true for all $f\dvtx S_d \rightarrow\re$ which are bounded,
continuous and vanishing on a neighborhood of $0$. We conclude that
$M(x,d\xi)$ has support in $S_d^+-x$. Therefore, the characteristic
function $\widehat{L}(x,u)$ admits an analytic extension to
$S_d^+\times\im S_d$, which then has to coincide with the Laplace
transform for $u\in S_d^+$. We conclude that, for all $x\in S_d^+$,
%
%e4.6 ###
%
\begin{eqnarray}\label{eq: levy khintchine R}
&&-\langle R(u)-\gamma,x\rangle\nonumber\\[-1pt]
&&\qquad=\frac{1}{2}\langle u, A(x)
u \rangle-\langle B(x),u\rangle\\[-1pt]
&&\qquad\quad{}+\int_{S_d\setminus\{0\}}
\bigl(e^{-\langle u, \xi\rangle}-1+\langle
\chi(\xi), u \rangle\bigr)M(x, d
\xi),\qquad u\in S_d^+.\nonumber
\end{eqnarray}

As the left-hand side of (\ref{eq: levy khintchine R}) is linear in the
components of $x$, it follows that $x \mapsto A(x)$, $x\mapsto B(x)$
as well as $x \mapsto\int_E(\|\xi\|^2 \wedge1) M(x,d\xi)$ for
every $E \in\mathcal{B}(S_d\setminus\{0\})$ are homogeneous
additive maps on $S_d^+$ in the sense of Lemma \ref{lem: linear map
cone}. This then implies that they are restrictions of linear maps
on $S_d$, such that we can write
\begin{eqnarray*}
A(x)&=&\sum_{i,j}a_{ij}x_{ij},\qquad B(x)=\sum_{i,j}\beta^{ij}x_{ij},\\[-1pt]
\int_E(\|\xi\|^2 \wedge1) M(x,d\xi)&=&\langle x,
\mu(E)\rangle=\sum_{i,j}\mu_{ij}(E)x_{ij},
\end{eqnarray*}
where (recall that $c^{ij}$ denotes the standard basis of $S_d$
defined in Section \ref{subsec: notation}):
\begin{eqnarray*}
a_{ij}&=&a_{ji}=(1+\delta_{ij})\frac{A(c^{ij})}{2}\dvtx S_d\to
S_d \qquad\mbox{linear},\\[-1pt]
\beta^{ij}&=&\beta^{ji}=(1+\delta_{ij})\frac{B(c^{ij})}{2}\in S_d
\end{eqnarray*}
and\vspace*{-2pt}
\[
E\mapsto\mu_{ij}(E)=\mu_{ji}(E)=(1+\delta_{ij})\frac{\int_E(\|\xi
\|^2
\wedge1)M(c^{ij},d\xi)}{2}
\]
are finite signed measures on $S_d\setminus\{0\}$. The fact that
$M(x,\cdot)$ is a nonnegative measure for each $x\in S_d^+$ implies
immediately that $\mu(E)$ is a positive semidefinite matrix.\vadjust{\goodbreak}

In (\ref{eq: Levymeasure limit}), take now $x=\frac{1}{n}e^{ij}$ and
nonnegative functions $f=f_n \in C_b(S_d)$ with
$f_n=0$ on $S_d^+-\frac{1}{n}e^{ij}$.
Then for each $n$ the left-hand side of (\ref{eq: Levymeasure limit})
is zero since the
$p_t(\frac{1}{n}e^{ij},d\xi+\frac{1}{n}e^{ij})$ is concentrated on
$S_d^+-\frac{1}{n}e^{ij}$. As $\supp(m)\subseteq S_d^+$, the first
integral on the right vanishes as well. Hence,
\begin{eqnarray*}
0&=&\int_{S_d\setminus\{0\}}f_n(\xi)M\biggl(\frac{1}{n}e^{ij},d\xi\biggr)=\int
_{S_d\setminus\{0\}}
\frac{f_n(\xi)}{\|\xi\|^2 \wedge1}\biggl\langle\frac{1}{n}e^{ij}, \mu
(d\xi)\biggr\rangle\\
&=&\frac{1}{n}\int_{S_d\setminus\{0\}}\frac{f_n(\xi)}{\|\xi\|^2
\wedge
1}\bigl(\mu_{ii}(d\xi)+(1-\delta_{ij})\bigl(\mu_{jj}(d\xi)+2\mu_{ij}(d\xi)\bigr)\bigr)
\end{eqnarray*}
for any nonnegative function $f_n \in C_b(S_d)$ with $f_n=0$ on
$S_d^+-\frac{1}{n}e^{ij}$ implies that $\supp(\mu_{ij})\subseteq
S_d^+-\frac{1}{n}e^{ij}$ for
each $n$. Thus, we can conclude that $\supp\mu_{ij} \subseteq
S_d^+$ for all $1\le i,j\le d$.

Now let $T\dvtx S_d \rightarrow S_d$ be any linear map with the property
$T(S_d^+-\re_+x)\subseteq S_d^+$. Then
$T(\supp(L(x,dy)))\subseteq S_d^+$. This
implies that the pushforward $T_{\ast}L(x,\cdot)$ of $L(x,dy)$ under
$T$ is an infinitely divisible distribution supported on $S_d^+$. By
the L\'evy--Khintchine formula on proper cones (see \cite{skorohod},
Theorem 3.21, and by \cite{sato}, Proposition 11.10) this implies
that for all $x \in S_d^+$
%
%e4.9 ###
%e4.8 ###
%e4.7 ###
%
\begin{eqnarray}
\label{eq: TAT}
TA(x)T^\top&=&0,\\
\label{eq: TBorig}
TB(x)+\int_{S_d^+\setminus\{0\}} \bigl(\widetilde{\chi}(T\xi)-T (\chi
(\xi))\bigr) M(x, d \xi) &\in& S_d^+,\\
\label{eq: projLevy}
\int_{S_d^+\setminus\{0\}} (\|\xi\| \wedge1)T_{\ast}M(x,d\xi)&<&
\infty,
\end{eqnarray}
where $\widetilde{\chi}$ denotes some truncation function associated
with $T_{\ast}L(x,\cdot)$ and $T_{\ast}M$ the pushforward of $M$
under $T$. Due to (\ref{eq: projLevy}), we can set
$\widetilde{\chi}=0$. Thus, (\ref{eq: TBorig}) becomes
%
%e4.10 ###
%
\begin{equation}\label{eq: TB}
TB(x)-\int T (\chi(\xi)) M(x, d \xi) \in S_d^+.
\end{equation}
Moreover, equations (\ref{eq: TAT}), (\ref{eq: TB}) and (\ref{eq:
projLevy}) are equivalent to
%
%e4.11 ###
%
\begin{eqnarray}\label{eq: projLevy 1}
\langle T^\top v, A(x) T^\top v\rangle&=&0 \qquad\mbox{for all $v \in
S_d$},\nonumber\\
\hspace*{32pt}\langle B(x), T^\top v\rangle-\int_{S_d^+\setminus\{0\}}
\langle(\chi(\xi)), T^\top v\rangle M(x, d
\xi) &\geq&0 \qquad\mbox{for all $v
\in S_d^+$},\nonumber\\
\int_{S_d^+\setminus\{0\}}(\|T\xi\| \wedge1)M(x,d\xi)&<&
\infty.
\end{eqnarray}
In particular, we claim that
%
%e4.14 ###
%e4.13 ###
%e4.12 ###
%
\begin{eqnarray}
\label{eq: zero divisor zero}
\langle u, A(x) u\rangle&=&0\nonumber\\[-8pt]\\[-8pt]
&&\eqntext{\mbox{for all } u \in S_d \mbox{ s.t.
} ux=xu=0,}\\
\label{eq: zero divisor zeroB}
\hspace*{-52pt}\langle B(x), u\rangle- \int_{S_d^+\setminus\{0\}}
\langle\chi(\xi),u\rangle M(x, d \xi) &\geq&
0\nonumber\\[-8pt]\\[-8pt]
\eqntext{\mbox{for all } u\in S_d^+ \mbox{ s.t.
}ux=xu=0,}\\
\label{eq: zero divisor levy}
\int_{S_d^+\setminus\{0\}}\langle\chi(\xi),u \rangle M(x,d\xi)&<&
\infty \nonumber\\[-8pt]\\[-8pt]
&&\eqntext{\mbox{for all } u \in S_d^+ \mbox{ s.t. }
ux=xu=0.}
\end{eqnarray}
Indeed, if $x$ is invertible then $ux=0$ is equivalent to $u=0$ and
the assertions are obvious. Otherwise, if $x$ is in $\partial
S_d^+$, the linear map $T_u$ defined in Lemma \ref{lem: linear map}
is self-adjoint and satisfies $T_u(S_d^+-\re_+x)\subseteq
S_d^+$. Furthermore, by Lemma \ref{lem: linear map}(ii),
there exists an element $v \in S_d$ such that $T_uv=u$.
Hence,
\[
\langle u, A(x) u\rangle=\langle T_u^\top v, A(x) T_u^\top
v\rangle=0.
\]
It follows from the proof of Lemma \ref{lem: linear
map} that for $u \in S_d^+$, $v$ is an element of $S_d^+$ as well
and we have $\langle B(x), u\rangle=\langle B(x),T^{\top}_u
v\rangle$ and
$\langle\chi(\xi),u\rangle=\langle(\chi(\xi)),
T^\top v\rangle$. Equation (\ref{eq: zero divisor levy}) is
obtained by choosing $T=T_{\sqrt{u}}$ in (\ref{eq: projLevy 1}).
Indeed,
\begin{eqnarray*}
\int_{S_d^+ \cap\{\|\xi\| \leq1\}}\langle\xi,u \rangle M(x,d\xi)&=&
\int_{S_d^+ \cap\{\|\xi\| \leq1\}}\langle I_d, \xi u \rangle
M(x,d\xi)\\
&\leq&\|I_d\|\int_{S_d^+ \cap\{\|\xi\| \leq1\}} \|\xi u\|M(x,d\xi
)\\
&=& \|I_d\|\int_{S_d^+ \cap\{\|\xi\| \leq1\}}
\bigl\|T_{\sqrt{u}}\xi\bigr\|M(x,d\xi)< \infty.
\end{eqnarray*}
From these arguments and Lemma \ref{lem: zero divisor},
properties (\ref{eq: betaij}) and (\ref{eq: mu}) can be derived so
far. Thus, only (\ref{eq: alpha}) remains to be shown.

\textit{Step} 3. \textit{Necessary admissibility condition for $\alpha$.}
Due to the linearity of $A(x)$, $\langle u, A(x) u\rangle$ can be
written as $4\langle x, \vartheta(u)\rangle$, where the $(ij)$th
component of $ \vartheta(u)\in S_d$ is defined by
$\vartheta_{ij}(u)=1/4\langle u, a_{ij} u\rangle$. Note that
$\vartheta$ is defined on all of $S_d$. Given that for all $x \in
S_d^+$, $A(x)$ is a positive semidefinite operator on $S_d$,
$\langle u, A(x) u\rangle\geq0$ for all $u \in S_d$ and therefore,
by the self duality of $S_d^+$, $\vartheta(u)\in S_d^+$.
By (\ref{eq: zero divisor zero}), we have for all $u$ such that
$ux=xu=0$
%
%e4.15 ###
%
\begin{equation}\label{eq: vartheta}
0=\langle u, A(x) u\rangle=4\langle x, \vartheta(u)\rangle.
\end{equation}

Next, we show that $\vartheta$ is quasi-constant, that is, $
\langle x , \vartheta(u+w) - \vartheta(u) \rangle= 0 $ for all $
x,u,w \in S_{d}^+ $ with $ \langle x, w \rangle= 0 $ (see
Definition \ref{def: quasimono}). Indeed, pick $ x,u,w \in S_{d}^+ $
with $ \langle x, w \rangle= 0 $. According to our assumptions,
$A(x)w=0$, due to (\ref{eq: vartheta}) and the positivity of $A$.
Hence,
\begin{eqnarray*}
4\langle x,\vartheta(u+w) - \vartheta(u)\rangle&=&\langle u+w,
A(x)(u+w)\rangle-\langle u, A(x)u\rangle\\&=& \langle u,
A(x)w\rangle+\langle A(x)w, u\rangle=0,
\end{eqnarray*}
where the second last equality holds in view of the symmetry of
$A(x)$.

We now claim that there exists some $\alpha\in S_d^+$ such that
$\vartheta(u)=u\alpha u$, for each \mbox{$u\in S_d$}. It is sufficient to
show that this statements holds for all orthogonal projectors $ p
\in S_d^+ $, that is, there exists some $\alpha\in S_d^+$ such that $
\vartheta(p) = p
\alpha p $ for all orthogonal projectors $p$. Indeed, if this is the
case, we can derive the general statement in the following way: take
$ u \in S_d^+ $, then---by spectral decomposition---there are
numbers $ \lambda_i \geq0 $ and orthogonal projectors $ p_i $,
which are mutually orthogonal, such that $ u = \sum_{i=1}^d
\lambda_i p_i $ (see, e.g., Kato \cite{kat95}, Section I.6.9). Since
the assertion holds for all orthogonal projectors, we have that
\[
2\vartheta(u) = \sum_{i,j=1}^d \lambda_i \lambda_j \bigl(\vartheta(p_i +
p_j) - \vartheta(p_i) - \vartheta(p_j)\bigr)
\]
by the property that $ \vartheta$ is quadratic. Since $ p_i + p_j $
is again an orthogonal projector, we obtain the result.

We prove the assertion on orthogonal projectors by quasi-constancy.
Take an arbitrary orthogonal projector $ p $ and define $ q = I_d -
p $. Additionally, we define $ \alpha= \vartheta(I_d) $. By
quasi-constancy, we obtain
\[
\langle x , \vartheta(p + q) - \vartheta(q) \rangle= \langle y ,
\vartheta(p + q) - \vartheta(p) \rangle= 0
\]
and
\[
\langle x , \vartheta(p) \rangle= \langle y , \vartheta(q) \rangle
= 0,
\]
for all $ x,y \in S_d^+ $ with $ \langle x , p \rangle= 0 $ and $
\langle y , q \rangle= 0 $. Therefore, $ \alpha- \vartheta(q) $ and
$\vartheta(p) $ are orthogonal to the orthogonal complement of $ p $
in $ S_d^+ $ (i.e., the positive symmetric matrices of the form $ q
u q $ by Lemma \ref{lemprojector}), and $ \alpha- \vartheta(p) $
and $\vartheta(q) $ are orthogonal to the orthogonal complement of $
q $ in $ S_d^+ $ (the positive symmetric matrices of the form $ p u
p $ by Lemma \ref{lemprojector}). This means that we can write
\[
\alpha= \vartheta(p) + \vartheta(q) + \beta,
\]
where the symmetric matrix $ \beta$ is orthogonal to all elements
which are orthogonal to $ p $ and $ q $ (in $ S_d^+$), that is, $ \beta
$ is orthogonal to the linear span of matrices of the form $ pup $
and $ q u q $. However, such a decomposition is unique, since all
vectors in the sum are mutually orthogonal, and the decomposition is
given by
\[
\alpha=(p+q)\alpha(p+q)= p \alpha p + q \alpha q + (p \alpha q + q
\alpha p).
\]
Therefore, we can conclude the assertion $ \vartheta(p) = p \alpha p
$. Since $ p $ was arbitrary the assertion is proved.

Finally, all the derived restrictions on the parameters together
with (\ref{eq: levy khintchine R}) then yield (\ref{eq: R}).
\end{pf*}
\begin{remark}
An alternative proof for the special form of the diffusion matrix
$A(x)$ can also be established by Stokes' theorem
\cite{stokes} on Riccati ODEs.
\end{remark}

%s4.2 ###
\subsection{Infinitesimal generator}\label{section: inf generator}

The aim of this section is to prove the form of the infinitesimal
generator as stated in (\ref{eq: generator}).
\begin{proposition}\label{th: generator}
The infinitesimal generator $\Acal$ of an affine process on $S_d^+$
satisfies $\Scal_+\subset\Dcal(\Acal)$ and is of the form (\ref
{eq: generator}) for all $f \in\Scal_+$ and $x \in S_d^+$.
\end{proposition}
\begin{pf}
As already mentioned in the proof of Proposition \ref{th: necessary
admissibility}, the $t$-derivative of $P_te^{-\langle u,x\rangle}$
at $t=0$ exists pointwise for all $ x,u \in S_d^+$ and is given by
(\ref{Asharpeux}). Furthermore, $x\mapsto(-F(u)-\langle
R(u),x\rangle)e^{-\langle u,x\rangle}\in C_0(S_d^+)$, for $u \in
S_d^{++}$. As $(P_t)$ is a Feller semigroup on $C_0(S_d^+)$, it
follows from \cite{sato}, Lemma 31.7, that $\{e^{-\langle
u,x\rangle} \mid u \in S_d^{++}\} \in D(\mathcal{A})$ and
\[
\mathcal{A}e^{-\langle u,x\rangle}=\bigl(-F(u)-\langle
R(u),x\rangle\bigr)e^{-\langle u,x\rangle}.
\]
Combined with Proposition \ref{th: necessary admissibility}, we thus
obtain
%
%e4.16 ###
%
\begin{eqnarray}\label{Aexpux}\qquad
\mathcal{A}e^{-\langle u, x\rangle}&=&\biggl(-\langle b, u\rangle- c + \int
_{S_d^+\setminus\{0\}}\bigl(e^{-\langle u, \xi\rangle}-1\bigr)m(d\xi
)\nonumber\\
&&\hspace*{4.5pt}{}+\biggl\langle2u\alpha
u-B^\top(u)-\gamma\nonumber\\
&&\hspace*{20pt}{}+\int_{S_d^+\setminus\{0\}}
\biggl(\frac{e^{-\langle u, \xi\rangle}-1+\langle
\chi(\xi), u \rangle}{\|\xi\|^2 \wedge1}\biggr)\mu(d \xi),x\biggr\rangle
\biggr)e^{-\langle u,x\rangle}\nonumber\\
&=&\frac{1}{2}\sum_{i,j,k,l}A_{ijkl}(x) u_{ij}u_{kl}e^{-\langle
u,x\rangle}+\bigl\langle b+B(x), \nabla e^{-\langle u,x\rangle}
\bigr\rangle\\
&&{} -
(c+\langle\gamma,x\rangle)e^{-\langle u,x\rangle}\nonumber\\
&&{}+ \int_{S_d^+\setminus\{0\}}\bigl(e^{-\langle u, x+\xi\rangle
}-e^{-\langle u, x\rangle}\bigr)m(d\xi)\nonumber\\
&&{}+\int_{S_d^+\setminus\{0\}}\bigl(e^{-\langle u,
x+\xi\rangle}-e^{-\langle u, x\rangle}
+\bigl\langle
\chi(\xi), \nabla e^{-\langle u,x\rangle}
\bigr\rangle\bigr)M(x,d \xi).\nonumber
\end{eqnarray}
Indeed, in order to obtain the form of the diffusion part, observe
that we have by symmetrization
\begin{eqnarray*}
2\langle u \alpha u, x\rangle&=&2\sum_{i,j,k,l}\alpha_{jk}x_{il}
u_{ij}u_{kl}\\
&=&\frac{1}{2}\sum_{i,j,k,l}(x_{ik}\alpha_{jl}+x_{il}\alpha
_{jk}+x_{jk}\alpha_{il}+x_{jl}\alpha_{ik})u_{ij}u_{kl}\\
&=&\frac{1}{2}\sum_{i,j,k,l}A_{ijkl}(x) u_{ij}u_{kl};
\end{eqnarray*}
see (\ref{eq: cijkl}).
\goodbreak
According to Theorem \ref{th: density}, the linear hull $\mathcal M$
of $ \{e^{-\langle u, \cdot\rangle} \mid u \in
S_d^{++}\} $ is dense in $\mathcal S_+$ with respect to the
family of seminorms $p_{k,+}$ defined in (\ref{eq: seminorm}).
Denoting the right-hand side of (\ref{eq: generator}) by
$\mathcal{A}^{\sharp}$, we now claim that for every $f \in\mathcal
S_+$
%
%e4.17 ###
%
\begin{equation}\label{eq: Asharp}
{\lim_{n\rightarrow\infty}}\|
\mathcal{A}^{\sharp}f_n-\mathcal{A}^{\sharp}f\|_{\infty}=0,
\end{equation}
where $f_n \in\mathcal{M}$ such that $\lim_{n\rightarrow\infty}
p_{k,+}(f-f_n)=0$ for every $k$. Indeed, this is obvious for the
differential operator part of $\mathcal{A}^{\sharp}$. By choosing
$\chi(\xi)=1_{\{\|\xi\|\leq1\}}\xi$ and by denoting
$g(x):=f_n(x)-f(x)$, we obtain %\vadjust{\goodbreak}
the following estimate for the
integral part:
\begin{eqnarray*}
&&\biggl\|\int_{S_d^+\setminus\{0\}} \biggl(\frac{g(x+\xi)-g(x)-\langle1_{\{\|
\xi\|\leq1\}}\xi, \nabla g(x)\rangle}{\|\xi\|^2\wedge1}\biggr)
x_{ij}\mu_{ij}(d\xi)\biggr\|_{\infty}\\
&&\qquad\leq\int_{S_d^+\setminus\{0\}\cap\{\|\xi\|\leq1\}}\biggl\|\biggl(\sum
_{k,l,m,n}\biggl(\int_0^1\frac{\partial^2 g(x+s\xi)}{\partial
x_{kl}\,\partial x_{mn}} (1-s)\,ds\biggr)\frac{\xi_{kl}\xi_{mn}}{\|\xi\|
^2}\biggr)x_{ij}\biggr\|_{\infty}\\
&&\qquad\quad\hspace*{64pt}{}\times\bigl(\mu^+_{ij}(d\xi)+\mu^-_{ij}(d\xi)\bigr)\\
&&\qquad\quad{}+\int_{S_d^+\setminus\{0\}\cap\{\|\xi\|>1\}}\bigl\|\bigl(g(x+\xi
)-g(x)\bigr)x_{ij}\bigr\|_{\infty}\bigl(\mu^+_{ij}(d\xi)+\mu^-_{ij}(d\xi)\bigr)\\
&&\qquad\leq C_1p_{3,+}(g)\int_{S_d^+\setminus\{0\}\cap\{\|\xi\|\leq1\}
}\frac{\|\xi\|^2}{\|\xi\|^2}\bigl(\mu^+_{ij}(d\xi)+\mu^-_{ij}(d\xi)\bigr)\\
&&\qquad\quad{} +\int_{S_d^+\setminus\{0\}\cap\{\|\xi\|>1\}}\bigl(\bigl\|g(x+\xi)(\|x+\xi\|
)\bigr\|_{\infty}+\|g(x)x_{ij}\|_{\infty}\bigr)\\
&&\qquad\quad\hspace*{77pt}{}\times\bigl(\mu^+_{ij}(d\xi)+\mu
^-_{ij}(d\xi)\bigr)\\
&&\qquad\leq C_2p_{3,+}(g)\bigl(\mu^+_{ij}(S_d^+)+\mu^-_{ij}(S_d^+)\bigr)\leq C_3,
p_{3,+}(g),
\end{eqnarray*}
where $C_1, C_2$ and $C_3$ denote some constants and $\mu_{ij}^+,
\mu_{ij}^-$ correspond to the Jordan decomposition
$\mu_{ij}=\mu_{ij}^+-\mu_{ij}^-$. In the second last inequality, we
use the estimate $x_{ij}\leq\|x+\xi\|$. The same as above can be
shown for the measure $m(d\xi)$, whence (\ref{eq: Asharp}) holds
true. As by the first part of the proof, we have
$\mathcal{A}^{\sharp}=\mathcal{A}$ for all elements of
$\mathcal{M}$, (\ref{eq: Asharp}) implies
\[
{\lim_{n\rightarrow\infty}}\|
\mathcal{A}f_n-\mathcal{A}^{\sharp}f\|_{\infty}=0.
\]
Since the infinitesimal generator of every Feller process is a
closed operator, it follows that $\mathcal S_+\subset\Dcal(\Acal)$
and $\mathcal{A}=\mathcal{A}^{\sharp}$ on $\mathcal S_+$.
\end{pf}

%s4.3 ###
\subsection{Linear transformations and canonical representation}

In this subsection, we shall deal with linear transformations of
affine processes. The proposition below states how the parameters of
an affine process on $S_d^+$ change under such linear maps, which
allows us to establish a canonical representation of an affine
process.
\begin{proposition}\label{prop: lin trafo}
Suppose $X$ is an affine process on $S_d^+$ with parameters $\alpha,
\beta^{ij}, c$, $\gamma, m, \mu$ as specified in
Definition \ref{def: necessary admissibility} and $b \in S_d^+$.
Furthermore, let $G\dvtx S_d^+\rightarrow S_d^+, x\mapsto g x g^\top$
be an automorphism, where $g\in M_d$ is invertible. Then, $Y:=g
X g^\top$ is an affine process on $S_d^+$, whose parameters, %\vadjust{\goodbreak}
denoted by $\widetilde{\cdot}$, are given as follows with respect to
the truncation function $\widetilde{\chi}=g\chi(g^{-1}\xi
(g^\top)^{-1})g^\top$:
\begin{eqnarray*}
\widetilde{b}&=&gbg^\top,\\
\widetilde{c}&=&c, \\
\widetilde{m}(d\xi)&=&G_{\ast} m(d\xi),\\
\widetilde{\alpha}&=&g \alpha g^\top,\\
\widetilde{\gamma}&=&(g^\top)^{-1}\gamma g^{-1},\\
\widetilde{\mu}(d\xi)&=&\biggl(\frac{\|\xi\|^2 \wedge1}{\|g^{-1}\xi
(g^\top)^{-1}\|^2 \wedge1}\biggr)(g^\top)^{-1}G_{\ast}\mu(d\xi)g^{-1},\\
\widetilde{B}^\top(u)&=&(g^\top)^{-1}B^\top(g^\top u g)g^{-1},
\end{eqnarray*}
where $G_{\ast}m$ ($G_{\ast}\mu$) is the pushforward of the measure
$m$ ($\mu$, resp.).
\end{proposition}
\begin{pf}
Let us consider the process
\[
Y_t^y=gX_t^{g^{-1}y(g^\top)^{-1}}g^\top,
\]
for which we have
\begin{eqnarray*}
&&\E[\exp(-\langle u, Y_t^y\rangle)]\\
&&\qquad=
\E\bigl[\exp\bigl(-\bigl\langle u,gX_t^{g^{-1}y(g^\top)^{-1}}g^\top\bigr\rangle \bigr)\bigr]\\
&&\qquad=\E\bigl[\exp\bigl(-\bigl\langle g^\top u g,X_t^{g^{-1}y(g^\top)^{-1}} \bigr\rangle\bigr)\bigr]\\
&&\qquad=\exp\bigl(-\phi(t, g^\top u g)-\langle\psi(t,g^\top u
g),g^{-1}y(g^\top)^{-1}\rangle\bigr)\\
&&\qquad=\exp\bigl(-\phi(t, g^\top u
g)-\langle(g^\top)^{-1}\psi(t,g^\top u
g)g^{-1},y\rangle\bigr).
\end{eqnarray*}
Define now $\widetilde{\phi}$ and $\widetilde{\psi}$ by
\[
\widetilde{\phi}(t,u)=\phi(t,g^\top u g) \quad\mbox{and}\quad
\widetilde{\psi}(t,u)=(g^\top)^{-1}\psi(t,g^\top ug)g^{-1},
\]
to see that $Y$ is an affine process on $S_d^+$. Using (\ref{eq:
F-Riccati}) and (\ref{eq: F}), we consequently obtain
\begin{eqnarray*}
\frac{\partial\widetilde{\phi}(t,u)}{\partial t}&=&\frac{\partial
\phi(t,g^\top ug)}{\partial t}\\
&=&F(\psi(t,g^\top ug))\\
&=&\langle b, \psi(t,g^\top ug)\rangle+ c - \int_{S_d^+\setminus\{0\}
}\bigl(e^{-\langle\psi(t,g^\top ug), \xi\rangle}-1\bigr)m(d\xi)\\
&=&\langle gbg^\top, \widetilde{\psi}(t,u)\rangle+ c - \int
_{S_d^+\setminus\{0\}}\bigl(e^{-\langle\widetilde{\psi}(t,u), g \xi
g^\top\rangle}-1\bigr)m(d\xi)\\
&=&\langle\widetilde b, \widetilde{\psi}(t,u)\rangle+ c
-\int_{S_d^+\setminus\{0\}}\bigl(e^{-\langle\widetilde{\psi}(t,u),
\xi\rangle}-1\bigr)G_{\ast}m(d\xi).
\end{eqnarray*}
Due to\vspace*{1pt} the uniqueness of the L\'evy--Khintchine decomposition, this
implies that $b$ transforms to $\widetilde{b}=g b g^\top$, $c$
remains constant and $m$ becomes %\vadjust{\goodbreak}
$\widetilde{m}(d\xi)=G_{\ast}
m(d\xi)$. For $\widetilde{\psi}$ we proceed similarly, that is, we
have
\begin{eqnarray*}
\frac{\partial\widetilde{\psi}(t,u)}{\partial t}&=&(g^\top
)^{-1}\frac{\partial\psi(t,g^\top ug)}{\partial t}g^{-1}=(g^\top
)^{-1}R(\psi(t,g^\top ug))g^{-1}\\
&=&(g^\top)^{-1}\biggl(-2\psi(t,g^\top ug)\alpha\psi(t,g^\top ug)+B^\top
(\psi(t,g^\top ug))+\gamma\\
&&{} -\int_{S_d^+\setminus\{0\}}
\biggl(\frac{e^{-\langle\psi(t,g^\top ug),
\xi\rangle}-1+\langle\chi(\xi), \psi(t,g^\top
ug) \rangle}{\|\xi\|^2 \wedge1}\biggr)
\mu(d \xi)\biggr)g^{-1},
\end{eqnarray*}
from which it can be seen that $\alpha$ transforms to
$\widetilde{\alpha}=g\alpha g^\top$, $\gamma$ becomes
$\widetilde{\gamma}=(g^\top)^{-1}\gamma g^{-1}$, and $\mu$ changes
to
\[
\widetilde{\mu}(E)=(g^\top)^{-1}\biggl(\int_E\biggl(\frac{\|\xi\|^2
\wedge1}{\|g^{-1}\xi(g^\top)^{-1}\|^2 \wedge
1}\biggr)G_{\ast}\mu(d\xi)\biggr)g^{-1}
\]
for every $E \in\mathcal{B}(S_d^+\setminus\{0\})$. Moreover, since
$\widetilde{\chi}=g\chi(g^{-1}\xi(g^\top)^{-1})g^\top$
%
%e4.18 ###
%
\begin{equation}\label{eq: Btrans}
\widetilde{B}^\top(u)=(g^\top)^{-1}B^\top(g^\top u g)g^{-1}.
\end{equation}
\upqed\end{pf}

By means of Proposition \ref{prop: lin trafo}, we can derive a
canonical representation for affine processes.
\begin{proposition} \label{prop: canonic}
Let $X$ be an affine process on $S_d^+$ with parameters $\alpha,
\beta^{ij}, c, \gamma$, $m, \mu$ as specified in
Definition \ref{def: necessary admissibility} and $b \in S_d^+$.
Then there exists an automorphism $G\dvtx S_d^+ \rightarrow S_d^+,
x\mapsto gxg^\top$ such that the parameters of the affine process
$Y=gXg^\top$, denoted by $\widetilde{\cdot}$, are as in
Proposition \ref{prop: lin trafo} with
\[
\widetilde{b}=\theta=\diag(\theta_{11},\ldots,\theta_{dd}),\qquad
\widetilde{\alpha}=I_r^d,
\]
where we define
\[
I_r^d=\pmatrix{I_r & 0\cr
0 & 0}.
\]
\end{proposition}
\begin{pf}
By Proposition \ref{prop: lin trafo}, the parameters of $Y=gXg^\top$
transform as
\[
\widetilde{\alpha}=g\alpha g^\top,\qquad \widetilde{b}=gbg^\top.
\]
Since $\alpha$ and $b \in S_d^+$, they are jointly diagonalizable
through an automorphism on~$S_d^+$. More precisely, there exists an
invertible matrix $g \in M_d$ such that
\[
g \alpha g^\top= I_r^d\qquad \mbox{with } r=\rk(\alpha) %\vadjust{\goodbreak}
\]
and
\[
gbg^\top=\diag(\theta_{11},\ldots,\theta_{dd})=:\theta,
\]
where $\rk$ denotes the rank of a matrix. For the proof of this
fact, we refer to \cite{golub}, Theorem 8.7.1.
\end{pf}

%s4.4 ###
\subsection{Condition on the constant drift}\label{sec: constant drift}

This subsection is devoted to show that condition (\ref{eq: b})
holds true for any affine process $X$ on $S_d^+$. Since the
automorphism $G\dvtx S_d^+ \rightarrow S_d^+$ in Proposition \ref{prop:
canonic} is order preserving, it suffices to consider affine
processes of the canonical form as specified in
Proposition \ref{prop: canonic}. The following result is a
consequence of the L\'evy--Khintchine formula on $\R_+$.
\begin{lemma}\label{lemprop: characteristicsdet}
Let $Y$ be an affine process of canonical form as specified in
Proposition \ref{prop: canonic} with parameters denoted by
$\widetilde{\cdot}$. Then, for any $y \in\partial S_d^+$, we have
%
%e4.19 ###
%
\begin{equation}\label{Mtild0inf}
\nabla\det(y)\in N_{S_d^+}(y),\qquad \int_{S_d^+\setminus\{0\}}
\langle\widetilde{\chi}(\xi), \nabla
\det(y)\rangle\widetilde{M}(y,d\xi)<\infty
\end{equation}
and
%
%e4.20 ###
%
\begin{eqnarray}\label{eq: inv det explicit}\qquad
&&\langle\theta, \nabla\det(y)\rangle+\langle\widetilde{B}(y),
\nabla\det(y)\rangle-\int_{S_d^+\setminus\{0\}}\langle\widetilde
{\chi}(\xi), \nabla\det(y)
\rangle\widetilde{M}(y,d\xi)\nonumber\\[-8pt]\\[-8pt]
&&\qquad{}+\frac{1}{2}\sum_{i,j,k,l}\widetilde{A}_{ijkl}(y)\,\partial
_{ij}\,\partial_{kl}\det(y)
\geq0.\nonumber
\end{eqnarray}
\end{lemma}
\begin{pf}
Let $y \in\partial S_d^+$ and let $f \in C_c^{\infty}(S_d^+)$ be a
function with $f \geq0$ and $f(x)=\det(x)$ for all $x$ in a
neighborhood of $y$. Then, for any $v \in\re_+$, the function
$x\mapsto e^{-vf(x)}-1$ lies in $C_c^{\infty}(S_d^+)$ and thus in
$\Dcal(\mathcal{\widetilde{A}})$, where $\widetilde{\Acal}$ denotes
the infinitesimal generator of $Y$. Note that $f(y)=0$. Hence, the
limit
\begin{eqnarray*}
\mathcal{\widetilde{A}}\bigl(e^{-vf(y)}-1\bigr) &=&\lim_{t\rightarrow
0^+}\frac{1}{t}\int_{S_d^+}\bigl(e^{-vf(\xi)}-1\bigr)\widetilde{p}_t(y,d\xi)\\
&=&\lim_{t\rightarrow
0^+}\frac{1}{t}\int_{\re_+}(e^{-vz}-1)p^f_t(y,dz),
\end{eqnarray*}
exists for any $v\in\R_+$, where $\widetilde{p}_t(y,d\xi)$ denotes
the transition function of $Y$, and
$p_t^f(y,dz)=f_{\ast}\widetilde{p}_t(y,dz)$ is the pushforward of
$\widetilde{p}_t(y,\cdot)$ under $f$, which is a probability measure
supported on $\re_+$.

Using the same arguments as in Proposition \ref{th: necessary
admissibility} [i.e., applying Lemma \ref{lem: FourierLaplace} as
done below (\ref{eq: FRinfdiv})], and noting that $f(y)=0$,
we conclude that
%
%e4.21 ###
%
\begin{eqnarray}\label{Atildelk}
v&\mapsto&\mathcal{\widetilde{A}}\bigl(e^{-vf(y)}-1\bigr)\nonumber\\
&=&\frac{1}{2}\sum_{i,j,k,l}\widetilde{A}_{ijkl}(y) \bigl(v^2\,
\partial_{ij} f(y)\,\partial_{kl}f(y)-v\,
\partial_{ij}\,\partial_{kl}f(y)\bigr)\\
&&{}-v \langle\theta+\widetilde{B}(y),\nabla f(y)\rangle\nonumber\\
&&{} +\int_{S_d^+\setminus\{0\}}
\bigl(e^{-vf(y+\xi)}-1 \bigr)\widetilde{m}(d\xi)\nonumber\\
&&{} +\int_{S_d^+\setminus\{0\}} \bigl(e^{-vf(y+\xi)}-1+v \langle
\widetilde{\chi}(\xi), \nabla
f(y)\rangle\bigr)\widetilde{M}(y,d\xi)\nonumber
\end{eqnarray}
is the logarithm of the Laplace transform of an infinitely divisible
distribution on $\R_+$. Note that
\[
\langle\nabla\det(y),x\rangle=\frac{d}{dt}
\det(y+tx)\bigg|_{t=0}\cases{\ge0,&\quad $x\in S_d^+$, \cr
=0,&\quad $x=y$.}
\]
Hence, $\nabla\det(y)\in N_{S_d^+}(y)$ and the
admissibility condition (\ref{eq: mu}) implies (\ref{Mtild0inf}). By
the L\'evy--Khintchine formula on $\R_+$ (see \cite{skorohod},
Theorem 3.21), the linear coefficient in $v$ in (\ref{Atildelk})
has to be nonpositive. But this is now just (\ref{eq: inv det
explicit}), whence the lemma is proved.
\end{pf}

It now remains to show that (\ref{eq: b}) follows from (\ref{eq: inv
det explicit}). For this purpose, it suffices to
evaluate (\ref{eq: inv det explicit}) at diagonal elements $y\in
\partial S_d^+$. Thus, we state the following lemma.
\begin{lemma}\label{lem: quadraticvar}
Let $y \in S_d^+$ be diagonal, and let $f\in C_c^{\infty}(S_d^+)$.
Then we
have
\begin{eqnarray*}
&&\frac{1}{2}\sum
_{i,j,k,l=1}^d\bigl(y_{ik}(I_r^d)_{jl}+y_{il}(I_r^d)_{jk}+y_{jk}(I_r^d)_{il}+y_{jl}(I_r^d)_{ik}\bigr)\,\frac
{\partial^2
f(x)} {\partial x_{ij}\,\partial x_{kl}}\bigg|_{x=y}\\
&&\qquad=\frac{1}{2}\sum_{i,j=1}^d\bigl(y_{ii}1_{\{j \leq r\}}+y_{jj}1_{\{i \leq
r\}}\bigr)\biggl(\frac{\partial^2 f(x)}{\partial x_{ij}^2}\bigg|_{x=y}
+\frac{\partial^2 f(x)}{\partial x_{ij}x_{ji}}\bigg|_{x=y}\biggr).
\end{eqnarray*}
\end{lemma}
\begin{pf}
Obvious.
\end{pf}

Next, we calculate the partial derivatives of the determinant.
\begin{lemma}\label{lem: detder}
Let $y\in S_d^+$ be diagonal, $y=\diag(y_{11},y_{22},\ldots,y_{dd})$.
Then we have
\[
\frac{\partial\det(x)}{\partial x_{ij}}\bigg|_{x=y}=\cases{
\displaystyle \prod_{k\neq i}y_{kk}, &\quad if $i=j$,\cr
0, &\quad else,}
\]
and
\begin{eqnarray*}
\frac{\partial^2 \det(x)}{\partial x_{ij}x_{ji}}\bigg|_{x=y}&=& -\prod
_{k=1,k\neq i,k\neq j}^d y_{kk}\qquad\mbox{for $1\leq i<j\leq d$},\\
\frac{\partial^2 \det(x)}{\partial x_{ij}^2}\bigg|_{x=y}&=& 0
\qquad\mbox{for $1\leq i\leq j\leq d$},
\end{eqnarray*}
where the empty product is defined to be $1$.
\end{lemma}
\begin{pf}
In dimension $d=2$, the assertion is easily checked, as
$\det(y)=y_{11}y_{22}-y_{12}y_{21}$. Therefore, we have
\[
\partial_{11} \det(y)=y_{22},\qquad \partial_{22} \det(y)=y_{11},\qquad
\partial_{12} \det(y)=\partial_{21} \det(y)=0
\]
as well as
\begin{eqnarray*}
\partial_{11}^2 \det(y)&=&\partial_{22}^2 \det(y)=\partial_{12}^2
\det(y)=\partial_{21}^2 \det(y)=0, \\
\partial_{12}\,\partial_{21}\det(y)&=&\partial_{21}\,\partial_{12}\det(y)=-1.
\end{eqnarray*}
For dimension strictly larger than $2$, we employ a combinatorial
argument. Recall Leibniz's definition of the determinant,
%
%e4.22 ###
%
\begin{equation}\label{eq: leibnizdet}
\det(x)=\sum_{\sigma\in\Sigma}\sgn(\sigma)\prod
_{k=1}^dx_{k\sigma
(k)},
\end{equation}
where $\sigma$ is an element of the permutation group $\Sigma$ on
the set $\{1,2,\ldots,d\}$ and $\sgn{}$ denotes the signum function
on $\Sigma$, that is, $\sgn=1$ if $\sigma$ is an even permutation and
$\sgn=-1$ if it is odd. Differentiation of (\ref{eq: leibnizdet})
with respect to $x_{ij}$ yields
\[
\frac{\partial\det(x)}{\partial
x_{ij}}\bigg|_{x=y}=\biggl(\sum_{\sigma\in\Sigma}\sgn(\sigma)1_{\{\sigma
(i)=j\}}\prod_{k\neq
i}x_{k\sigma(k)}\biggr)\bigg|_{x=y}= \cases{
\displaystyle \prod_{k\neq i}y_{kk}, &\quad if $i=j$,\cr
0, &\quad else.}
\]
Thus, for the second derivative we have
\[
\frac{\partial^2 \det(x)}{\partial x_{ij}\,\partial
x_{ji}}\bigg|_{x=y}=\biggl(\sum_{\sigma\in\Sigma}\sgn(\sigma)1_{\{\sigma
(i)=j\}}1_{\{\sigma(j)=i\}}\prod_{k\neq
i\neq j}x_{k\sigma(k)}\biggr)\bigg|_{x=y}=-\prod_{k\neq i\neq
j}y_{kk},
\]
where the last equality holds since $y$ is diagonal. For
$\partial^2_{ij}\det(x)$, the statement is obvious.
\end{pf}

We are prepared to prove the admissibility condition on the constant
drift.
\begin{proposition}\label{th: constant drift}
Let $X$ be an affine process on $S_d^+$, then (\ref{eq: b}) holds,
that is,
\[
b \succeq(d-1)\alpha.
\]
\end{proposition}
\begin{pf}
Since the automorphism $G\dvtx S_d^+ \rightarrow S_d^+$ in
Proposition \ref{prop: canonic} is order preserving, it suffices to
show that (\ref{eq: inv det explicit}) in Lemma \ref{lemprop:
characteristicsdet} implies
%
%e4.23 ###
%
\begin{equation}\label{eq: canonic drift}
\theta\succeq(d-1)I_r^d.
\end{equation}
We show that $\theta_{mm}\geq d-1$, if $r \geq m$. To this end, take
again some diagonal $y \in\partial S_d^+$ of form
$y=\diag(y_{11}>0,\ldots,y_{mm}=0,\ldots, y_{dd}>0)$. By
Lemmas \ref{lem: quadraticvar} and \ref{lem: detder}, we obtain
\begin{eqnarray*}
&&\sum_{i=1}^d \theta_{ii}\,\partial_{ii}\det(y)+\sum
_{i,j}(\widetilde{B}(y))_{ij}\,\partial_{ij}\det(y)\\
&&\quad{} -\int_{S_d^+\setminus\{0\}}\biggl(\sum_{i,j}(\widetilde{\chi}(\xi
))_{ij}\,\partial_{ij} \det(y)\biggr)\widetilde{M}(y,d\xi)\\
&&\quad{} +\frac{1}{2}\sum_{i,j=1}^d\bigl(\bigl(y_{ii}1_{\{j \leq r\}}+y_{jj}1_{\{i
\leq r\}}\bigr)\bigl(\partial^2_{ij}\det(y)+\partial_{ij}\,\partial_{ji}\det
(y)\bigr)\bigr)\\
&&\qquad=\sum_{i=1}^d\biggl( \theta_{ii}\prod_{k\neq i}y_{kk}\biggr)+\sum_{l\neq
m}\biggl(\widetilde{\beta}_{mm}^{ll}y_{ll}\prod_{k\neq m}y_{kk}\biggr)\\
&&\qquad\quad{} -\sum_{l\neq m}\int_{S_d^+\setminus\{0\}}\frac{(\widetilde{\chi
}(\xi))_{mm}y_{ll}\prod_{k\neq m}y_{kk}}{\|\xi\|^2 \wedge
1}\widetilde{\mu}_{ll}(d\xi)\\
&&\qquad\quad{} -\frac{1}{2}\sum_{i \neq j}\biggl(\prod_{k\neq j}y_{kk}1_{\{j \leq r\}
}+\prod_{k\neq i}y_{kk}1_{\{i \leq r\}}\biggr)\\
&&\qquad=\theta_{mm}\prod_{k\neq m}y_{kk}+ \prod_{k\neq m}y_{kk}\biggl(\sum
_{l\neq m}\biggl(\widetilde{\beta}_{mm}^{ll}y_{ll}-y_{ll}\int
_{S_d^+\setminus\{0\}}\frac{(\widetilde{\chi}(\xi))_{mm}}{\|\xi\|
^2 \wedge1}\widetilde{\mu}_{ll}(d\xi)\biggr)\biggr)\\
&&\qquad\quad{} -(d-1)\prod_{k\neq m}y_{kk}1_{\{m \leq r\}} \geq0.
\end{eqnarray*}
As $\prod_{k\neq m}y_{kk}>0$ and by (\ref{eq: betaij}) also
\[
\biggl(\widetilde{\beta}_{mm}^{ll}y_{ll}-y_{ll}\int_{S_d^+\setminus\{0\}
}\frac{(\widetilde{\chi}(\xi))_{mm}}{\|\xi\|^2
\wedge1}\widetilde{\mu}_{ll}(d\xi)\biggr)\geq0
\]
for $l\neq m$, letting $y_{ll}\rightarrow0, l \neq m$ yields
$\theta_{mm}\geq d-1$ for $r \geq m$. Relabeling of indices then
proves (\ref{eq: canonic drift}).
\end{pf}

%s5 ###
\section{Sufficient conditions for the existence and uniqueness of
affine processes}\label{sec: existence}

In this section, we prove that for a given admissible parameter set
$\alpha$, $b$, $\beta^{ij}$, $c$, $\gamma$, $m$, $\mu$ satisfying the
conditions of Definition \ref{def: necessary admissibility}, there
exists a unique affine process on~$S_d^+$, whose infinitesimal
generator $\Acal$ is of form (\ref{eq: generator}). Our approach to
derive this result is to consider the martingale problem for the
operator $\mathcal{A}$. In order to prove uniqueness for this
martingale problem, we shall need the following existence and
uniqueness result for the generalized Riccati differential
equations (\ref{eq: F-Riccati}) and (\ref{eq: R-Riccati}).

%s5.1 ###
\subsection{Generalized Riccati differential equations}\label
{section: Riccati}

We first derive some properties of the function $R$ given in
(\ref{eq: R}).
\begin{lemma}\label{lem: R quasi monotone}
$R$ is analytic on $S_d^{++}$ and quasi-monotone increasing on
$S_d^+$.
\end{lemma}
\begin{pf}
That $R$ is analytic on $S_d^{++}$ follows by dominated convergence
(see, e.g., \cite{dfs}, Lemma A.2).

Now let $\delta>0$, and define
\begin{eqnarray*}
R^{\delta}(u)&=&-2u\alpha u+B^\top(u)+\gamma-\int_{\{\|\xi\|\ge
\delta\}} \biggl(\frac{e^{-\langle u,
\xi\rangle}-1+\langle\chi(\xi), u
\rangle}{\|\xi\|^2 \wedge
1}\biggr)\mu(d \xi)\\
&=&-2u\alpha u+\gamma+\biggl(B^\top(u)-\int_{\{\|\xi\|\ge
\delta\}}\frac{\langle\chi(\xi), u
\rangle}{\|\xi\|^2 \wedge1}\mu(d \xi)\biggr)\\
&&{} +\int_{\{\|\xi\|\ge\delta\}} \biggl(\frac{1-e^{-\langle
u, \xi\rangle}}{\|\xi\|^2 \wedge1}\biggr)\mu(d \xi).
\end{eqnarray*}
Now, the map $u\mapsto-2u\alpha u+\gamma$ is quasi-monotone
increasing, as it is shown in Step~3 of the proof of
Proposition \ref{th: necessary admissibility}. Furthermore, it
follows from the admissibility condition (\ref{eq: betaij}) that
\[
u\mapsto B^\top(u)-\int_{\{\|\xi\|\ge\delta\}} \frac{\langle
\chi(\xi), u \rangle}{\|\xi\|^2 \wedge1}\mu(d \xi)
\]
is a quasi-monotone increasing linear map on $S_d^+$. Finally, the
quasi-monotonicity of
\[
u\mapsto\int_{\{\|\xi\|\ge\delta\}} \biggl(\frac{1-e^{-\langle
u, \xi\rangle}}{\|\xi\|^2 \wedge1}\biggr)\mu(d \xi)
\]
is a consequence of the monotonicity of the exponential and that
$\supp(\mu)\subseteq S_d^+$.

By dominated convergence, we have $\lim_{\delta\rightarrow0}
R^\delta(u)= R(u)$ pointwise for each $u\in S_d^+$. Hence, the
quasi-monotonicity carries over to $R$. Indeed, choose $x, u, v\in
S_d^+$ such that $u\preceq v$ and $\langle v-u,x\rangle=0$. Then we
have for all $\delta$, $\langle R^\delta(v)-R^\delta(u)$, $x\rangle
\geq
0$. Thus,
\[
\langle R^\delta(v)-R^\delta(u),x\rangle\rightarrow\langle
R(v)-R(u),x\rangle\geq0,
\]
as $\delta\rightarrow0$, which proves that $R$ is quasi-monotone
increasing.
\end{pf}
\begin{lemma}\label{lemKestR}
There exists a constant $K$ such that
%
%e5.1 ###
%
\begin{equation}\label{eq: above estimate R(u) on cone}
\langle u,R(u)\rangle\leq\frac{K}{2} (\|u\|^2+1),\qquad u\in S_d^+.
\end{equation}
\end{lemma}
\begin{pf}
We may assume, without loss of generality, that the truncation
function in Definition \ref{def: necessary admissibility} takes the
form $\chi(\xi)=1_{\{\|\xi\| \leq1\}}\xi$ [otherwise adjust $B(u)$
accordingly]. Then, for all $u\in S_d^+$ we have
%
%e5.2 ###
%
\begin{eqnarray}\label{est: crucial for upper bounds}
R(u)&=&-2u\alpha
u+B^\top(u)+\gamma-\int_{S_d^+\setminus\{0\}\cap\{\|\xi\| \leq
1\}}\underbrace{\biggl(\frac{e^{-\langle
u,\xi\rangle}-1+\langle\xi,u\rangle}{\|\xi\|^2}\biggr)}_{\geq0}\mu
(d\xi)\hspace*{-22pt}\nonumber\\
&&{} -\int_{S_d^+\setminus\{0\}\cap\{\|\xi\| >
1\}}\bigl(e^{-\langle
u,\xi\rangle}-1\bigr)\mu(d\xi)\nonumber\\[-8pt]\\[-8pt]
&\preceq&-2u\alpha u+B^\top(u)+\gamma+\mu(S_d^+\cap\{\|\xi\| > 1\}
)\nonumber\\
&\preceq& B^\top(u)+\gamma+\mu(S_d^+\cap\{\|\xi\| > 1\}),\nonumber
\end{eqnarray}
where we use that
\[
-\int_{S_d^+\setminus\{0\}\cap\{\|\xi\| > 1\}} \bigl(e^{-\langle
u,\xi\rangle}-1\bigr)\mu(d\xi)\preceq
\int_{S_d^+\setminus\{0\}\cap\{\|\xi\| > 1\}} \mu(d\xi).
\]
Set now
\[
\overline{\gamma}:=\gamma+ \mu(S_d^+\cap\{\|\xi\| > 1\}) \in S_d^+.
\]
By (\ref{est: crucial for upper bounds}), we obtain, for $u\in
S_d^+$, that
\[
\langle u,R(u)\rangle\leq\langle u, B^\top(u)\rangle+\langle
u,\overline{\gamma}\rangle,
\]
from which we derive the existence of a positive constant $K$ such
equation (\ref{eq: above estimate R(u) on cone}) holds.
\end{pf}

Here is our main existence and uniqueness result for the generalized
Riccati differential equations (\ref{eq: F-Riccati}) and (\ref{eq: R-Riccati}).
\begin{proposition}\label{prop_ricc_sol}
For every $u\in S_d^{++}$, there exists a unique global $\R_+\times
S_d^{++}$-valued solution $(\phi,\psi)$ of (\ref{eq:
F-Riccati}) and (\ref{eq: R-Riccati}). Moreover, $\phi(t,u)$ and
$\psi(t,u)$ are analytic in $(t,u)\in\R_+\times S_d^{++}$.
\end{proposition}
\begin{pf}
We only have to show that, for every $u\in S_d^{++}$, there exists a
unique global $S_d^{++}$-valued solution $\psi$ of (\ref{eq:
R-Riccati}), as then $\phi$ is uniquely determined by integrating
(\ref{eq: F-Riccati}) and has the desired properties by
admissibility of the parameter set.

Let $u\in S_d^{++}$. Since $R$ is analytic on $S_d^{++}$, standard
ODE results (e.g., \cite{dieu69}, Theorem 10.4.5) yield there
exists a unique local $S_d^{++}$-valued solution $\psi(t,u)$ of
(\ref{eq: R-Riccati}) for $t\in[0,t_+(u))$, where
\[
t_+(u)=\lim\inf_{n\to\infty} \{ t\ge0\mid\|\psi(t,u)\|\ge n\mbox
{ or }
\psi(t,u)\in\partial S_d^+\}\le\infty.
\]
It thus remains to show
that $t_+(u)=\infty$. That $\psi(t,u)$, and hence $\phi(t,u)$, is
analytic in $(t,u)\in\R_+\times S_d^{++}$ then follows from
\cite{dieu69}, Theorem 10.8.2.

Since $R$ may not be Lipschitz continuous at $\partial S_d^+$ (see
Remark \ref{reminwpointing} below), we first have to regularize it.
We thus define
\[
\widetilde{R}(u)=-2u\alpha u+B^\top(u)+\gamma
-\int_{S_d^+\setminus\{0\}\cap\{\|\xi\| \leq
1\}}\biggl(\frac{e^{-\langle
u,\xi\rangle}-1+\langle\xi,u\rangle}{\|\xi\|^2}\biggr)\mu(d\xi).
\]
It then follows as in Lemmas \ref{lem: R quasi monotone} and
\ref{lemKestR} that $\widetilde{R}$ is quasi-monotone increasing on
$S_d^+$ and that (\ref{eq: above estimate R(u) on cone}) holds for
some constant $\widetilde{K}$. Moreover, $\widetilde{R}$ is analytic
on $S_d$. Hence, for all $u\in S_d$, there exists a unique local
$S_d$-valued solution $\widetilde{\psi}$ of
\[
\frac{\partial\widetilde{\psi}(t,u)}{\partial
t}=\widetilde{R}(\widetilde{\psi}(t,u)),\qquad
\widetilde{\psi}(0,u)=u ,
\]
for all $t\in[0,\widetilde{t}_+(u))$ with maximal lifetime
\[
\widetilde{t}_+(u)=\lim\inf_{n\to\infty} \{ t\ge0\mid\|
\widetilde{\psi}(t,u)\|\ge n\}\le\infty.
\]
From (\ref{eq: above estimate R(u) on cone}), we infer that for all
$u\in S_d^+$ and $t<\widetilde{t}_+(u)$,
\[
\partial_t \|\widetilde{\psi}(t,u)\|^2=2\langle
\widetilde{\psi}(t,u),\partial_t{\widetilde{\psi}}(t,u)\rangle
\leq
\widetilde{K} \bigl(\|\widetilde{\psi}(t,u)\|^2+1\bigr).
\]
Gronwall's inequality (e.g., \cite{dieu69}, (10.5.1.3)) implies
%
%e5.3 ###
%
\begin{equation}\label{eq: aprioribound}
\|\widetilde{\psi}(t,u)\|^2\leq e^{\widetilde{K}t}(\|u\|^2+1) ,\qquad
t<\widetilde{t}_+(u).
\end{equation}
Hence, $\widetilde{t}_+(u)=\infty$ for $u\in S_d^{+}$.
As $\widetilde{R}$ is quasi-monotone increasing on
$S_d^+$, Volkmann's comparison Theorem \ref{th: Volkmann} now yields
\[
0\preceq\widetilde{\psi}(t,u)\preceq\widetilde{\psi}(t,v),\qquad t\geq
0\mbox{, for all $0\preceq u\preceq v$.}
\]
Therefore and since $\widetilde{\psi}(t,u)$ is
also analytic in $u$, Lemma \ref{prop: feller prop} implies that
$\widetilde{\psi}(t,u) \in S_d^{++}$ for all $(t,u) \in\re_+ \times
S_d^{++}$.

We now carry this over to $\psi(t,u)$ and assume without loss of
generality, as in the proof of Lemma \ref{lemKestR}, that the
truncation function in Definition \ref{def: necessary admissibility}
takes the form $\chi(\xi)=1_{\{\|\xi\| \leq1\}}\xi$. Then
\[
R(u)-\widetilde{R}(u)=-\int_{S_d^+\setminus\{0\}\cap\{\|\xi\| > 1\}
}\bigl(e^{-\langle
u,\xi\rangle}-1\bigr)\mu(d\xi)\succeq0,\qquad u\in S_d^+.
\]
Hence,
for $u\in S_d^{++}$ and $t<t_+(u)$, we have
\[
\frac{\partial\widetilde{\psi}(t,u)}{\partial
t}-\widetilde{R}(\widetilde{\psi}(t,u)) = \frac{\partial
\psi(t,u)}{\partial t}-R(\psi(t,u)) \preceq\frac{\partial
\psi(t,u)}{\partial t}-\widetilde R(\psi(t,u)).
\]
Theorem \ref{th: Volkmann} thus implies
\[
\psi(t,u)\succeq\widetilde{\psi}(t,u) \in S_d^{++},\qquad t\in
[0,t_+(u)).
\]
Hence, $t_+(u)=\lim\inf_{n\to\infty} \{ t\ge0\mid
\|\psi(t,u)\|\ge n\}$. Using (\ref{eq: above estimate R(u) on cone})
again, we now can show as for $\widetilde{\psi}$ that
\[
\|\psi(t,u)\|^2\leq e^{Kt}(\|u\|^2+1) ,\qquad t<t_+(u).
\]
Hence $t^+(u)=\infty$, as desired.
\end{pf}
\begin{remark}\label{reminwpointing}
Lemma \ref{lem: R quasi monotone} states that the admissibility of
the parameters $\alpha,\beta^{ij},\gamma,\mu$ implies
quasi-monotonicity of $R$ on $S_d^+$.\footnote{We conjecture that
the converse also holds: $R$ is quasi-monotone on $S_d^+$ and
$\supp(\mu)\subseteq S_d^+$ if and only if the parameters
$\alpha,\beta^{ij},\gamma,\mu$ are admissible.} Moreover,
quasi-monotonicity just means that $R$ is ``inward pointing'' close
to the boundary $S_d^+$. Indeed, let $u,x\in S_d^+$ with $\langle
u,x\rangle=0$. Then $\langle R(u),x\rangle\geq\langle
\gamma,x\rangle\ge0$. Hence, if $R$ were Lipschitz continuous on
$S_d^+$, a deterministic variant of Theorem \ref{th: martingale
problem convex} would imply the invariance of $S_d^+$ with respect
to (\ref{eq: R-Riccati}) right away. However, the map $R$ might fail
to be Lipschitz at $\partial S_d^+$ (see the one-dimensional
counterexample~\cite{dfs}, Example 9.3), even though it is analytic
on the interior $S_d^{++}$. Here, quasi-monotonicity plays the decisive role.
It leads to the phenomenon that $\psi(t,u)$ stays away from the
boundary $\partial
S_d^+$ for $u\in S_d^{++}$, which is of crucial importance in our
analysis.
\end{remark}

%s5.2 ###
\subsection{The martingale problem for $\mathcal{A}$}

We are now prepared to study the martingale problem for the operator
$\mathcal{A}$ given by (\ref{eq: generator}). For the notion of
martingale problems, we refer to \cite{ethier}, Chapter 4. We shall
proceed in four steps. First, we approximate $\Acal$ by regular
operators $\Acal^{\varepsilon,\delta,n}$ on the space $\Scal_+$
of rapidly decreasing $C^\infty$-functions on $S_d^+$, defined in
(\ref{defS+}).
Second, using Theorem \ref{th: martingale
problem convex} below, we show that there exists an $S_d^+$-valued
c\`adl\`ag solution of the martingale problem for
$\Acal^{\varepsilon,\delta,n}$. Third, a subsequence of these solutions
is shown to converge to an $S_d^+\cup\{\Delta\}$-valued c\`adl\`ag
solution of the martingale problem for $\Acal$. Finally, we show
that this solution is unique, Markov and affine, as desired.

Note that we cannot employ Stroock's \cite{stroock1} seminal
existence and uniqueness results for martingale problems, since
those are solved on $\R^n$ and require uniform elliptic diffusion
parts. Neither of these is satisfied in our case.

Now let $(\alpha,b,\beta^{ij},c=0,\gamma=0,m,\mu)$ be some
admissible parameter set. Fix some $\varepsilon, \delta>0$ and
$n\in\N$. In order to bound the coefficients and cut off the small
jumps, we let
%
%e5.4 ###
%
\begin{equation}\label{eq: varphi}\hspace*{28pt}
\varphi_n\in C^\infty_b(S_d),\qquad 0\le\varphi_n\le
1,\qquad \varphi_n(x)=\cases{1,&\quad $\|x\|\leq n$,\vspace*{2pt}\cr
\dfrac{n}{\|x\|}, &\quad $\|x\|\geq n+1$.}
\end{equation}
We then define the bounded and smooth parameters
\begin{eqnarray*}
B^n(x)&=&B(\varphi_n(x)x),\\
m^{\delta}(d\xi)&=&m(d\xi)1_{\{\|\xi\| >\delta\}},\\
M^{\delta,n}(x,d\xi)&=& \biggl\langle\varphi_n(x)x,
\frac{\mu(d\xi)}{\|\xi\|^2 \wedge1}1_{\{\|\xi\|
>\delta\}}\biggr\rangle.
\end{eqnarray*}

Concerning the diffusion function $A_{ijkl}(x)$ given by (\ref{eq:
cijkl}), we first find an appropriate factorization which will allow
us to write the continuous martingale part of $X$ as a stochastic
integral. Thereto observe that any $S_d^+$-valued solution, presumed
that it exists, of the following symmetric matrix-valued diffusion
SDE:
%
%e5.5 ###
%
\begin{equation}\label{eq: Wishart diff}
dZ_t= \sqrt{Z_t}\,dW_t\,\Sigma+\Sigma^{\top}\,dW_t^{\top}\,\sqrt{Z_t},
\end{equation}
where $W$ is a standard $d \times d$-matrix Brownian motion and
$\Sigma\in M_d$ with $\Sigma^{\top}\Sigma=\alpha$, has quadratic
variation $d\langle Z_{ij}, Z_{kl}\rangle_t=A_{ijkl}(Z_t)$. Define
now $\sigma^{kl}(x) \in S_d$ by
%
%e5.6 ###
%
\begin{equation}\label{sigmakldef}
\sigma^{kl}(x)=\sqrt{x}M^{kl}\Sigma+\Sigma^{\top}M^{lk}\sqrt{x},
\end{equation}
where $M^{kl}_{ij}=\delta_{ik}\delta_{jl}$. Then (\ref{eq: Wishart
diff}) can be written as
\[
dZ_t= \sum_{k,l=1}^d \sigma^{kl}(Z_t)\,dW_{t,kl}
\]
and
$A_{ijkl}(x)=\sum_{m,n=1}^d\sigma^{mn}_{ij}(x)\sigma^{mn}_{kl}(x)$.

Since $\sigma^{kl}(x)$ involves the matrix square root, which is
neither Lipschitz continuous nor bounded nor globally defined, we
need to introduce some approximating regularization in order to meet
the assumptions of Theorem \ref{th: martingale problem convex} below.
Thereto fix some truncation function
\[
\eta_{\varepsilon}\in C^{\infty}_b(S_d),\qquad
\eta_{\varepsilon}(x)=\cases{1, &\quad $x\in S_d^+$,\cr
0,&\quad $x\notin
S_d^+ -\varepsilon I_d$,}
\]
and define
%
%e5.7 ###
%
\begin{equation}\label{eq: sigmavareps}
s_{\varepsilon,n}(x)=\cases{
\eta_{\varepsilon}(\varphi_n(x)x)\bigl(\sqrt{\varphi_n(x)x+\varepsilon
I_d}-\sqrt{\varepsilon I_d}\bigr), &\quad if $x \in S_d^+
-\varepsilon I_d$,\cr
0, &\quad otherwise.}\hspace*{-28pt}
\end{equation}
Note that $s_{\varepsilon,n}$ satisfies:
\begin{itemize}
\item$s_{\varepsilon,n}\in C^{\infty}_b(S_d,S_d)$,
\item$s_{\varepsilon,n}(x)=(\sqrt{\varphi_n(x)x+\varepsilon
I_d}-\sqrt{\varepsilon I_d})$ on $S_d^+$,
\item$\lim_{\varepsilon\rightarrow0^+}s_{\varepsilon,n}(x)=\sqrt
{\varphi_n(x)x}$.
\end{itemize}
With this, we can now define the regularization of $\sigma^{kl}$ by
%
%e5.8 ###
%
\begin{equation}\label{eq: reg sigmakl}
\sigma^{kl}_{\varepsilon,n}(x)=s_{\varepsilon,n}(x)M^{kl}\Sigma
+\Sigma^{\top}M^{lk}s_{\varepsilon,n}(x),
\end{equation}
which then satisfies the smoothness condition of Theorem \ref{th:
martingale problem convex}. Finally, we set
%
%e5.9 ###
%
\begin{eqnarray}\label{eq: approx c}
A^{\varepsilon,n}_{ijkl}(x)&=&\sum_{m,n}^d(\sigma^{mn}_{\varepsilon
,n}(x))_{ij}(\sigma^{mn}_{\varepsilon,n}(x))_{kl}\nonumber\\[-1.5pt]
&=&(s^2_{\varepsilon,n}(x))_{ik}\alpha_{jl}+(s^2_{\varepsilon
,n}(x))_{il}\alpha_{jk}\\[-1.5pt]
&&{}+(s^2_{\varepsilon,n}(x))_{jk}\alpha
_{il}+(s^2_{\varepsilon,n}(x))_{jl}\alpha_{ik},\nonumber
\end{eqnarray}
and define the corresponding regularized operator on $C_0(S_d)$
%
%e5.10 ###
%
\begin{eqnarray}
\label{eq: approx operator}
\mathcal{A}^{\varepsilon,\delta,n}f(x)&=&\frac{1}{2}\sum
_{i,j,k,l}A^{\varepsilon,n}_{ijkl}(x)\,
\frac{\partial^2
f(x)} {\partial x_{ij}\,\partial x_{kl}}\nonumber\\
&&{}+\sum
_{i,j}\bigl(b_{ij}+B^{n}_{ij}(x)\bigr)\,\frac{\partial f(x)}{\partial
x_{ij}}\nonumber\\[-8pt]\\[-8pt]
&&{} +\int_{S_d^+\setminus\{0\}}
\bigl(f(x+\xi)-f(x)\bigr)m^{\delta}(d\xi)\nonumber\\
&&{} +\int_{S_d^+\setminus\{0\}}\bigl(f(x+\xi)-f(x)-\langle
\chi(\xi), \nabla f(x)\rangle\bigr)
M^{\delta,n}(x,d\xi).\nonumber
\end{eqnarray}

We now show that $\mathcal{A}^{\varepsilon,\delta,n}$ approximates
$\Acal$. We let $\Scal=\Scal(S_d)$ and $\Scal_+$ denote the locally
convex spaces of rapidly decreasing $C^\infty$-functions on $S_d$
and $S_d^+$ defined in (\ref{defS+}) below, respectively.
\begin{lemma}\label{lemAconvA}
$\Scal\subset\Dcal(\mathcal{A}^{\varepsilon,\delta,n})$ and, for
every $f \in\Scal_+$,
%
%e5.11 ###
%
\begin{equation}\label{eq: gen convergence}
{\lim_{{\varepsilon, \delta, n}}} \|\mathcal{A}^{\varepsilon, \delta,
n}f-\mathcal{A}f\|_{\infty}=0.
\end{equation}
\end{lemma}
\begin{pf}
Since $\varphi_n$ as defined in (\ref{eq: varphi}) converges
uniformly on compact sets to~1, this is clear for the differential
operator part. Concerning the integral part, we have
\begin{eqnarray*}
&&\biggl\|\int_{S_d^+ \setminus\{0\}}\bigl(f(x+\xi)-f(x)-\bigl\langle1_{\{\|\xi\|
\leq1\}}\xi, \nabla f(x)\bigr\rangle\bigr) \bigl(M^{\delta,n}(x,d\xi)-M(x,d\xi
)\bigr)\biggr\|\\
&&\qquad\leq\biggl\|\sum_{i,j}\int_{S_d^+ \setminus\{0\}}\biggl(\frac{f(x+\xi
)-f(x)-\langle1_{\{\|\xi\|\leq1\}}\xi, \nabla f(x)\rangle}{\|\xi\|
^2 \wedge1}\biggr)\\
&&\qquad\quad\hspace*{126.6pt}{}\times x_{ij}\bigl(\varphi_n(x)-1\bigr)\mu_{ij}^{\delta}(d\xi)\biggr\|\\
&&\qquad\quad{}+\biggl\|\sum_{i,j}\int_{S_d^+
\setminus\{0\}}\biggl(\frac{f(x+\xi)-f(x)-\bigl\langle1_{\{\|\xi\|\leq
1\}}\xi, \nabla f(x)\rangle}{\|\xi\|^2 \wedge
1}\biggr)\\
&&\qquad\quad\hspace*{128pt}{}\times x_{ij}\bigl(1_{\{\|\xi\|>\delta\}}-1\bigr)\mu_{ij}(d\xi)\biggr\|.
\end{eqnarray*}
By dominated convergence the second term goes uniformly in $x$ to
$0$, thus we only have to consider the first one. By splitting the
first integral into $\int_{\{\|\xi\|\leq1\}}+\int_{\{\|\xi\|>
1\}}$, we note that $\|\int_{\{\|\xi\|\leq1\}}\|$ converges
uniformly in $x$ to $0$. Hence, it remains to analyze
\[
\biggl\|\sum_{i,j}\int_{\{\|\xi\|> 1\}}\bigl(f(x+\xi)-f(x)\bigr)
x_{ij}\bigl(\varphi_n(x)-1\bigr)\mu_{ij}(d\xi)\biggr\|,
\]
which can be estimated by
\begin{eqnarray*}
&&\sum_{i,j}\biggl(\int_{\{\|\xi\|> 1\}}\bigl\|f(x+\xi)x_{ij}\bigl(\varphi_n(x)-1\bigr)\bigr\|
\bigl(\mu_{ij}^+(d\xi)+\mu_{ij}^-(d\xi)\bigr)\\
&&\qquad\hspace*{4.5pt}{}+\int_{\{\|\xi\|>
1\}}\bigl\|f(x)x_{ij}\bigl(\varphi_n(x)-1\bigr)\bigr\|\bigl(\mu_{ij}^+(d\xi)+\mu_{ij}^-(d\xi)\bigr)\biggr),
\end{eqnarray*}
where $\mu_{ij}^+, \mu_{ij}^-$ correspond to the Jordan
decomposition of $\mu_{ij}=\mu_{ij}^+-\mu_{ij}^-$. As $f$ lies in
$\Scal_+$, the second term converges uniformly to $0$. For the first
one, observe that for every $n$
\[
\bigl\|f(x+\xi)x_{ij}\bigl(\varphi_n(x)-1\bigr)\bigr\|\leq\|f(x+\xi)x_{ij}\|\leq
\|f(x+\xi)\|\|x+\xi\|,
\]
such that we can apply dominated convergence. Again, since $f$ lies
in $\Scal_+$, the first integral converges uniformly in $x$ to $0$
as well. Hence (\ref{eq: gen convergence}) holds true, and
$\Scal\subset\Dcal(\mathcal{A}^{\varepsilon,\delta,n})$ follows
similarly.
\end{pf}

We now establish existence for the martingale problem for
$\mathcal{A}^{\varepsilon,\delta,n}$.
\begin{lemma}\label{lemprop: martingale problem approx}
For every $x\in S_d^+$ there exists an $S_d^+$-valued c\`adl\`ag
solution $X$ to the martingale problem for
$\mathcal{A}^{\varepsilon,\delta,n}$ with $X_0=x$. That is,
\[
f(X_t)-\int_0^t \mathcal{A}^{\varepsilon,\delta,n}f(X_s) \,ds
\]
is a martingale, for all $f\in\Scal$.
\end{lemma}
\begin{pf}
Consider the following SDE of type (\ref{eq: jump diffusion}):
%
%e5.12 ###
%
\begin{eqnarray}\label{eq: approx SDE}\qquad
X^{\varepsilon, \delta, n}_t&=&x+\int_0^t \biggl(b+B^n(X^{\varepsilon,
\delta, n}_s)
- \int_{{S_d^+}\setminus\{0\}} \chi(\xi)M^{\delta
,n}(X^{\varepsilon, \delta, n}_s,d\xi)\biggr) \,ds
\nonumber\\[-8pt]\\[-8pt]
&&{} + \sum_{k,l}^d
\int_0^t\sigma^{kl}_{\varepsilon,n}(X^{\varepsilon, \delta,
n}_s)\,dW_{s,kl}+J_t,\nonumber
\end{eqnarray}
where $W$ is a $d \times d$-matrix of standard Brownian motions
and $J$ a finite activity jump process with compensator
$m^{\delta}(d\xi)+M^{\delta,n}(X^{\varepsilon, \delta, n}_t,d\xi)$.
Note that the quadratic variation of the continuous martingale part
of $X^{\varepsilon, \delta, n}_t$ is given by
$A^{\varepsilon,n}_{ijkl}(x)$ as defined in (\ref{eq: approx c}). It
thus follows by inspection that any c\`adl\`ag solution
$X^{\varepsilon, \delta, n}$ of (\ref{eq: approx SDE}) solves the
martingale problem for $\mathcal{A}^{\varepsilon,\delta,n}$.

Hence, it remains to show that there exists an $S_d^+$-valued
c\`adl\`ag solution of (\ref{eq: approx SDE}). Let us recall the
normal cone (\ref{eq: normal cone Sd+}) to $S_d^+$. As $b+B^n(x)-
\int_{{S_d^+}\setminus\{0\}} \chi(\xi)\times M^{\delta,n}(x,d\xi)$,
$\sigma^{kl}_{\varepsilon,n}(x)$ and
$m^{\delta}(d\xi)+M^{\delta,n}(x,d\xi)$ are designed to satisfy the
assumptions of Theorem \ref{th: martingale problem convex} and since
$x+\supp(m^{\delta}(\cdot)+M^{\delta,n}(x,\cdot))\subseteq S_d^+$
for all $x \in S_d^+$, we only have to show that for all $x \in
\partial S_d^+$ and $u \in N_{S_d^+}(x)$
%
%e5.14 ###
%e5.13 ###
%
\begin{eqnarray}\label{eq: volapara}
&&\hspace*{125pt}\langle\sigma_{\varepsilon,n}^{kl}(x),u\rangle=0,\\
\label{eq: stratodrift cond}
&&\Biggl\langle b+B^n(x)- \int_{{S_d^+}\setminus\{0\}}
\chi(\xi)M^{\delta,n}(x,d\xi)\nonumber\\[-8pt]\\[-8pt]
&&\hspace*{45pt}{}-\frac{1}{2}\sum_{k,l=1}^d
D\sigma_{\varepsilon,n}^{kl}(x)\sigma_{\varepsilon
,n}^{kl}(x),u\Biggr\rangle
\geq0.\nonumber
\end{eqnarray}
Due to the definition of $\sigma_{\varepsilon,n}^{kl}(x)$,
respectively, the definition of $s_{\varepsilon,n}(x)$ given
in (\ref{eq: sigmavareps}), condition (\ref{eq: volapara}) is
satisfied. Concerning (\ref{eq: stratodrift cond}), we have
by (\ref{eq: betaij})
\[
\biggl\langle B^n(x)- \int_{{S_d^+}\setminus\{0\}}
\chi(\xi)M^{\delta,n}(x,d\xi),u\biggr\rangle\geq0.
\]
Moreover, it is shown in Lemma \ref{prop: strato} below that
\[
\Biggl\langle b-\frac{1}{2}\sum_{k,l=1}^d
D\sigma_{\varepsilon,n}^{kl}(x)\sigma_{\varepsilon
,n}^{kl}(x),u\Biggr\rangle
\geq0.
\]
 The lemma now follows from Theorem \ref{th: martingale
problem convex}.
\end{pf}
\begin{lemma}\label{prop: strato}
Let $x=O\Lambda O^\top\in S_d^+$ where
$\Lambda=\diag(\lambda_1,\ldots, \lambda_d)$ contains the
eigenvalues in decreasing order and let
$\sigma_{\varepsilon,n}^{kl}$ be defined by (\ref{eq: reg sigmakl}).
Then, for all \mbox{$x \in S_d^+$},
%
%e5.15 ###
%
\begin{eqnarray}\qquad
\label{eq: Strato corr}
\frac{1}{2}\sum_{k,l=1}^d
D\sigma_{\varepsilon,n}^{kl}(x)\sigma_{\varepsilon
,n}^{kl}(x)&=&\frac{1}{2}\sum_{i=1}^d
\frac{\varphi_n(x)(\sqrt{\lambda_i\varphi_n(x)+\varepsilon}-\sqrt
{\varepsilon})}{\sqrt{\lambda_i\varphi_n(x)+\varepsilon
}}U^i\nonumber\\
&&{} +\frac{1}{2}\sum_{i\neq j}
\frac{\varphi_n(x)(\sqrt{\lambda_j\varphi_n(x)+\varepsilon}-\sqrt
{\varepsilon})}{\sqrt{\lambda_i\varphi_n(x)+\varepsilon}+\sqrt
{\lambda_j\varphi_n(x)+\varepsilon}}U^i\\
&&{} +\frac{1}{2}\sum_{i,k,l}\frac{\lambda_i}{2\sqrt{\lambda_i\varphi
_n(x)+\varepsilon}}
\langle\nabla\varphi_n(x), \sigma_{\varepsilon,n}^{kl}\rangle
Z^{ikl},\nonumber
\end{eqnarray}
where $U^i_{mn}=((\Sigma^{\top}\Sigma)O)_{mi}O_{ni}+((\Sigma^{\top
}\Sigma)O)_{ni}O_{mi}$ and
$Z^{ikl}_{mn}=O_{mi}O_{ki}\Sigma_{ln}+O_{ni}O_{ki}\Sigma_{lm}$.

Furthermore, if
%
%e5.16 ###
%
\begin{equation}\label{cond: drift eps}
b \succeq(d-1)\Sigma^{\top}\Sigma,
\end{equation}
then
%
%e5.17 ###
%
\begin{equation}\label{eq: inward cond}
\Biggl\langle b-\frac{1}{2}\sum_{k,l=1}^d
D\sigma_{\varepsilon,n}^{kl}(x)\sigma_{\varepsilon
,n}^{kl}(x),u\Biggr\rangle
\geq0
\end{equation}
for all $x \in\partial S_d^+$ and for all $u \in N_{S_d^+}(x)$.
\end{lemma}
\begin{pf}
Let us denote
\[
C^{\varepsilon,n}(x) =\frac{1}{2}\sum_{k,l=1}^d
D\sigma_{\varepsilon,n}^{kl}(x)\sigma_{\varepsilon,n}^{kl}(x),
\]
and
notice that
\begin{eqnarray*}
C^{\varepsilon,n}(x)&=&\frac{1}{2}\sum_{k,l}\biggl(\frac
{d}{dt}s_{\varepsilon,n}\bigl(x+t\sigma_{\varepsilon
,n}^{kl}(x)\bigr)\big|_{t=0}M^{kl}\Sigma
\\
&&\hspace*{29pt}{} +\Sigma^{\top}(M^{kl})^\top\,\frac{d}{dt}s_{\varepsilon
,n}\bigl(x+t\sigma_{\varepsilon,n}^{kl}(x)\bigr)\big|_{t=0}\biggr).
\end{eqnarray*}
We now use the following formula from \cite{horn}, Theorem 6.6.30:
\[
\frac{d}{dt} f(V(t))=O(t)\biggl(\sum_{i,j} \Delta f(\lambda_i(t),
\lambda_j(t))M^{ii}[O(t)^\top V'(t)O(t)]M^{jj}\biggr) O(t)^\top,
\]
where $V(t)=O(t)\diag(\lambda_1(t),\ldots,\lambda_d(t))O(t)^\top$ is
a family of symmetric matrices and $\Delta
f(u,v)=\frac{(f(u)-f(v))}{(u-v)}$ for\vspace*{1pt} $u\neq v$ and $\Delta
f(u,u)=f'(u)$. This holds true if $V(\cdot)$ is continuously
differentiable for $t \in(a,b)$ and $f(\cdot)$ is continuously
differentiable on an open real interval which contains all
eigenvalues of $V(t)$ for all $t\in(a,b)$.

We now apply this formula to our case, where $f(t)=\sqrt{t}$ and
\[
V(t)=\varphi_n\bigl(x+t\sigma^{kl}_{\varepsilon,n}(x)\bigr)\bigl(x+t\sigma
^{kl}_{\varepsilon,n}(x)\bigr)+\varepsilon
I_d.
\]
Since we take the derivative at $t=0$, we only have to consider
\[
V(0)=O\bigl(\varphi_n(x)\Lambda+\varepsilon I_d\bigr)O^\top,
\]
where $O$ is the orthogonal matrix diagonalizing $x$ and
\[
V'(0)=\langle\nabla\varphi_n(x),\sigma^{kl}_{\varepsilon
,n}(x)\rangle
x+\varphi_{n}(x)\sigma^{kl}_{\varepsilon,n}(x).
\]
Note that we do not have an explicit contribution of
$\eta_{\varepsilon}$ which is part of the definition of
$s_{\varepsilon,n}$, since $\eta_{\varepsilon}(S_d^+)=1$ and
$\nabla\eta_{\varepsilon}(S_d^+)=0$. Some lines of calculations
then yield (\ref{eq: Strato corr}).

Let us now verify (\ref{eq: inward cond}). Take an arbitrary
$x=O\Lambda O^\top\in\partial S_d^+$ and assume first that it has
rank $d-1$, that is, $\lambda_d=0$ and all other eigenvalues are
strictly positive. By Lemma \ref{lem: zero divisor} and (\ref{eq:
normal cone Sd+}), the elements of $N_{S_d^+}(x)$ can then be written
as $u=OKO^\top$, where $K=\diag(0,\ldots,0,k)$ with $k \geq0$.
Thus, (\ref{eq: inward cond}) now reads
\[
\langle b-C^{\varepsilon,n}(x),OKO^\top\rangle=k[O^\top b O-O^\top
C^{\varepsilon,n}(x)O]_{dd}.
\]
As $[O^\top U^iO]_{dd}=2 \delta_{id}(O^{\top}\Sigma^{\top}\Sigma
O)_{id}$ and $O^{\top}Z^{ikl}O=2 \delta_{id}O_{ki}(\Sigma O)_{ld}$,
we have
\[
[O^\top C^{\varepsilon,n}(x)O]_{dd}=\sum_{j\neq
d}\frac{\varphi_n(x)(\sqrt{\lambda_j\varphi_n(x)
+\varepsilon}-\sqrt{\varepsilon})}{\sqrt{\lambda_j\varphi
_n(x)+\varepsilon}+\sqrt{\varepsilon}}[O^{\top}\Sigma^{\top
}\Sigma
O]_{dd}.
\]
Since $ \sum_{j\neq d}\frac{\varphi_n(x)(\sqrt{\lambda_j\varphi_n(x)
+\varepsilon}-\sqrt{\varepsilon})}{\sqrt{\lambda_j\varphi
_n(x)+\varepsilon}+\sqrt{\varepsilon}}\leq
d-1$, we obtain by condition (\ref{cond: drift eps})
\[
[O^\top b O-O^{\top}C^{\varepsilon,n}(x)O]_{dd}\geq
\bigl[O^{\top}\bigl(b-(d-1)\Sigma^{\top}\Sigma\bigr)O\bigr]_{dd} \geq0,
\]
which proves (\ref{eq: inward cond}) for $x \in\partial S_d^+$ with
$\rk=d-1$. In the general case, we can proceed similarly. For $x \in
\partial S_d^+$ with $\rk=r\leq d-1$, the elements of $N_{S_d^+}(x)$
are given by $u=OKO^\top$, where
\[
K=\pmatrix{
0 & 0\cr
0 & k}
\]
with $k \in S_{d-r}^+$. This follows again from Lemma \ref{lem: zero
divisor} and (\ref{eq: normal cone Sd+}). Now, (\ref{eq: inward
cond}) can be written as
\begin{eqnarray*}
&&\langle b-C^{\varepsilon,n}(x),O K O^\top\rangle\\
&&\qquad =\biggl\langle O^{\top} \biggl(b-\sum_{j\leq r}\frac{\varphi_n(x)(\sqrt
{\lambda_j\varphi_n(x)+\varepsilon}-\sqrt{\varepsilon})}{\sqrt
{\lambda_j\varphi_n(x)+\varepsilon}+\sqrt{\varepsilon}}\Sigma
^{\top}\Sigma\biggr)O,K\biggr\rangle\\
&&\qquad \geq\langle O^{\top} (b-r\Sigma^{\top}\Sigma)O,K \rangle\geq
0,
\end{eqnarray*}
which proves the assertion.
\end{pf}

Combining Lemmas \ref{lemAconvA} and \ref{lemprop: martingale
problem approx}, we obtain the announced existence result for the
martingale problem for $\Acal$.
\begin{lemma}\label{lemth: stoch exist}
For every $x\in S_d^+$, there exists an $S_d^+\cup\{\Delta\}$-valued
c\`adl\`ag solution $X$ to the martingale problem for $\mathcal{A}$
with $X_0=x$. That is,
\[
f(X_t)-\int_0^t \mathcal{A}f(X_s) \,ds
\]
is a martingale, for all $f\in\Scal_+$.
\end{lemma}
\begin{pf}
By Lemma \ref{lemprop: martingale problem approx}, there exists a
solution $X^{\varepsilon, \delta, n}$ to the martingale problem for
$\mathcal{A}^{\varepsilon, \delta, n}$ with sample paths in
$\mathbb{D}(S_d^+)$ (the space of $S_d^+$-valued c\`adl\`ag paths),
and hence also in $\mathbb{D}(S_d^+ \cup\{\Delta\})$.
We now claim that $(X^{\varepsilon, \delta, n})$ is
relatively compact considered as a sequence of processes with sample
paths in $\mathbb{D}(S_d^+ \cup\{\Delta\})$.\footnote{This means that
the family of probability distributions associated to $(X^{\varepsilon
, \delta, n})$
is relatively compact, that is, the closure of $(\mathbb
{P}^{\varepsilon, \delta, n})$
in $\mathcal{P}(\mathbb{D}(S_d^+ \cup\{\Delta\}))$ is compact. Here,
$\mathcal{P}(\mathbb{D}(S_d^+ \cup\{\Delta\}))$ denotes the family of
probability distributions on $\mathbb{D}(S_d^+ \cup\{\Delta\})$ and
$\mathbb{P}^{\varepsilon, \delta, n}$ the distribution of
$X^{\varepsilon, \delta, n}$.} For the proof of this
assertion, we shall make use of Theorems 9.1 and 9.4 in Chapter 3
of \cite{ethier}. In order to meet the assumption of
\cite{ethier}, Chapter 3, Theorem 9.4, we take $C_c^{\infty}(S_d^+)$
as subalgebra of
$C_b(S_d^+)$. Then, for every $T>0$ and $f \in C_c^{\infty}(S_d^+)$,
we have
\[
\sup_{\varepsilon, \delta,n} \mathbb{E}_x\Bigl[{\mathop{\esssup}_{t \in
[0,T]}}|\mathcal{A}^{\varepsilon, \delta,
n}f(X_{t}^{\varepsilon, \delta, n})|\Bigr]< \infty,
\]
since there exists a constant $C$ such that
$\|\mathcal{A}^{\varepsilon, \delta, n}f\|_{\infty}\leq C
p_{3,+}(f)< \infty$ for all $n, \varepsilon, \delta$, where
$p_{k,+}$ are the semi-norms
as defined in (\ref{eq: seminorm}) (see also the proof of Proposition
\ref{th: generator}). Thus, the requirements of
\cite{ethier}, Chapter 3, Theorem 9.4, are satisfied. Note that $Y$ in
the notation of \cite{ethier}, Chapter 3, Theorem 9.4, corresponds in
our case to $f(X)$ such that \cite{ethier},
Chapter 3, Condition
(9.17), is automatically fulfilled. It then follows by the
conclusion of \cite{ethier}, Chapter 3, Theorem 9.4, that
$(f(X_{t}^{\varepsilon, \delta, n}))$ is relatively compact [as
family of processes with sample paths in $\mathbb{D}(\re)$] for each
$f \in
C_c^{\infty}(S_d^+)$. Furthermore, since we consider $S_d^+\cup
\{\Delta\}$, the compact containment condition is always satisfied,
that is, for every $\eta>0$ and $T>0$, there exists a compact set
$\Gamma_{\eta, T}\subset(S_d^+\cup\{\Delta\})$ for which
\[
\inf_{\varepsilon, \delta, n}\mathbb{P}_x\bigl[X_{t}^{\varepsilon,
\delta, n}\in\Gamma_{\varepsilon, T} \mbox{ for $t \in
[0,T]$}\bigr] \geq1-\eta
\]
holds true. By \cite{ethier}, Chapter 3, Theorem 9.1, and the fact that
$\{1,C_c^{\infty}(S_d^+)\}$ is dense in $C(S_d^+\cup
\{\Delta\})$, we therefore obtain that $(X^{\varepsilon, \delta,
n})$ is relatively compact in $\mathbb{D}(S_d^+ \cup\{\Delta\})$.
Thus, there exists a subsequence $(\mathbb{P}^{\varepsilon_k, \delta
_k, n_k})$
of the probability distributions associated to $(X^{\varepsilon,
\delta, n})$
which converges in the Prohorov metric to some limit probability distribution.
By \cite{ethier}, Chapter 3, Theorem 3.1, this implies weak
convergence of
$(\mathbb{P}^{\varepsilon_k, \delta_k, n_k})$ and hence the
subsequence $(X^{\varepsilon_k, \delta_k, n_k})$ converges in
distribution to some limit process $X$ in $\mathbb{D}(S_d^+ \cup
\{\Delta\})$.

Combining this with Lemma \ref{lemAconvA} and \cite{ethier},
Chapter 4, Lemma 5.1, we conclude that $X$ is a solution to the
martingale problem for $\mathcal{A}$. Hence, the lemma is proved.
\end{pf}

We can now prove the existence and uniqueness of an affine process
for any admissible parameter set.
\begin{proposition}\label{th: existence Markov}
Let $(\alpha, b,\beta^{ij}, c, \gamma, m, \mu)$ be an admissible
parameter set. Then there exists a unique affine process on $S_d^+$
with infinitesimal generator (\ref{eq: generator}), and (\ref{def:
affine process}) holds for all $(t,u) \in\re_+\times S_d^+$, where
$\phi(t,u)$ and $\psi(t,u)$ are given by (\ref{eq: F-Riccati})
and (\ref{eq: R-Riccati}).
\end{proposition}
\begin{pf}
Suppose first that $c=0$ and $\gamma=0$. Let $x \in S_d^+$. Then
Lemma \ref{lemth: stoch exist} implies the existence of an
$S_d^+\cup\{\Delta\}$-valued c\`adl\`ag solution $X$ of the
martingale problem for $\Acal$ with $X_0=x$. We now show that $X$
is unique in distribution.

Thereto, note that by \cite{ethier}, Chapter 4, Theorem 7.1,
%
%e5.18 ###
%
\begin{equation}\label{eqmartftxt}
f(t,X_t)-\int_0^t \bigl( \Acal f(s,X_s)+\partial_s
f(s,X_s)\bigr)\,ds
\end{equation}
is a martingale for all rapidly decreasing functions $f\in
\Scal(\R_+\times S_d^+)$, similarly defined as $\Scal_+$ in
(\ref{defS+}). Now let $\phi$ and $\psi$ be the unique solutions of
the generalized Riccati differential equations (\ref{eq:
F-Riccati}) and (\ref{eq: R-Riccati}), given by Proposition
\ref{prop_ricc_sol}. Fix $t>0$, $u\in S_d^{++}$, and some
$f\in\Scal(\R_+\times S_d^+)$ such that
\[
f(s,x)=e^{-\phi(t-s,u)-\langle\psi(t-s,u),x\rangle},\qquad 0\le
s\le t, x\in S_d^+.
\]
Then
\[
\Acal f(s,x)+\partial_s f(s,x)=0,\qquad 0\le
s\le t, x\in S_d^+.
\]
In view of (\ref{eqmartftxt}), the Laplace
transform of $X_t$ at $u$ is thus given by
%
%e5.19 ###
%
\begin{equation}\label{eqexiproof}\quad
\E_x\bigl[e^{-\langle
u,X_t\rangle}\bigr]=\E_x[f(t,X_t)]=f(0,x)-0=e^{-\phi(t,u)-\langle\psi
(t,u),x\rangle}.
\end{equation}
Since $u\in S_d^{++}$ was arbitrary, we conclude that the
distribution of $X_t$ is uniquely determined for all $t>0$. From
\cite{ethier}, Chapter 4, Theorem 4.2, we infer that $X$ is a Markov
process with generator $\Acal$ on $\Scal_+$ and thus unique in law
as solution of the martingale problem for $\Acal$. Moreover, by
(\ref{eqexiproof}), $X$ is stochastically continuous and affine.
Thus, the proposition is proved under the premise that $c=0$ and
$\gamma=0$.

For general parameters $c$ and $\gamma$, we employ a Feynman--Kac
argument. Denote by $\Bcal$ and $(Q_t)$ the affine generator and
corresponding Feller semigroup associated with
$(\alpha,b,\beta^{ij},c=0,\gamma=0,m,\mu)$ from the first part of
the proof, respectively. Since $x\mapsto c+\langle\gamma,x\rangle$
is nonnegative on $S_d^+$, it follows along the lines of
\cite{dfs}, Proposition 11.1, that
\[
P_t f(x)=\mathbb E_x\bigl[e^{-\int_0^t c+\langle\gamma,X_s\rangle
\,ds}f(X_t)\bigr]
\]
defines a Feller semigroup $(P_t)$ on $C_0(S_d^+)$ with
infinitesimal generator $\mathcal Af(x)=\mathcal B
f(x)-(c+\langle\gamma,x\rangle) f(x)$ for $f\in\mathcal{S}_+$,
which is the desired solution.
\end{pf}

%s5.3 ###
\subsection{An alternative existence proof for jump processes}\label{subsec:
alternative existence}

For affine processes without diffusion component (i.e., the
admissible parameter $\alpha$ vanishes), the existence question can
be handled entirely as in the case of affine processes on $\mathbb
R_+^m\times\mathbb R^n$ \cite{dfs}, Section 7. In this section, we
elaborate an alternative existence proof in this specific case, by
following the lines of \cite{dfs}. Note that the OU-type processes
driven by matrix L\'evy subordinators \cite{barndorffstelzer} are
contained in the class of pure jump processes of this section.

We call a function $f\dvtx S_d^+\rightarrow\mathbb R$ of
L\'evy--Khintchine form on $S_d^+$, if
\[
f(u)=\langle b_0,u\rangle-\int_{S_d^+\setminus\{0\}}\bigl(e^{-\langle u,
\xi\rangle}-1\bigr)m_0(d\xi),
\]
where $b_0 \in S_d^+$ and $m_0$ is a Borel measure supported on
$S_d^+$ such that
\[
\int_{S_d^+\setminus\{0\}}(\|\xi\|\wedge1)m_0(d\xi) <
\infty.
\]

Once again, we recall that a distribution on $S_d^+$ is infinitely
divisible if and only if its Laplace transform takes the form
$e^{-f(u)}$, where $f$ is of the above form (see also Step 1 in the
proof of Proposition \ref{th: necessary admissibility}).

Similarly to \cite{dfs}, we introduce the sets
\begin{eqnarray*}
\mathcal
C&:=&\{f+c \mid f\dvtx S_d^+\rightarrow\mathbb R\mbox{{ is of
L\'evy--Khintchine form on }} S_d^+, c\in\mathbb R_+\},\\
\mathcal C_{S}&:=&\{\psi\mid u\mapsto\langle\psi(u),
x\rangle\in\mathcal C \mbox{{ for all }} x\in S_d^+ \}.
\end{eqnarray*}
The following technical statement can be obtained easily by
mimicking the proofs of the corresponding statements in
\cite{dfs}, Proposition 7.2 and Lemma 7.5:
\begin{lemma}\label{lemma: three dfs statements}
We have:
\begin{longlist}[(iii)]
\item[(i)]$\mathcal C$, $\mathcal C_{S}$ are convex cones in
$C(S_d^+)$.
\item[(ii)] $\phi\in\mathcal C$, $\psi\in\mathcal C_{S}$ imply $\phi
(\psi)\in\mathcal C$.
\item[(iii)] $\psi,\psi_1\in\mathcal C_{S}$ imply $\psi_1(\psi)\in
\mathcal C_{S}$.
\item[(iv)] If $\phi_k\in\mathcal C$ converges to a continuous function
$\phi$ on $S_d^+$, then $\phi\in\mathcal C$. A similar statement
holds for
sequences in $\mathcal C_{S}$.
\item[(v)] Let $(\alpha=0, b,\beta^{ij}, c, \gamma, m,
\mu)$ be an
admissible parameter set. Then $R^\delta\rightarrow R$ locally
uniformly as
$\delta\rightarrow0$, where $R^\delta$ corresponds to the
admissible parameter set $(\alpha=0, b,\beta^{ij}, c, \gamma, m, \mu
1_{\{\|\xi\|\geq\delta\}})$. (Note that there is one fixed
truncation function.)
\end{longlist}
\end{lemma}
\begin{proposition}\label{prop: C CS semiflow}
Let $(\alpha=0, b,\beta^{ij}, c, \gamma, m, \mu)$ be an admissible
parameter set. Then for all $t\geq0$, the solutions $(\phi(t,\cdot
),\psi(t,\cdot))$ of
(\ref{eq: F-Riccati}) and (\ref{eq: R-Riccati}) lie in $(\mathcal
C,\mathcal C_S)$.
\end{proposition}
\begin{pf}
Suppose first that\footnote{According to our conjecture in Section~\ref{sec2.1.4}, this would
already cover all possible jump measures if $d \geq 2$.}
%
%e5.20 ###
%
\begin{equation}\label{finite var}
\int_{S_d^+ \setminus\{0\}} \frac{\mu_{ij}(d\xi)}{\|\xi\| \wedge
1} < \infty
\end{equation}
for all $i\le j$. Then equation
(\ref{eq: R-Riccati}) is equivalent to the integral equation
%
%e5.21 ###
%
\begin{equation}\label{eq: var of const}
\psi(t,u)=e^{\widetilde B^{\top} t}(u)+\int_0^te^{\widetilde B^{\top}
(t-s)}\widetilde R(\psi(s,u))\,ds,
\end{equation}
where $R(u)=\widetilde R(u)+\widetilde B^\top(u)$ and
$\widetilde B^\top\in\mathcal L(S_d)$ is given by
\[
\widetilde B^\top(u):=B^\top(u)-\int
_{S_d^+\setminus\{0\}}\frac{\langle
\chi(\xi),u\rangle}{\|\xi\|^2\wedge1}\mu(d\xi).
\]
Here, $e^{\widetilde{B}^{\top} t}(u)$ is the notation for the semi-group
induced by $\partial_t x(t,u)=\widetilde B^\top(x(t$, $u))$,
$x(0,u)=u$. Hence, the variation of constants formula yields
(\ref{eq: var of const}).

Due to admissibility condition (\ref{eq: betaij}), we have that
$\widetilde B^\top$ is a linear drift which is ``inward pointing'' at
the boundary of $S_d^+$, which is equivalent to $e^{\widetilde
{B}^{\top}
t}$ being a positive semi-group, that is, $e^{\widetilde{B}^{\top}
t}$ maps
$S_d^+$ into $S_d^+$. Therefore, $e^{\widetilde{B}^{\top} t} \in
\mathcal{C}_S$
and since $\widetilde R(u)$ is given by
\[
\widetilde R(u)=\gamma-\int_{S_d^+\setminus\{0\}}\frac{(e^{-\langle
u,\xi\rangle}-1)}{\|\xi\|^2 \wedge1}\mu(d\xi)
\]
with $\mu$ satisfying (\ref{finite var}), we also have
%
%e5.22 ###
%
\begin{equation}\label{tilRinCS}
\widetilde R \in\mathcal C_S.
\end{equation}

Using Picard's iteration and Lemma \ref{lemma: three dfs statements},
it follows that the sequence $\psi^{(k)}$ defined as
\begin{eqnarray*}
\psi^{(0)}(t,u)&:=&u,\\
\psi^{(k+1)}(t,u)&:=&e^{\widetilde B^{\top} t}(u)+\int
_0^te^{\widetilde B^{\top}
(t-s)}\widetilde R\bigl(\psi^{(k)}(s,u)\bigr)\,ds,
\end{eqnarray*}
lies in $\mathcal C_S$, for each $t\geq0$, hence so does its limit
$\psi(t,\cdot)$. Since $F \in\mathcal{C}$, we have again by Lemma
\ref{lemma: three dfs statements}
$\phi(t,\cdot)=\int_0^t F(\psi(s,\cdot))\,ds\in\mathcal C$.

By an application of Lemma \ref{lemma: three dfs statements}(v),
the general case is then reduced to the former, since
$R^\delta$ clearly satisfies (\ref{finite var}).
\end{pf}

We are prepared to provide an alternative proof of
Proposition \ref{th: existence Markov} under the additional
assumption $\alpha=0$: by Proposition \ref{prop: C CS semiflow},
$(\phi(t,\cdot),\psi(t,\cdot))$ lie in $(\mathcal C,\mathcal C_S)$.
Hence for all $t\geq0$, $x\in S_d^+$, there exists an infinitely
divisible sub-stochastic
kernel $p_t(x,{d}\xi)$ with Laplace-transform $e^{-\phi(t,u)-\langle
\psi(t,u),x\rangle}$. The Chapman--Kolmogorov equations hold in view
of properties (\ref{eq: flow phi}) and (\ref{eq: flow psi}). Whence,
Proposition \ref{th: existence Markov} follows.
\begin{remark}\label{remapproxinfdiv}
We note that the proof of statement (v) in Lemma \ref
{lemma: three dfs statements} is much easier
than the one of \cite{dfs}, Lemma 7.5, because $\alpha=0$.
However, for $\alpha \neq 0$ and $d \geq 2$, $R$ cannot be locally
uniformly approximated by functions $R^{\delta}$ of a pure jump type
such as in Lemma~\ref{lemma: three dfs statements}.
 Indeed, otherwise one could infer as above the
existence of an affine process which is
infinitely decomposable and has nonvanishing diffusion component.
This is in contradiction with Proposition \ref{th: char infdec
infdiv} and in the case of pure diffusions it contradicts
Example \ref{counterex: Bru}.
\end{remark}

%s6 ###
\section{Proof of the main results}\label{secproofs}

%s6.1 ###
\subsection{\texorpdfstring{Proof of Theorem \protect\ref{th: main theorem}}{Proof of Theorem 2.4.}}
The first part is a summary of Propositions \ref{th: regularity},
\ref{th: necessary admissibility}, \ref{th: generator}
and \ref{th: constant drift}. The second part follows from
Proposition \ref{th: existence Markov}.

%s6.2 ###
\subsection{\texorpdfstring{Proof of Theorem \protect\ref{thmsemim}}{Proof of Theorem 2.6.}}
Let $X$ be a conservative affine process. It is shown in
Proposition \ref{th: generator} that $\{e^{-\langle u, \cdot
\rangle} \mid u \in S_d^{++}\} \subset D(\mathcal{A})$. Hence,
\[
e^{-\langle u, X_t \rangle}-e^{-\langle u, x \rangle}-\int_0^t
\mathcal{A}e^{-\langle u, X_s \rangle} \,ds
\]
is a $(\mathcal{\widetilde{F}}_t,\mathbb{P}_x)$-martingale with
$\mathcal{\widetilde{F}}_t$ defined in (\ref{eq: filt}).
From \cite{jacod}, Theorem II.2.42, combined with (\ref{Aexpux}) and
Remark \ref{remcons}, it then follows that $X$ is a semimartingale
with characteristics (\ref{eq: char c})--(\ref{eq: char nu}). The
canonical semimartingale representation (\cite{jacod}, Theorem
II.2.34) of $X$ is thus given by
\begin{eqnarray*}
X_t&=&x+B_t+X_t^c+\int_0^t\int_{{S_d^+}\setminus\{0\}}\chi(\xi)\bigl(\mu
^X(ds,d\xi)-\nu(ds,d\xi)\bigr)\\
&&{} +\int_0^t\int_{{S_d^+}\setminus\{0\}}\bigl(\xi-\chi(\xi)\bigr)\mu
^X(ds,d\xi),
\end{eqnarray*}
where $X^c$ denotes the continuous martingale part, and $\mu^X$ the
random measure associated with the jumps of $X$. In order to
establish representation (\ref{eq: BM rep}), we find it convenient
to consider the vectorization, $\vecnew(X^c) \in\re^{d^2}$, of
$X^c$. The aim is now to find a $d^2$-dimensional Brownian motion
$\widetilde{W}$ on a possibly enlarged probability space and a $d^2
\times d^2$-matrix-valued function $\sigma$ such that
%
%e6.1 ###
%
\begin{equation}\label{eq: cont mart rep}
\vecnew(X^c_t)=\int_0^t \sigma(X_s)\,d\widetilde{W}_s.
\end{equation}
Thus, $\sigma$ has to fulfill
%
%e6.2 ###
%
\begin{eqnarray}\label{eq: quad var brow}
d\langle X^c_{ij}, X^c_{kl}
\rangle_t&=&X_{t,ik}\alpha_{jl}+X_{t,il}\alpha_{jk}+X_{t,jk}\alpha
_{il}+X_{t,jl}\alpha_{ik}\nonumber\\[-8pt]\\[-8pt]
&=&(\sigma(X_t)\sigma^{\top}(X_t))_{ijkl}.\nonumber
\end{eqnarray}

As suggested by (\ref{sigmakldef}), we define the entries of the
$d^2\times d^2$-matrix $\sigma(x)$ in terms of $\sigma^{kl}(x)$
given in (\ref{sigmakldef}) by
%
%e6.3 ###
%
\begin{equation}\label{eq: sigma matrix}
\sigma_{ijkl}(x)=\sigma^{kl}_{ij}(x)=\sqrt{x}_{ik}\Sigma
_{lj}+\Sigma^{\top}_{il}\sqrt{x}_{jk}.
\end{equation}
Note that the $(kl)$th column of $\sigma(x)$ is just the
vectorization of the matrix $\sigma^{kl}(x)$. We thus obtain
$A_{ijkl}(x)= (\sigma(x)\sigma^{\top}(x))_{ijkl}$. Hence, $\sigma(x)$
satisfies (\ref{eq: quad var brow}). Analogous to the proof
of \cite{rogers}, Theorem 20.1, we can now build a $d^2$-dimensional
Brownian motion $\widetilde{W}$ on an enlargement of the probability
space such that (\ref{eq: cont mart rep}) holds true. As the
$(ij)$th entry of $X^c$ is given by
\begin{eqnarray*}
X^c_{t,ij}&=&\vecnew(X^c_{t})_{ij}=\int_0^t\sum_{k,l=1}^d\sigma
_{ijkl}(X_s)\,d\widetilde{W}_{s,kl}\\
&=&\int_0^t\bigl(\sqrt{X_s}\,dW_s\Sigma+\Sigma^{\top}\,dW_s^{\top}\sqrt{X_s}\bigr)_{ij},
\end{eqnarray*}
where $W$ is the $d \times d$-matrix Brownian motion satisfying
$\vecnew(W)=\widetilde{W}$, we obtain the desired representation.

%s6.3 ###
\subsection{\texorpdfstring{Proof of Theorem \protect\ref{th: char infdec infdiv}}{Proof of Theorem 2.9.}}

We first prove some technical lemmas.
\begin{lemma}\label{lem: Cauchy Functional equation}
Let $g\dvtx S_d^+\rightarrow\mathbb R$ be an additive function, that is,
$g$ satisfies Cauchy's functional equation
%
%e6.4 ###
%
\begin{equation}
g(x+y)=g(x)+g(y),\qquad x,y\in S_d^+.
\end{equation}
Then $g$ can be extended to an additive function $f\dvtx
S_d\rightarrow\mathbb R$. Moreover, if $g$ is measurable on $S_d^+$
then $f$ is measurable on $S_d$. In that case, $f$ is a continuous
linear functional, that is, $f(x)=\langle c,x\rangle$ for some $c\in
S_d$.
\end{lemma}
\begin{pf}
The first part follows from Lemma \ref{lem: linear map cone}.

Concerning measurability, let $E\in\mathcal B(\mathbb R)$, a Borel
measurable set. Then we have by the additivity of $f$,
\begin{eqnarray*}
f^{-1}(E)&=&\bigcup_{n=1}^\infty B_n=\bigcup_{n=1}^\infty\{x+n
I_d \mid x\in S_d, f(x)\in E, \|x\|\leq
n\}-nI_d\\&=&\bigcup_{n=1}^\infty\{y\in S_d \mid f(y)\in
E+f(nI_d), \|y-nI_d\|\leq n\}-nI_d\\&=&\bigcup_{n=1}^\infty\{y\in
S_d^+ \mid g(y)\in E+g(nI_d), \|y-nI_d\|\leq n\}-n I_d,
\end{eqnarray*}
which is again a measurable set, in view of the measurability of $g$
on $S_d^+$.

For $x\in S_d$ we write $x=(x_i)_{i}$, where $1\leq i\leq
\frac{d(d+1)}{2}$. We introduce the additive functions $f_i\dvtx \mathbb
R\rightarrow\mathbb R$ via $f_i(x_i)=f(0,\ldots,0,x_i,0,\ldots,0)$. By
the just proved measurability of $f$, we infer that all $f_i$ are
measurable functions on $\mathbb R$. By \cite{aczel}, Chapter 2,
Theorem 8, any additive measurable function on the real line is a
continuous linear functional. Hence for each $i$, we infer the
existence of $c_i\in\mathbb R$ such that $f_i(x_i)=c_i x_i$ holds.
Since $f(x)=\sum_{i} f_i(x_i)$, it follows that $f(x)=\langle
c,x\rangle$ for some $c\in S_d$.
\end{pf}

Also, we consider Cauchy's exponential equation for $h\dvtx
S_d^+\rightarrow\mathbb R_+$, that is,\vspace*{2pt}
%
%e6.5 ###
%
\begin{equation}\label{eq: multiplicative functionals}
h(x+y)=h(x)h(y),\qquad x,y\in S_d^+.
\end{equation}
\begin{lemma}\label{lem: Cauchy Exp Eq.}
Suppose $h\dvtx S_d^+\rightarrow\mathbb R_+$ is measurable, strictly
positive, and satisfies (\ref{eq: multiplicative functionals}). Then
$h(x)=e^{-\langle c,x\rangle}$, for some $c\in S_d$. If $h\leq1$,
then $c\in S_d^+$.
\end{lemma}
\begin{pf}
Since $h$ is strictly positive, its logarithm yields the well
defined function $g\dvtx S_d^+\rightarrow\mathbb R$, $g(x):=\log
h(x)$. Clearly $g$ is additive, hence by the first part of
Lemma \ref{lem: Cauchy Functional equation}, there exists a unique
additive extension $f\dvtx S_d\rightarrow\mathbb R$. Also, $f$ is
measurable on $S_d^+$, hence by the second assertion of
Lemma \ref{lem: Cauchy Functional equation} we have $f(x)=-\langle
c,x\rangle$, for some $c\in S_d$. The last statement follows from
the monotonicity of the exponential and the self duality of $S_d^+$.
\end{pf}
\begin{remark}\label{rem:
missinglink} The assumption of strict positivity of $h$ in the
preceding lemma is essential. Otherwise, there exist solutions $h$
which are not of the asserted form.
\end{remark}

Lemma \ref{lem: Cauchy Exp Eq.} is the main ingredient of the proof
of the following characterization concerning $k$-fold convolutions
of Markov processes.
\begin{lemma}\label{lema: dfslemma10.3}
Let $(\mathbb P^{(i)}_x)_{x\in S_d^+}\in\mathcal P$
($i=0,1,\ldots,k$). Then\vspace*{2pt}
%
%e6.6 ###
%
\begin{equation}\label{eq: kfoldconv}\quad
\mathbb P^{(1)}_{x^{(1)}}*\cdots*\mathbb P^{(k)}_{x^{(k)}}=\mathbb
P^{(0)}_x\qquad \forall x^{(i)}\in S_d^+,\qquad x=x^{(1)}+\cdots+x^{(k)},
\end{equation}
if and only if for all $\mathbf{t}= (t_1,\ldots, t_N)\in\mathbb R_+^N$
and $\mathbf{u}=(u^{(1)},\ldots, u^{(N)})\in(S_d^+)^N$, \mbox{$N\in\mathbb
N_0$}, there exists $0 < \rho^{(i)}(\mathbf{t},\mathbf{u})\leq1$ and
$\psi(\mathbf{t},\mathbf{u})\in S_d^+$ such that
$\prod_{i=1}^k\rho^{(i)}(\mathbf{t}$, $\mathbf{u})=\rho^{(0)}(\mathbf{t},\mathbf{u})$ and\vspace*{2pt}
%
%e6.7 ###
%
\begin{eqnarray}\label{eq: char kfoldconv}
&&\mathbb E_x^{(j)}\bigl[e^{-\sum_{i=1}^N\langle
u^{(i)},X_{t_i}\rangle}\bigr]=\rho^{(j)}(\mathbf{t},\mathbf{u})
e^{-\langle\psi(\mathbf{t},\mathbf{u}),x\rangle}\nonumber\\[-8pt]\\[-8pt]
&&\eqntext{\forall x\in
S_d^+, j=0,1,\ldots,k.}
\end{eqnarray}
\end{lemma}
\begin{pf}
We proceed similarly as in the proof of \cite{dfs}, Lemma 10.3. Fix
$k>1$, $N>1$, $\mathbf{t}, \mathbf{u}$ and set
\[
g^{(j)}(x):=\mathbb E_x^{(j)}\bigl[e^{-\sum_{i=1}^N\langle
u^{(i)},X_{t_i}\rangle}\bigr].
\]
By the definition of the convolution, (\ref{eq: kfoldconv}) is
equivalent to the following:
%
%e6.8 ###
%
\begin{eqnarray}\label{eq: equivalent to kfoldconv}
g^{(1)}\bigl(x^{(1)}\bigr)\cdot\cdots\cdot g^{(k)}\bigl(x^{(k)}\bigr)=g^{(0)}(x)\nonumber\\[-8pt]\\[-8pt]
&&\eqntext{\forall
x^{(i)}\in S_d^+,\qquad x=x^{(1)}+\cdots+x^{(k)}.}
\end{eqnarray}
Hence, the implication (\ref{eq: char kfoldconv})${}\Rightarrow{}$(\ref{eq: equivalent to kfoldconv}) is obvious. For the converse
direction, we observe that $g^{(i)}$ are strictly positive on all of
$S_d^+$. Thus, by (\ref{eq: equivalent to kfoldconv}) we have
\[
g:=g^{(1)}/g^{(1)}(0)=\cdots=g^{(k)}/g^{(k)}(0)=g^{(0)}/g^{(0)}(0)
\]
and $g$ is a measurable, strictly positive function on $S_d^+$
satisfying (\ref{eq: multiplicative functionals}). Hence, an
application\vspace*{1pt} of Lemma \ref{lem: Cauchy Exp Eq.} yields the validity
of (\ref{eq: char kfoldconv}), where $\rho^{(i)}(\mathbf{t},\mathbf{u})=g^{(i)}(0)$.
By the definition of $g^{(i)}$, it follows that $0 <\rho^{(i)}(\mathbf{t},\mathbf{u})\leq1$ and
\mbox{$\psi(\mathbf{t},\mathbf{u}) \in S_d^+$}.
\end{pf}

We are prepared to prove Theorem \ref{th: char infdec infdiv}:

(i)${}\Rightarrow{}$(ii): due to
Lemma \ref{lema: dfslemma10.3}, infinite decomposability implies that
$X$ is affine. Also, by the definition of infinite decomposability
and by Lemma \ref{lema: dfslemma10.3} we have that the $k$th root
$(\mathbb P^{(k)}_x)$ for each $k\geq1$ is an affine process with
state space $S_d^+$ with exponents $\psi(t,u)$ and $\phi(t,u)/k$.
This implies that $(\mathbb P^{(k)}_x)_{x\in S_d^+}$ has admissible
parameters $(\alpha, b/k, \beta^{ij}, c/k, \gamma, m/k, \mu)$.
Hence, the admissibility condition proved in Proposition \ref{th:
constant drift} implies $b/k\succeq(d-1)\alpha\succeq0 $, for each
$k$, which is impossible, unless $\alpha=0$ or $d=1$.

(ii)${}\Rightarrow{}$(iii): follows from
Proposition \ref{prop: C CS semiflow}, in view of the
L\'evy--Khintchine form of
$-\phi(t,\cdot)-\langle\psi(t,\cdot),x\rangle$, for each $t>0$.

(iii)${}\Rightarrow{}$(i): by definition,
every transition kernel $p_t(x,d\xi)$ of $X$ is infinitely divisible
with Laplace transform $P_te^{-\langle u,x \rangle}=e^{-\phi
(t,u)-\langle
x,\psi(t,u)\rangle}$. For each $k\geq1$, the maps
$\phi^{(k)}:=\frac{\phi}{k}$, $\psi^{(k)}:=\psi$ satisfy the
properties (\ref{eq: flow phi}) and (\ref{eq: flow psi}). Also,
infinite divisibility implies that for each $(t,x)\in\mathbb R_+\times
S_d^+$,
\[
Q^{(k)}_te^{-\langle u,x \rangle}:=e^{-\phi^{(k)}(t,u)-\langle
\psi^{(k)}(t,u),{x/k}\rangle}
\]
is the Laplace transform of a sub-stochastic measure on $S_d^+$. In
conjunction with Properties (\ref{eq: flow phi}) and (\ref{eq: flow
psi}) we may conclude that $Q^{(k)}_t$ gives rise to a Feller
semigroup on $C_0(S_d^+)$, which is affine in $y=x/k$. Hence, we have
constructed for each $k\geq1$  a $k$th root of $X$  which is
stochastic continuous  by the definition of its characteristic
exponents $\phi^{(k)},\psi^{(k)}$. Thus Theorem \ref{th: char
infdec infdiv} is proved.

\begin{appendix}

%s7 ###
\section{Existence and viability of a class of jump-diffusions}\label
{app: viability}

In this section, we study existence and viability in a nonempty
closed convex set $D \subset\re^n$ of solutions to the equation
%
%e7.1 ###
%
\begin{equation}\label{eq: jump diffusion}
X_t=x+\int_0^t b(X_s)+\int_0^t \sigma(X_s) \,dW_s + J_t,
\end{equation}
where $b(x) \in C_b(\re^n, \re^n)$, $\sigma(x) \in
C_b(\re^n,\re^{n\times m})$ are Lipschitz continuous maps, $W$ a
standard $m$-dimensional Brownian motion and $J$ a finite
activity jump process with state-dependent, absolutely continuous
compensator $K(X_t,d\xi)\,dt$. We further assume that $ x\mapsto
K(x,\re^n)$ is bounded.

We tackle this problem in three steps. First, we derive some
regularity and existence results for diffusion SDEs. These results
are not in the standard literature, we thus provide full proofs.
Second, we prove existence of a c\`adl\`ag solution $X$ for
(\ref{eq: jump diffusion}). Finally, we provide sufficient
conditions for $X$ to be $D$-valued.

%s7.1 ###
\subsection{Diffusion stochastic differential equations}

Let $(\Omega,\mathcal{F},(\mathcal{F}_t),\mathbb{P})$ be a filtered
probability space satisfying the usual conditions and carrying an
$m$-dimensional standard Brownian motion $W$. We consider the
following diffusion SDE:
%
%e7.2 ###
%
\begin{equation}\label{sde}
X_t=x+ \int_{0}^t b(X_s)1_{\{\theta\le s\}} \,ds+\int_{0}^t
\sigma(X_s)1_{\{\theta\le s\}} \,dW_s,
\end{equation}
where $(\theta,x)\in[0,\infty]\times\re^n$ and $b$ and $\sigma$ are
as above. Recall that $X$ is a solution of (\ref{sde}) if $X$ is
continuous and (\ref{sde}) holds for all $t \geq0$ a.s. In
particular, note that this null set depends on $(\theta,x)$.
\begin{lemma}
Fix $T >0$ and let $p\geq2$. Furthermore, let $\Theta_1, \Theta_2$
be stopping times and for $i=1,2$, $U_i$,
$\mathcal{F}_{\Theta_i}$-measurable random variables. Consider the
following equations:
\begin{eqnarray*}
X_t&=&U_1+ \int_{0}^t b(X_s)1_{\{\Theta_1\le s\}} \,ds+\int_{0}^t
\sigma(X_s)1_{\{\Theta_1\le s\}} \,dW_s,\\
Y_t&=&U_2+ \int_{0}^t b(Y_s)1_{\{\Theta_2\le s\}} \,ds+\int_{0}^t
\sigma(Y_s)1_{\{\Theta_2\le s\}} \,dW_s.
\end{eqnarray*}
Then there exists a constant $C$ depending only on $p$, $T$, $n$,
the Lipschitz constants of $b$ and $\sigma$ and
$\|b\|_{\infty}, \|\sigma\|_{\infty}$ such that for $0\leq t\leq
T$,
%
%e7.3 ###
%
\begin{eqnarray}\label{eq: estimate}
&&\mathbb{E}\Bigl[\sup_{s\leq t}\|X_s-Y_s\|^p\Bigr]\nonumber\\
&&\qquad\leq C\mathbb{E}\biggl[\|U_1-U_2\|^p+|\Theta_1 \wedge t-\Theta_2
\wedge t|^{{p/2}}\\
&&\qquad\quad\hspace*{70.5pt}{}+{\int_0^t\sup_{u\leq s}}\|X_u-Y_u\|^p
\,ds\biggr].\nonumber
\end{eqnarray}
\end{lemma}
\begin{pf}
By the same arguments as in the proof of \cite{rogers}, Lemma 11.5,
we first obtain the following estimate:
\begin{eqnarray*}
&&\sup_{s\leq t}\|X_s-Y_s\|^p\\
&&\qquad\leq 3^{p-1}\biggl(\|U_1-U_2\|^p+\biggl(\int_0^t\|
b(X_s)1_{\{\Theta_1\le s\}}-b(Y_s)1_{\{\Theta_2\le s\}}\| \,ds\biggr)^p\\
&&\qquad\quad\hspace*{58.7pt}{} +\sup_{s\leq
t}\biggl\|\int_0^s\bigl(\sigma(X_u)1_{\{\Theta_1\le
u\}}-\sigma(Y_u)1_{\{\Theta_2\le u\}}\bigr) \,dW_u \biggr\|^p\biggr).
\end{eqnarray*}
Moreover,
\begin{eqnarray*}
&&\biggl(\int_0^t\bigl\|b(X_s)1_{\{\Theta_1\le s\}}-b(Y_s)1_{\{\Theta_2\le s\}
}\bigr\| \,ds\biggr)^p\\
&&\qquad\leq 2^{p-1}\biggl(\biggl(\int_{\Theta_1\wedge t}^{(\Theta_1 \vee\Theta
_2)\wedge t}\|b(X_s)\| \,ds\biggr)^p+\biggl(\int_{\Theta_2\wedge t}^{(\Theta_1
\vee\Theta_2)\wedge t}\|b(Y_s)\| \,ds\biggr)^p\\
&&\qquad\quad\hspace*{117.2pt}{} +\biggl(\int_{(\Theta_1\vee\Theta_2)\wedge t}^t\|b(X_s)-b(Y_s)\|
\,ds\biggr)^p\biggr)\\
&&\qquad\leq 2^{p-1}\biggl(K|\Theta_1\wedge t-\Theta_2 \wedge t|^p+t^{p-1}\int
_{0}^t\|b(X_s)-b(Y_s)\|^p \,ds\biggr)\\
&&\qquad\leq K\biggl(t^{{p/2}}|\Theta_1\wedge t-\Theta_2 \wedge
t|^{p/2}+\int_{0}^t\sup_{u\leq s}\|X_u-Y_u\|^p \,ds\biggr).
\end{eqnarray*}
For the stochastic integral part, we apply the
Burkholder--Davis--Gundy inequality
\begin{eqnarray*}
&&\mathbb{E}\biggl[\sup_{s\leq t}\biggl\|\int_0^s\bigl(\sigma(X_u)1_{\{\Theta_1\le
u\}}-\sigma(Y_u)1_{\{\Theta_2\le u\}}\bigr) \,dW_u\biggr\|^p\biggr]\\
&&\qquad\leq K\mathbb{E}\biggl[\biggl(\int_0^t\bigl\|\sigma(X_u)1_{\{\Theta_1\le u\}
}-\sigma(Y_u)1_{\{\Theta_2\le u\}}\bigr\|^2 \,du\biggr)^{{p/2}}\biggr]\\
&&\qquad\leq K\mathbb{E}\biggl[\biggl(\int_{\Theta_1 \wedge t}^{(\Theta_1 \vee\Theta
_2)\wedge t}\|\sigma(X_u)\|^2 \,du\biggr)^{p/2}+\biggl(\int_{\Theta_2
\wedge t}^{(\Theta_1 \vee\Theta_2)\wedge t}\|\sigma(Y_u)\|^2
\,du\biggr)^{p/2}\\
&&\qquad\quad\hspace*{125.4pt}{} +\biggl(\int_{(\Theta_1 \vee\Theta_2)\wedge t}^t\|\sigma(X_u)-\sigma
(Y_u)\|^2 \,du\biggr)^{{p/2}}\biggr]\\
&&\qquad \leq K\mathbb{E}\biggl[|\Theta_1\wedge t-\Theta_2 \wedge t|^{
{p/2}}+\int_{0}^t\|\sigma(X_s)-\sigma(Y_s)\|^p \,ds\biggr]\\
&&\qquad \leq K\mathbb{E}\biggl[|\Theta_1 \wedge t-\Theta_2 \wedge
t|^{{p/2}}+\int_{0}^t\sup_{u\leq s}\|X_u-Y_u\|^p \,ds\biggr],
\end{eqnarray*}
where $K$ always denotes a constant which varies from line to line.
The last estimate in both inequalities follows from the the
Lipschitz continuity of $b$ and $\sigma$. By assembling these
pieces, the proof is complete.
\end{pf}

Here is a fundamental existence result, which is not stated in this
general form in the standard literature. Therefore, we provide a full
proof.
\begin{theorem}\label{th: regulardiff}
There exists a function $Z\dvtx [0,\infty]\times\re^n \times
\Omega\times\re_+ \rightarrow\re^n$ with the following properties:
\begin{longlist}[(iii)]
\item[(i)] $Z(\theta,x,\omega,t)$ is continuous in $(\theta
,x,t)$ for
all $\omega$.
\item[(ii)] $Z$ is $\mathcal{B}([0,\infty]\times\re
^n)\otimes\mathcal{P}$-measurable.\footnote{Here, $\mathcal{P}$
denotes the predictable $\sigma$-field.}
\item[(iii)] $Z(\theta,x,\omega,t)$ solves (\ref{sde}) for all
$(\theta,x)$.
\item[(iv)] Let $\Theta$ be a stopping time and $U$ an
$\mathcal{F}_{\Theta}$ measurable random variable, then $X_t=Z(\Theta
,U,t)$ solves
%
%e7.4 ###
%
\begin{equation}\label{sde2}
X_t=U+ \int_{0}^t b(X_s)1_{\{\Theta\le s\}} \,ds+\int_{0}^t
\sigma(X_s)1_{\{\Theta\le s\}} \,dW_s.
\end{equation}
\end{longlist}
\end{theorem}
\begin{pf}
For every $(\theta,x)\in[0,\infty]\times\re^n$, there exists a
unique solution $X_t(\omega)=\widetilde{Z}(\theta,x,\omega,t)$
of (\ref{sde}), which is continuous in $t$.
This is a consequence of the Lipschitz continuity
of $x\mapsto b(x)1_{\{\theta\le s\}}$ and $x
\mapsto\sigma(x)1_{\{\theta\le s\}}$. Uniqueness is meant modulo
indistinguishability. From estimate (\ref{eq: estimate}), we can
deduce for $p \geq2$, $x,y \in[-T,T]^n$, $0 \leq\theta_1, \theta_2
\leq T$
and $0\leq t\leq T$,
\begin{eqnarray*}
&&\mathbb{E}\Bigl[{\sup_{s\leq t}}\|\widetilde{Z}(\theta_1,x,s)-\widetilde
{Z}(\theta_2,y,s)\|^p\Bigr]\\
&&\qquad\leq
K\biggl(\|x-y\|^{p/2}+|\theta_1-\theta_2|^{p/2}\\
&&\qquad\quad\hspace*{17.3pt}+\int_0^t
\mathbb{E}\Bigl[\sup_{u\leq
s}\|\widetilde{Z}(\theta_1,x,u)-\widetilde{Z}(\theta_2,y,u)\|^p\Bigr]
\,ds\biggr)
\end{eqnarray*}
for some constant $K$. Hence, by Gronwall's lemma,
\begin{eqnarray*}
\mathbb{E}\Bigl[{\sup_{s\leq t}}\|\widetilde{Z}(\theta_1,x,s)-\widetilde
{Z}(\theta_2,y,s)\|^p\Bigr]&\leq& Ke^{KT}(\|x-y\|^{p/2}+|\theta
_1-\theta_2|^{p/2})\\
&\leq& C\|(\theta_1,x)-(\theta_2,y)\|^{p/2}.
\end{eqnarray*}
Let now $\Dya=\{j2^{-k},j \in\mathbb{Z}, k \in\mathbb{N}\}$ be the
set of dyadic rational numbers in $\re$ and $\Dya^n=\Dya
\times\cdots\times\Dya$ the set of dyadic rational numbers in
$\re^n$. Furthermore, we define $M$ by $M=\Dya^{n+1} \cap([0,T]
\times[-T,T]^n)$. By setting $p=2n+4$, we can apply Kolmogorov's
lemma. Indeed, analogous to the proof of \cite{karatzas},
Theorem 2.8, we derive for all $(\theta_1,x), (\theta_2,y) \in M$
with $0 <\|(\theta_1,x)-(\theta_2,y)\|< h(\omega)$, where $h$ is a
positive valued random variable, and for all $\omega\in
\Omega_T^{\ast}$, where $\Omega_T^{\ast} \in\mathcal{F}$ is some
set depending on $T$ with $\mathbb{P}(\Omega_T^{\ast})=1$, the
following estimate:
%
%e7.5 ###
%
\begin{equation}\label{eq: Kolmogorov estimate}
{\sup_{s\leq t} }\|\widetilde{Z}(\theta_1,x,\omega,
s)-\widetilde{Z}(\theta_2,y,\omega, s)\|\leq
\delta\|(\theta_1,x)-(\theta_2,y)\|^{\gamma}.
\end{equation}
Here, $\gamma\in(0, \frac{1}{p})$ and $\delta$ is some positive
constant. Let
us now define $Z$: if $\omega\notin\Omega_T^{\ast}$, then
$Z(\theta,x,\omega,t)=x$ for $0 \leq t \leq T$. For $\omega\in
\Omega_T^{\ast}$ and $(\theta,x) \in M$, $Z(\theta,x, \omega,
t)=\widetilde{Z}(\theta,x,\omega,t)$ for $0 \leq t \leq T$. If
$(\theta,x) \in M^c$, we choose a sequence $(\theta_n,
x_n)_{n\in\mathbb{N}} \subseteq M$ such that $(\theta_n,
x_n)\rightarrow(\theta,x)$. By estimate (\ref{eq: Kolmogorov
estimate}), $\widetilde{Z}(\theta_n,x_n, \omega, t)$ is a Cauchy
sequence converging with respect to \mbox{${\sup_{s\leq t}}\|\cdot\|$}. We
can therefore set $Z(\theta,x, \omega, t)=\lim_{n\rightarrow
\infty}\widetilde{Z}(\theta_n,x_n,\omega,t)$. As we have uniform
convergence in $t$ and as $\widetilde{Z}$ is continuous in $t$, the
resulting process $Z$ is jointly continuous
in $(\theta,x,t)$. Furthermore, for every $(\theta,x)$, $Z$ is
indistinguishable from $\widetilde{Z}$, that is,
%
%e7.6 ###
%
\begin{equation}\label{eq: cont version}
\mathbb{P}[Z(\theta,x,t)=\widetilde{Z}(\theta,x,t) \mbox{ for all
} 0 \leq t\leq T]=1.
\end{equation}
Indeed, for $(\theta,x) \in M$, this is clear and for $(\theta,x)
\in M^c$, we have for $(\theta_n,x_n)_{n\in\mathbb{N}} \subseteq M$
with $(\theta_n,
x_n)\rightarrow(\theta,x)$
\[
\mathbb{P}\Bigl[{\sup_{s\leq
t}}\|\widetilde{Z}(\theta_n,x_n,s)-\widetilde{Z}(\theta,x,s)\|\geq
\varepsilon\Bigr]\leq
C\varepsilon^{-p}\|(\theta_n,x_n)-(\theta,x)\|^{p/2},
\]
which implies that $\widetilde{Z}(\theta_n,x_n,t)\rightarrow
\widetilde{Z}(\theta,x,t)$ in probability, uniformly in $t$. As
$\widetilde{Z}(\theta_n,x_n, t)\rightarrow Z(\theta,x, t)$ a.s., and
thus in particular in probability, it follows that
$Z(\theta,x,t)=\widetilde{Z}(\theta,x,t)$ a.s. for all $0 \leq t\leq
T$. Letting $T \rightarrow\infty$ proves assertion (i).

Statement (ii) is then a consequence of (i) and the
$\mathcal{F}_t$-measurability of $\omega\mapsto
Z(\theta,x,t,\omega)$, which is satisfied since $\mathcal{F}_0$
contains all null sets of $\mathcal{F}$.

Furthermore, property (\ref{eq: cont version}) implies that
$Z(\theta,x,t)$ is a solution of (\ref{sde}) for all $(\theta,x)$,
which yields assertion (iii).

In order to prove (iv), we proceed in two steps:

\textit{Step} 1. We first assume that $\Theta$ and $U$ take finitely
many values $\theta_1,\ldots,\theta_k\in[0,\infty]$ and
$x_1,\ldots,x_l\in\re^n$, respectively. Denote
\[
A_j=\{ \Theta=\theta_j\},\qquad B_h=\{ U=x_h\} .
\]
Then
\[
Z(\Theta,U,t)=\sum_{j,h} 1_{A_j\cap B_h} Z(\theta_j,x_h,t)
\]
does the job. Indeed, as $A_j\cap B_h$ are disjoint and $A_j\cap
B_h\in\mathcal{F}_{\theta_j}$ for all $j,h$, we have (see,
e.g., \cite{lamberton}, page 39)
\begin{eqnarray*}
&&U+\int_0^t b(Z(\Theta,U,s))1_{\{\Theta\le s\}} \,ds +\int_0^t \sigma
(Z(\Theta,U,s))1_{\{\Theta\le s\}} \,dW_s \\
&&\qquad=U+\int_0^t \sum_{j,h} 1_{A_j\cap B_h} b(Z(\theta_j,x_h,s))1_{\{
\theta_j\le
s\}} \,ds\\
&&\qquad\quad{}+\int_0^t \sum_{j,h} 1_{A_j\cap B_h} \sigma(Z(\theta_j,x_h,s))1_{\{
\theta_j\le
s\}} \,dW_s\\
&&\qquad=\sum_{j,h} 1_{A_j\cap B_h} \biggl( x_h+\int_0^t b(Z(\theta_j,x_h,s))1_{\{
\theta_j\le
s\}} \,ds\\
&&\qquad\quad\hspace*{55pt}{}+\int_0^t \sigma(Z(\theta_j,x_h,s))1_{\{\theta_j\le
s\}} \,dW_s\biggr)\\
&&\qquad=\sum_{j,h} 1_{A_j\cap B_h} Z(\theta_j,x_h,t)=Z(\Theta,U,t)
\end{eqnarray*}
for all $t\ge0$ a.s.

\textit{Step} 2. For general $\Theta$, $U$, approximate
$\Theta^{(k)}\downarrow\Theta$ by the simple stopping times
\[
\Theta^{(k)} =\cases{
j 2^{-k} ,&\quad $(j-1)2^{-k}\le\Theta<j 2^{-k}, j=1,\ldots,k 2^k$,\cr
\infty,&\quad $k\le\Theta$.}
\]
Let $U^{(l)}$ be a sequence of $\mathcal{F}_\Theta$-measurable
random variables, each $U^{(l)}$ taking finitely many values, and
$U^{(l)}\to U$ in $L^2$ (such $U^{(l)}$ obviously exists). Moreover,
$\{\Theta^{(k)}= \theta_j\}\cap\{U^{(l)}=x_h\} \in
\mathcal{F}_{\theta_j}$ for all $j,h$ (see \cite{karatzas},
Chapter 1, Problem 2.24).

By Step 1, each $Z(\Theta^{(k)},U^{(l)})$ satisfies the respective
SDE. Moreover, from estimate (\ref{eq: estimate}) and Grownwall's
lemma we deduce that for any $T>0$, there exists a constant $C$ such
that
\begin{eqnarray*}
&&\E\Bigl[ \sup_{t\le
T}\bigl\|Z\bigl(\Theta^{(k)},U^{(l)},t\bigr)-Z\bigl(\Theta^{(k')},U^{(l')},t\bigr)\bigr\|^2 \Bigr] \\
&&\qquad\le C
e^{CT}\mathbb{E}\bigl[\bigl\|U^{(k)}-U^{(k')}\bigr\|^2+\bigl|\Theta^{(k)}\wedge
T-\Theta^{(k')}\wedge T\bigr|\bigr].
\end{eqnarray*}
Hence, $Z(\Theta^{(k)},U^{(l)})$ is a Cauchy sequence and thus
converging with respect to $\E[{\sup_{t\le
T}}\|\cdot\|^2]$, for all $T>0$, to some continuous process $X$
satisfying (\ref{sde2}). On the other hand, by the continuity of
$(\theta,x) \mapsto Z(\theta,x,t)$, we know that
\[
Z\bigl(\Theta^{(k)},U^{(l)},t\bigr)\to Z(\Theta,U ,t)
\]
for all $\omega$ and $t\ge0$. Again, by continuity of $t \mapsto
Z(\Theta,U,t)$, we conclude that $Z(\Theta,U )=X$ up to
indistinguishability, which proves the claim.
\end{pf}

%s7.2 ###
\subsection{Existence of jump-diffusions}

We now provide a constructive proof for the existence of a solution
of (\ref{eq: jump diffusion}) on a specific stochastic basis which
is defined as follows:
\begin{itemize}
\item$(\Omega, \mathcal{F}, (\mathcal{F}_t)_{t\geq0})$ is a
filtered space, where $\Omega:=\Omega_1 \times\Omega_2$, $\mathcal
{F}_t:=\mathcal{G}_t\otimes\mathcal{H}_t$ and $\mathcal{F}=\mathcal
{G}\otimes\mathcal{H}$ are precisely defined below. Note that we do
not have a measure on $(\Omega, \mathcal{F})$ for the
moment. The generic sample element will be denoted by
$\omega=(\omega_1,\omega_2)\in\Omega$.
\item$(\Omega_1, \mathcal{G}, (\mathcal{G}_t)_{t\geq0}, \mathbb
{P}_1)$ is some filtered probability space satisfying the usual
conditions and carrying an $m$-dimensional standard Brownian motion
$W$. We shall consider the above diffusion SDE (\ref{sde}) on
$\Omega_1$ and thus obtain the respective solution
$Z(\theta,x,\omega_1,t)$ satisfying the regularity properties of
Theorem \ref{th: regulardiff}.
\item$(\Omega_2, \mathcal{H})$ is the canonical space for $\re
^n$-valued marked point processes
(see, e.g., \cite{jacodj}): $\Omega_2$ consists of all c\`adl\`ag,
piecewise constant functions $\omega_2\dvtx [0$, $T_{\infty}(\omega_2))
\rightarrow\re^n$ with $\omega_2(0)=0$ and
$T_{\infty}(\omega_2)=\lim_{n\rightarrow\infty} T_n(\omega_2)\leq
\infty$, where $T_n(\omega_2)$, defined by $T_0=0$ and
\[
T_n(\omega_2):=\inf\{t > T_{n-1}(\omega_2) \mid
\omega_2(t)\neq\omega_2(t-)\}\wedge\infty,\qquad n\ge1,
\]
are the successive jump times of $\omega_2$. We denote by
\[
J_t(\omega)=J_t(\omega_2)=\omega_2(t) \qquad\mbox{on } [0,
T_{\infty}(\omega_2))
\]
the canonical jump process, and let $\mathcal{H}_t = \sigma(J_s \mid
s \leq t)$ be its natural filtration with $\mathcal{H} =
\mathcal{H}_{\infty}$. Note that $T_n$ are $(\Hcal_t)$ and
$(\Fcal_t)$-stopping times if interpreted as
$T_n(\omega)=T_n(\omega_2)$.
\end{itemize}

The following statement is meant to be pointwise, referring to the
filtered measure space $(\Omega,\mathcal{F},(\mathcal{F}_t))$
without reference to a probability measure.
\begin{lemma}\label{lemma: 1}
Let $Z(\theta, x,\omega_1,t)$ be as of Theorem \ref{th:
regulardiff}. Then for an $\mathcal{F}_{T_n}$-measurable random
variable $U(\omega_1,\omega_2)$ the process
$Z(T_n(\omega_2),U(\omega_1,\omega_2),\omega_1,t)$ is:
\begin{longlist}[(ii)]
\item[(i)] continuous in $t$ for all $(\omega_1,\omega_2)$,
\item[(ii)] $\mathcal{F}_t$-adapted on $\{T_n \leq t\}$.
\end{longlist}
\end{lemma}
\begin{pf}
The first assertion is a consequence of Theorem \ref{th:
regulardiff}(i). The second one follows from the
$\mathcal{B}([0,\infty]\times\re^n)\otimes
\mathcal{P}$-measurability of $Z(\theta,x,\omega_1,t)$, as stated
in Theorem \ref{th: regulardiff}(ii), and the fact that $T_n$ and
$U$ are $\mathcal{F}_t$-measurable on $\{T_n \leq t\}$.
\end{pf}

Here is our existence result for (\ref{eq: jump diffusion}).
\begin{theorem}
There exists a c\`adl\`ag $\Fcal_t$-adapted process $X$ and a
probability measure $\Pa$ on $(\Omega,\Fcal)$ with
$\Pa|_\Gcal=\Pa_1$, such that $X$ is a solution of (\ref{eq: jump
diffusion}) on $(\Omega,\Fcal,(\Fcal_t),\Pa)$.
\end{theorem}
\begin{pf}
We follow the arguments in the proof of
\cite{filipovicoverbeck}, Theorem 5.1, which is based on
\cite{jacodj}, Theorem 3.6, and proceed in three steps.

\textit{Step} 1. We start by solving (\ref{eq: jump diffusion}) along
every path $\omega_2$. To this end, let us define recursively:
$\Delta\omega_2(0)=\Delta\omega_2(\infty)=0$, $X^{(0)}_0=x$, and for
$n\ge1$:
\[
X^{(n)}_t(\omega_1,\omega_2)=\cases{
Z\bigl(T_{n-1}(\omega_2),X^{(n-1)}_{T_{n-1}}(\omega_1,\omega_2)\cr
\hspace*{37.7pt}{}+\Delta
\omega_2(T_{n-1}),\omega_1,t\bigr),&\quad
$t\in[0,\infty)$,\cr
x_0,&\quad $t=\infty$,}
\]
where $x_0$ is any fixed point in $D \subset\mathbb{R}^n$ and $Z$
satisfies the
properties of Theorem \ref{th: regulardiff}. By Lemma \ref{lemma: 1}, every
$X^{(n)}$ is continuous in $t$ for all $(\omega_1,\omega_2)$ and
$\mathcal{F}_t$-adapted on $\{T_n \leq t\}$
since $X^{(n-1)}_{T_{n-1}}(\omega_1,\omega_2)+\Delta
\omega_2(T_{n-1})$ is $\mathcal{F}_{T_n}$-measurable. Thus, the
process
%
%e7.7 ###
%
\begin{equation}\label{Xdef}
X_t(\omega_1,\omega_2) = \sum_{n\ge1} X^{(n)}_t(\omega_1,\omega
_2) 1_{\{T_{n-1}\le t<T_{n}\}}
\end{equation}
is c\`adl\`ag $\mathcal{F}_t$-adapted and solves (\ref{eq: jump
diffusion}) on $(\Omega_1,\Gcal,(\Gcal_t),\Pa_1)$ for $t\in
[0,T_\infty(\omega_2))$ and any fixed path $\omega_2$.

\textit{Step} 2. It remains to show that there exists a probability
measure $\mathbb{P}$ such that $K(X_t,d\xi)$ is the compensator of
$J$ and $\mathbb{P}|_{\mathcal{G}}=\mathbb{P}_1$ holds true. For
this purpose, we shall make use of \cite{jacodj}, Theorem 3.6. Let
us define the following random measure $\nu$ by
\[
\nu(dt,d\xi)=\cases{K(X_t,d\xi)\,dt, &\quad $t < T_{\infty}$,\cr
0, &\quad $t \geq T_{\infty}$.}
\]
Observe that $\nu$ is predictable, since $X_t$ is c\`adl\`ag and
$\mathcal{F}_t$-adapted. Theorem 3.6 in \cite{jacodj} now implies
that there exists a unique probability kernel $\mathbb{P}_2$ from
$\Omega_1$ to~$\mathcal{H}$, such that $\nu$ is the compensator of
the random measure $\mu$ associated to the jumps of $J$. On
$(\Omega, \mathcal{F})$ we then define the probability measure
$\mathbb{P}$ by
$\mathbb{P}(d\omega)=\mathbb{P}_1(d\omega_1)\mathbb{P}_2(\omega
_1,d\omega_2)$
whose restriction to $\mathcal{G}$ is equal to $\mathbb{P}_1$.

\textit{Step} 3. We finally show that $X$ defined by (\ref{Xdef})
solves (\ref{eq: jump diffusion}) on $(\Omega, \mathcal{F},
(\mathcal{F}_t), \mathbb{P})$ for all $t\ge0$. Note that
$W(\omega)=W(\omega_1)$ is an
$(\Omega,\mathcal{F},(\mathcal{F}_t),\mathbb{P})$-Brownian motion.
This implies that $Z(\theta,x,\omega,t)=Z(\theta,x,\omega_1,t)$ is a
solution of (\ref{sde}) on
$(\Omega,\mathcal{F},(\mathcal{F}_t),\mathbb{P})$, satisfying the
properties of Theorem \ref{th: regulardiff}. It thus remains to show
that $T_{\infty}=\infty$ $\Pa$-a.s. Let $\mu$ be the random measure
associated with the jumps. As $x\mapsto K(x,\re^n)$ is bounded, we
have for all $T \geq0$,
\[
\mathbb{E}_{\mathbb{P}}\bigl[\mu([0,T]\times
\re^n)\bigr]=\mathbb{E}_{\mathbb{P}}\bigl[\nu([0,T]\times
\re^n)\bigr]=\mathbb{E}_{\mathbb{P}}\biggl[\int_0^T
K(X_t,\re^n)\,dt\biggr]\leq C T
\]
for some constant $C$. This implies that $\mu([0,T]\times\re^n)<
\infty$ a.s. for all $T \geq0$ and hence $\mathbb{P}[T_{\infty}<
\infty]=0$ or equivalently $T_\infty=\infty$ a.s.
\end{pf}

%s7.3 ###
\subsection{Viability of jump-diffusions}

Consider a nonempty closed convex set $D \subset\re^n$. We now
provide sufficient conditions for the solution $X$ in (\ref{Xdef})
to be $D$-valued. This result is based on \cite{daprato},
Theorem 4.1. We recall the notion of the \textit{normal cone}
%
%e7.8 ###
%
\begin{equation}\label{defnorcone}
N_{D}(x)=\{u \in\re^n \mid\langle u,y-x\rangle\geq0\mbox{,  for
all }y \in D \}
\end{equation}
of $D$ at $x\in D$, consisting of inward pointing vectors. See, for example,
\cite{hiriart}, Definition III.5.2.3, except for a change of the
sign.
\begin{theorem}\label{th: martingale problem convex}
Assume that $\sigma$ also has a Lipschitz continuous derivative.
Suppose furthermore that
%
%e7.11 ###
%e7.10 ###
%e7.9 ###
%
\begin{eqnarray}
\label{eq: cond jumps}
x+\supp(K(x,\cdot)) &\subseteq& D,\\
\label{eq: invariance cond1}
\langle\sigma^i(x),u\rangle&=&0,\\
\label{eq: invariance cond2}
\Biggl\langle b(x)-\frac{1}{2}\sum_{i=1}^n
D\sigma^i(x)\sigma^i(x),u\Biggr\rangle&\geq&0,
\end{eqnarray}
for all $u \in N_D(x)$ and $x \in D$, where $\sigma^i$ denotes the
$i$th column of $\sigma$. Then, for every initial point $x \in
D$, the process $X$ defined in (\ref{Xdef}) is a $D$-valued solution
of (\ref{eq: jump diffusion}).
\end{theorem}
\begin{pf}
We have to show that $X_t=\sum_{n\ge1} X^{(n)}_t(\omega_1,\omega_2)
1_{\{T_{n-1}\le t<T_{n}\}} \in D$ a.s. for all $t\ge0$. We proceed
by induction on $n$. For $n=1$, $X^{(1)}_t$, is simply given by
\[
X^{(1)}_t=x+\int_{0}^{t}b\bigl(X_s^{(1)}\bigr)\,ds+\int_{0}^{t}\sigma\bigl(X_s^{(1)}\bigr)\,dW_s.
\]
Due to \cite{daprato}, Theorem 4.1, conditions (\ref{eq: invariance
cond1}) and (\ref{eq: invariance cond2}) imply that for all \mbox{$t \geq
0$}, $X_t^{(1)} \in D$ a.s. Let us now assume that for all $t \geq
0$, $X^{(n-1)}_t \in D$ a.s., thus in particular
$X^{(n-1)}_{T_{n-1}}=X_{T_{n-1}-} \in D$ a.s. If $T_{n-1}=\infty$,
then we immediately obtain
\[
X^{(n)}_t=X^{(n-1)}_{T_{n-1}}+\Delta J_{T_{n-1}}=x_0\in D.
\]
Otherwise, let $f \in C_b(\re^n, \re_+)$ satisfy $\supp(f)\subseteq
D^c$. Then,
\begin{eqnarray*}
\mathbb{E}\bigl[f\bigl(X^{(n-1)}_{T_{n-1}}+\Delta J_{T_{n-1}}\bigr)\bigr]&=&\mathbb
{E}[f(X_{T_{n-1}-}+\Delta J_{T_{n-1}})]\\
&=&\mathbb{E}\biggl[\int_{\re^n\setminus\{0\}}f(X_{T_{n-1}-}+\xi
)K(X_{T_{n-1}-},d\xi)\biggr]=0,
\end{eqnarray*}
since by (\ref{eq: cond jumps}),
$X_{T_{n-1}-}+\supp(K(X_{T_{n-1}-},\cdot))\subseteq D$ a.s. and
$f(D)=0$. Hence, $f(X^{(n-1)}_{T_{n-1}}+\Delta J_{T_{n-1}})=0$ a.s.,
implying that $X^{(n-1)}_{T_{n-1}}+\Delta J_{T_{n-1}} \notin
\supp(f)$ a.s. As this holds true for all $f\in C_b(\re^n, \re_+)$
with $\supp(f)\subseteq D^c$, it follows that
$X^{(n-1)}_{T_{n-1}}+\Delta J_{T_{n-1}} \in D$ a.s. Thus,\vspace*{1pt} again by
\cite{daprato}, Theorem 4.1, and conditions (\ref{eq: invariance
cond1}) and (\ref{eq: invariance cond2})
\[
X^{(n)}_t=X^{(n-1)}_{T_{n-1}}+\Delta J_{T_{n-1}} +
\int_{0}^{t}b\bigl(X_s^{(n)}\bigr)1_{\{T_{n-1}\le s\}}\,ds
+\int_0^t\sigma\bigl(X^{(n)}_s\bigr)1_{\{T_{n-1}\le s\}} \,dW_s
\]
a.s. takes values in $D$, which proves the induction hypothesis. The
definition of $X$ then yields the assertion.
\end{pf}

%s8 ###
\section{An approximation lemma on the cone of positive semidefinite
matrices}\label{secstoneweier}

In this section, we deliver a differentiable variant of the
Stone--Weierstrass theorem for $C^\infty$-functions on $S_d^+$. This
approximation statement is essential for the description of the
generator of an affine semigroup, as is elaborated in Section
\ref{section: inf generator}.

We employ multi-index notation in the sequel. For $n\geq1$, a
multi-index is an element
${\bolds{\alpha}}=(\alpha_1,\ldots,\alpha_n)\in\mathbb N_0^n$ having
length $\vert{\bolds{\alpha}}\vert:=\alpha_1+\cdots+\alpha_n$. The
factorial is defined by ${\bolds{\alpha}}!:=\prod_{i=1}^n \alpha_i!$.
The partial order $\leq$ is understood componentwise, and so are the
elementary operations $+,-$. That is, ${\bolds{\alpha}}\geq
{\bolds{\beta}}$ if and only if $\alpha_i\geq\beta_i$ for
$i=1,\ldots,n$. Moreover, for
${\bolds{\alpha}}\geq{\bolds{\beta}}$, the multinomial coefficient is
defined by
\[
\pmatrix{\bolds{\alpha}\cr\bolds{\beta}}:=\frac{{\bolds{\alpha
}}!}{(\bolds{\alpha}-\bolds\beta)!{\bolds{\beta}}!}.
\]
We define the monomial $x^{{\bolds{\alpha}}}:=\prod_{i=1}^n
x_i^{\alpha_i}$, and the differential operator
$\partial^{{\bolds{\alpha}}}:=\frac{\partial^{\vert{\bolds
{\alpha}}
\vert}}{\partial_{ x_1}^{\alpha_1}\cdots\partial_{
x_n}^{\alpha_n}}$. Corresponding\vspace*{-2pt} to a polynomial
$P(x)=\sum_{\vert{\bolds{\alpha}}\vert\leq k} a_{\bolds{\alpha}}
x^{\bolds{\alpha}}$, we introduce the differential operator
$P(\partial):=\sum_{\vert{\bolds{\alpha}}\vert\leq k} a_{\bolds
{\alpha}}
\frac{\partial^{\vert{\bolds{\alpha}}\vert}}{\partial
x^{\bolds{\alpha}}}$.\vspace*{1pt}

Let $\Scal=\Scal(S_d)$ denote the locally convex space of rapidly
decreasing $C^\infty$-functions on $S_d$ (see \cite{rudinfana},
Chapter 7), and define the
space of rapidly decreasing $C^\infty$-functions on $S_d^+$ via the
restriction
%
%e8.1 ###
%
\begin{equation}\label{defS+}
\Scal_+=\{ f=F\mid_{S_d^+} \dvtx F\in\Scal\} .
\end{equation}
Equipped with the increasing family of semi-norms
%
%e8.2 ###
%
\begin{equation}\label{eq: seminorm}
p_{k,+}(f):={\sup_{x\in S_d^+,
\vert{\bolds{\alpha}}+{\bolds{\beta}}\vert\leq k}}\vert x^{\bolds
{\alpha}}\,
\partial^{\bolds{\beta}} f(x)\vert,
\end{equation}
$\Scal_+$ becomes a locally convex vector space (see \cite{rudinfana},
Theorem 1.37).

For technical reasons, we also introduce for $\varepsilon\geq0$ the
semi-norms
\[
p_{k,\varepsilon}(f):={\sup_{x\in S_d^+ + B_{\leq
\varepsilon}(0), \vert{\bolds{\alpha}}+{\bolds{\beta}}\vert\leq
k}}\vert
x^{\bolds{\alpha}}
\partial^{\bolds{\beta}} f(x)\vert
\]
on $C^\infty(S_d)$, where $B_{\le r}(y)=\{z \in S_d \mid
\|z-y\|\leq r\}$ denotes the closed ball with radius $r$ and center
$y$. Note that $p_{k,+}=p_{k,0}$. We first give an alternative
description of $\Scal_+$.
\begin{lemma}\label{lemC1}
We have
\[
\mathcal S_+=\{f=G\mid_{S_d^+} \dvtx G\in C^\infty(S_d)\mbox{ and }
\exists\varepsilon>0\mbox{ such that }p_{k,\varepsilon}(G)<\infty
\ \forall k\geq0\}.
\]
\end{lemma}
\begin{pf}
The inclusion $\subseteq$ is trivial. Hence, we prove $\supseteq$. So
let $f=G\mid_{S_d^+}$ for some $G\in C^\infty(S_d)$ with
$p_{k,\varepsilon}(G)<\infty$ for all $k\geq0$ and some
$\varepsilon>0$.

We choose a standard mollifier $\rho\in C^\infty_c(S_d)$ supported
in $B_{\leq\varepsilon/2}(0)$ and satisfying $\rho\geq0$, $\int
\rho=1$. For $\delta>0$ we introduce the neighborhoods
$K_\delta:=S_d^+ + B_{\leq\delta}(0)$ of $S_d^+$. The convolution
$\varphi:=\rho*1_{K_{\varepsilon/2}}\in C^\infty(S_d)$ of the
indicator function for $K_{\varepsilon/2}$ with $\rho$ satisfies
$\varphi=1$ on $S_d^+$ and it vanishes outside $K_{\varepsilon}$.
Furthermore, all derivatives of $\varphi$ are bounded, since
\begin{eqnarray*}
|\partial^{\bolds{\alpha}}\varphi(x)|&=& \biggl|\int_{
K_{\varepsilon/2}}\partial^{\bolds{\alpha}}\rho(y-x)\,dy
\biggr|
=\biggl|\int_{
K_{\varepsilon/2}-x}\partial^{\bolds{\alpha}}\rho(z)\,dz\biggr|\\
&\leq&
\int_{
B_{\leq\varepsilon/2}(0)}|\partial^{\bolds{\alpha}}\rho(z)|\,
dz<\infty,
\end{eqnarray*}
where the last estimate holds because $\supp\rho\subseteq B_{\leq
\varepsilon/2}(0)$.

Now we set $F:=G\cdot\varphi$. By construction $F\in
C^\infty(S_d)$, $F\mid_{S_d^+}=f$ and $F$ vanishes outside
$K_\varepsilon$, because $\varphi$ does. What is left to show is
that $F\in\Scal$. Since $F$ vanishes outside $K_\varepsilon$, it is
sufficient to deliver all estimates of its derivatives on
$K_\varepsilon$.

Let ${\bolds{\alpha}},{\bolds{\beta}}\in\mathbb N_0^{d(d+1)/2}$,
then we
have by the Leibniz rule
\begin{eqnarray*}
x^{\bolds{\alpha}}\,\partial^{\bolds{\beta}} F(x)&=&
x^{\bolds{\alpha}}\sum_{0\leq
{\bolds{\gamma}}\leq{\bolds{\beta}}}\pmatrix{\bolds{\beta}\cr
\bolds{\gamma}}(\partial^{{\bolds{\beta}}-{\bolds{\gamma
}}}\varphi(x)) (\partial^{\bolds{\gamma}}
G(x))\\
&=&\sum_{0\leq{\bolds{\gamma}}\leq{\bolds{\beta}}}
\pmatrix{\bolds{\beta}\cr\bolds{\gamma}}(\partial^{{\bolds
{\beta}}-{\bolds{\gamma}}}\varphi(x)) (x^{\bolds{\alpha}}\,\partial
^{\bolds{\gamma}}
G(x)).
\end{eqnarray*}
By assumption $x^{\bolds{\alpha}}\,\partial^{\bolds{\gamma}} G$ is bounded
on $K_\varepsilon$, and
$(\partial^{{\bolds{\beta}}-{\bolds{\gamma}}}\varphi(x))$ is
bounded on
all of~$S_d$. Hence, by the last equation, we have ${\sup_{x\in S_d
, \vert{\bolds{\alpha}}+{\bolds{\beta}}\vert\leq k}}\vert
x^{\bolds{\alpha}}\,
\partial^{\bolds{\beta}} F(x)\vert<\infty$ for all
\mbox{$k\in\mathbb N_0$}, which by definition means $F\in\mathcal S $.
\end{pf}
\begin{lemma}\label{lemC2}
Let $u\in S_d^{++}$. Then for each $\varepsilon\geq0$, and for all
$k\geq0$ we have $p_{k,\varepsilon}(\exp(-\langle
u,\cdot\rangle))<\infty$. In particular, we have
\[
f_u:=\exp(-\langle u,\cdot\rangle)|_{S_d^+}\in\mathcal S_+.
\]
That
is, $f_u=F_u\mid_{S_d^+}$ for some $F_u\in\mathcal S$.
\end{lemma}
\begin{pf}
Since $u\in S_d^{++}$, there exists a positive constant $c$ such
that $\langle u,x\rangle\geq c\|x\|$, for all $x\in S_d^+$. Hence, we
obtain by a straightforward calculation, $p_{k,+}(\exp(-\langle
u,\cdot\rangle))<\infty$, for all $k\geq0$.

Next, let $\varepsilon>0$, and write $x=y+z$, where $y\in S_d^+$ and
$z\in B_{\leq\varepsilon}(0)$ and pick multi-indices
${\bolds{\alpha}},{\bolds{\beta}}\in\mathbb N_0^{d(d+1)/2}$. Then
we have
by the binomial formula
\begin{eqnarray*}
x^{\bolds{\alpha}}\,\partial^{\bolds{\beta}} e^{-\langle u,x\rangle}&=&
x^{\bolds{\alpha}} (-1)^{\vert\bolds{\beta}\vert}
u^{\bolds{\beta}}e^{-\langle u,x\rangle}\\
&=&(y+z)^{\bolds{\alpha}} (-1)^{\vert\bolds{\beta}\vert}
u^{\bolds{\beta}}e^{-\langle u,y+z\rangle}\\
&=&(-1)^{\vert\bolds{\beta}\vert} u^{\bolds{\beta}} \sum_{0\leq
\bolds{\gamma}\leq\bolds{\alpha}}\pmatrix{\bolds{\alpha}\cr
\bolds{\gamma}}
\bigl(y^{\bolds{\alpha}}e^{-\langle
u,y\rangle}\bigr)\bigl(z^{\bolds{\alpha}-\bolds{\gamma}}e^{-\langle
u,z\rangle}\bigr).
\end{eqnarray*}
Now since $z$ ranges in a compact set, and since
$p_{k,+}(\exp(-\langle u,\cdot\rangle))<\infty$ we see that
$x^{\bolds{\alpha}}\partial^{\bolds{\beta}} e^{-\langle u,x\rangle
}$ must
be bounded uniformly in $x\in S_d^+ +B_{\leq\varepsilon}(0)$. Hence,
$p_{k,\varepsilon}(\exp(-\langle u,\cdot\rangle))<\infty$, for all
$k\geq0$.

Together with Lemma \ref{lemC1}, this implies $f_u \in\mathcal{S}_+$.
\end{pf}

We are now prepared to deliver the following density result for the
$\mathbb R$-linear hull $\mathcal M$ of $\{f_u=\exp(-\langle u,\cdot
\rangle)|_{S_d^+}, u\in S_d^{++}\}$
in $\Scal_+$.
\begin{theorem}\label{th: density}
$\mathcal M$ is dense in $\mathcal S_+$.
\end{theorem}
\begin{pf}
Denote by $\mathcal S' =\mathcal S'(S_d)$ and $\Scal'_+$ the
topological dual of $\Scal$ and $\Scal_+$, respectively. The former,
$\Scal'$, is known as the space of tempered distributions. The
distributional action is denoted by $\langle\cdot,\cdot\rangle$ and
$\langle\cdot,\cdot\rangle_+$ for $\mathcal S'$ and $\mathcal S_+'$,
respectively.

Now suppose by contradiction, that $\mathcal M$ is not dense in
$\mathcal S_+$. Then by \cite{rudinfana}, Theorem 3.5, there exists some
$T_1\in\mathcal S'_+\setminus\{0\}$ such that $T_1=0$ on $\mathcal
M$. Hence, $\langle T_1,f_u\rangle_+=0$, for all $u\in S_d^{++}$. The
restriction $F\mapsto F|_{S_d^+}$ yields a continuous linear
embedding $\mathcal S\hookrightarrow\mathcal S_+$. Hence, the
restriction $T$ of $T_1$ to $\mathcal S$, given by
\[
\langle T,\varphi\rangle:=\langle
T_1,\varphi|_{S_d^+}\rangle_+,\qquad \varphi\in\mathcal S(S_d),
\]
yields an element of $\mathcal S'$ with $\supp(T)\subseteq S_d^+$.
Pick an $F_u\in\mathcal S$ according to Lem\-ma~\ref{lemC2}. By the
definition of $T$, we have $\langle T,F_u\rangle=\langle
T_1,f_u\rangle_+=0$, for all $u\in S_d^{++}$. By the
Bros--Epstein--Glaser theorem (see \cite{reedsimonII}, Theorem IX.15),
there exists a function $G\in C(S_d)$ with $\supp(G)\subseteq
S_d^+$, polynomially bounded [i.e., for suitable constants $C,N$
we have $\vert G(x)\vert\leq C(1+\|x\|)^N$, for all $x\in S_d^+$]
and a real polynomial $P(x)$ such that $P(\partial)G=T$ in $\mathcal
S'$. Hence, we obtain for any $u\in S_d^{++}$
\begin{eqnarray*}
0&=&\langle T,F_u\rangle=\langle P(\partial)G,F_u\rangle=\langle
G,P(-\partial)F_u\rangle\\
&=&\int_{S_d ^+}G(x)P(-\partial)F_u(x)\,
dx=P(u)\int_{S_d^+}G(x)\exp(-\langle u,x\rangle) \,dx.
\end{eqnarray*}
But the last factor is just the Laplace transform of $G$. This
implies $G=0$, hence $T=0$, which in turn implies that $T_1$
vanishes on all of $\mathcal S_+$, a contradiction.
\end{pf}
\end{appendix}

\section*{Acknowledgments}
We thank Martin Keller-Ressel and Alexander Smirnov for discussions and
helpful comments. %Authors Christa Cuchiero and Josef Teichmann
%gratefully acknowledge the support from the FWF-Grant Y328 (START prize
%from the Austrian Science Fund). Authors Damir Filipovi\'c and Eberhard
%Mayerhofer gratefully acknowledge the support from WWTF (Vienna Science
%and Technology Fund) and Swissquote.

% imsref loaded by lrinkeviciute, 2010-10-22 15:49:37
%

\printaddresses


\begin{thebibliography}{52}

%b1 ###
\bibitem{aczel}
%
\begin{bbook}[vtex]
\bauthor{\bsnm{Acz{\'e}l},~\bfnm{J.}\binits{J.}} \AND
\bauthor{\bsnm{Dhombres},~\bfnm{J.}\binits{J.}}
(\byear{1989}).
\btitle{Functional Equations in Several Variables}.
\bseries{Encyclopedia of Mathematics and Its Applications}
\bvolume{31}.
\bpublisher{Cambridge Univ. Press}, \baddress{Cambridge}.
%and social sciences}.
\bid{mr={1004465}}
\end{bbook}
%
\endbibitem

%b2 ###
\bibitem{bns}
%
\begin{bincollection}[mr]
\bauthor{\bsnm{Barndorff-Nielsen},~\bfnm{Ole~E.}\binits{O.~E.}}
\AND
\bauthor{\bsnm{Shephard},~\bfnm{Neil}\binits{N.}}
(\byear{2001}).
\btitle{Modelling by {L}\'evy processes for financial econometrics}.
In \bbooktitle{L\'evy Processes}
\bpages{283--318}.
\bpublisher{Birkh\"auser}, \baddress{Boston, MA}.
\bid{mr={1833702}}%
\end{bincollection}%
%
\endbibitem%

%b3 ###
\bibitem{barndorffstelzer}
%
\begin{barticle}[mr]
\bauthor{\bsnm{Barndorff-Nielsen},~\bfnm{Ole~Eiler}\binits{O.~E.}}
\AND
\bauthor{\bsnm{Stelzer},~\bfnm{Robert}\binits{R.}}
(\byear{2007}).
\btitle{Positive-definite matrix processes of finite variation}.
\bjournal{Probab. Math. Statist.}
\bvolume{27}
\bpages{3--43}.
\bid{mr={2353270}}
\end{barticle}
%
\endbibitem

%b4 ###
\bibitem{bau96}
%
\begin{bbook}[vtex]
\bauthor{\bsnm{Bauer},~\bfnm{Heinz}\binits{H.}}
(\byear{1996}).
\btitle{Probability Theory}.
\bseries{de Gruyter Studies in Mathematics}
\bvolume{23}.
\bpublisher{de Gruyter}, \baddress{Berlin}.
%and revised by the author}.
\bid{mr={1385460}}
\end{bbook}
%
\endbibitem

%b5 ###
\bibitem{bru89}
%
\begin{barticle}[mr]
\bauthor{\bsnm{Bru},~\bfnm{Marie-France}\binits{M.-F.}}
(\byear{1989}).
\btitle{Diffusions of perturbed principal component analysis}.
\bjournal{J. Multivariate Anal.}
\bvolume{29}
\bpages{127--136}.
\bid{doi={10.1016/0047-259X(89)90080-8}, mr={0991060}}
\end{barticle}
%
\endbibitem

%b6 ###
\bibitem{bru}
%
\begin{barticle}[mr]
\bauthor{\bsnm{Bru},~\bfnm{Marie-France}\binits{M.-F.}}
(\byear{1991}).
\btitle{Wishart processes}.
\bjournal{J. Theoret. Probab.}
\bvolume{4}
\bpages{725--751}.
\bid{doi={10.1007/BF01259552}, mr={1132135}}
\end{barticle}
%
\endbibitem

%b7 ###
\bibitem{burcietro07}
%
\begin{bmisc}[vtex]
\bauthor{\bsnm{Buraschi},~\bfnm{B.}\binits{B.}},
\bauthor{\bsnm{Cieslak},~\bfnm{A.}\binits{A.}} \AND
\bauthor{\bsnm{Trojani},~\bfnm{F.}\binits{F.}}
(\byear{2007}).
\bhowpublished{Correlation risk and the term structure of interest rates.
Working paper, Univ. St. Gallen}.
\end{bmisc}
%
\endbibitem

%b8 ###
\bibitem{buraschiporchiatrojani}
%
\begin{barticle}[vtex]
\bauthor{\bsnm{Buraschi},~\bfnm{B.}\binits{B.}},
\bauthor{\bsnm{Porchia},~\bfnm{P.}\binits{P.}} \AND
\bauthor{\bsnm{Trojani},~\bfnm{F.}\binits{F.}}
(\byear{2010}).
\btitle{Correlation risk and optimal portfolio choice}.
\bjournal{J. Finance}
\bvolume{65}
\bpages{393--420}.
\end{barticle}
%
\endbibitem

%b9 ###
\bibitem{fonsecagrasselliielpo1}
%
\begin{bmisc}[vtex]
\bauthor{\bsnm{Da~Fonseca},~\bfnm{J.}\binits{J.}},
\bauthor{\bsnm{Grasselli},~\bfnm{M.}\binits{M.}} \AND
\bauthor{\bsnm{Ielpo},~\bfnm{F.}\binits{F.}}
(\byear{2008}).
\bhowpublished{Hedging (co)variance risk with variance swaps. Working paper
ESILV RR-37, Ecole Sup\'{e}rieure d'Ing\'{e}nierie
L\'{e}onard de Vinci}.
\end{bmisc}
%
\endbibitem

%b10 ###
\bibitem{fonsecagrasselliielpo2}
%
\begin{bmisc}[vtex]
\bauthor{\bsnm{Da~Fonseca},~\bfnm{J.}\binits{J.}},
\bauthor{\bsnm{Grasselli},~\bfnm{M.}\binits{M.}} \AND
\bauthor{\bsnm{Ielpo},~\bfnm{F.}\binits{F.}}
(\byear{2008}).
\bhowpublished{Estimating the Wishart affine stochastic correlation
model using
the empirical characteristic function. Working paper ESILV RR-35, Ecole Sup\'{e}rieure d'Ing\'{e}nierie
L\'{e}onard de Vinci}.
\end{bmisc}
%
\endbibitem

%b11 ###
\bibitem{fonsecaetal2}
%
\begin{barticle}[auto:SpringerTagBib|2010-03-24|17:41:21]
\bauthor{\bsnm{Da~Fonseca},~\bfnm{J.}\binits{J.}},
\bauthor{\bsnm{Grasselli},~\bfnm{M.}\binits{M.}} \AND
\bauthor{\bsnm{Tebaldi},~\bfnm{C.}\binits{C.}}
(\byear{2007}).
\btitle{Option pricing when correlations are stochastic: An analytical
framework}.
\bjournal{Review of Derivatives Research}
\bvolume{10}
\bpages{151--180}.
\end{barticle}
%
\endbibitem

%b12 ###
\bibitem{fonsecaetal1}
%
\begin{barticle}[mr]
\bauthor{\bsnm{Da~Fonseca},~\bfnm{Jos{\'e}}\binits{J.}},
\bauthor{\bsnm{Grasselli},~\bfnm{Martino}\binits{M.}} \AND
\bauthor{\bsnm{Tebaldi},~\bfnm{Claudio}\binits{C.}}
(\byear{2008}).
\btitle{A multifactor volatility {H}eston model}.
\bjournal{Quant. Finance}
\bvolume{8}
\bpages{591--604}.
\bid{doi={10.1080/14697680701668418}, mr={2457710}}
\end{barticle}
%
\endbibitem

%b13 ###
\bibitem{daprato}
%
\begin{barticle}[mr]
\bauthor{\bsnm{Da~Prato},~\bfnm{Giuseppe}\binits{G.}} \AND
\bauthor{\bsnm{Frankowska},~\bfnm{H{\'e}l{\`e}ne}\binits{H.}}
(\byear{2004}).
\btitle{Invariance of stochastic control systems with deterministic arguments}.
\bjournal{J. Differential Equations}
\bvolume{200}
\bpages{18--52}.
\bid{doi={10.1016/j.jde.2004.01.007}, mr={2046316}}
\end{barticle}
%
\endbibitem

%b14 ###
\bibitem{dieu69}
%
\begin{bbook}[vtex]
\bauthor{\bsnm{Dieudonn{\'e}},~\bfnm{J.}\binits{J.}}
(\byear{1969}).
\btitle{Foundations of Modern Analysis}.
\bpublisher{Academic Press}, \baddress{New York}.
%10-I}.
\bid{mr={0349288}}
\end{bbook}
%
\endbibitem

%b15 ###
\bibitem{donatimartin}
%
\begin{barticle}[vtex]
\bauthor{\bsnm{Donati-Martin},~\bfnm{Catherine}\binits{C.}},
\bauthor{\bsnm{Doumerc},~\bfnm{Yan}\binits{Y.}},
\bauthor{\bsnm{Matsumoto},~\bfnm{Hiroyuki}\binits{H.}} \AND
\bauthor{\bsnm{Yor},~\bfnm{Marc}\binits{M.}}
(\byear{2004}).
\btitle{Some properties of the {W}ishart processes and a matrix
extension of
the {H}artman--{W}atson laws}.
\bjournal{Publ. Res. Inst. Math. Sci.}
\bvolume{40}
\bpages{1385--1412}.
\bid{mr={2105711}}
\end{barticle}
%
\endbibitem

%b16 ###
\bibitem{dfs}
%
\begin{barticle}[mr]
\bauthor{\bsnm{Duffie},~\bfnm{D.}\binits{D.}},
\bauthor{\bsnm{Filipovi{\'c}},~\bfnm{D.}\binits{D.}} \AND
\bauthor{\bsnm{Schachermayer},~\bfnm{W.}\binits{W.}}
(\byear{2003}).
\btitle{Affine processes and applications in finance}.
\bjournal{Ann. Appl. Probab.}
\bvolume{13}
\bpages{984--1053}.
\bid{doi={10.1214/aoap/1060202833}, mr={1994043}}
\end{barticle}
%
\endbibitem

%b17 ###
\bibitem{ethier}
%
\begin{bbook}[vtex]
\bauthor{\bsnm{Ethier},~\bfnm{Stewart~N.}\binits{S.~N.}} \AND
\bauthor{\bsnm{Kurtz},~\bfnm{Thomas~G.}\binits{T.~G.}}
(\byear{1986}).
\btitle{Markov Processes: Characterization and Convergence}.
%and Mathematical Statistics}.
\bpublisher{Wiley}, \baddress{New York}.
\bid{doi={10.1002/9780470316658}, mr={0838085}}
\end{bbook}
%
\endbibitem

%b18 ###
\bibitem{filaff05}
%
\begin{barticle}[mr]
\bauthor{\bsnm{Filipovi{\'c}},~\bfnm{Damir}\binits{D.}}
(\byear{2005}).
\btitle{Time-inhomogeneous affine processes}.
\bjournal{Stochastic Process. Appl.}
\bvolume{115}
\bpages{639--659}.
\bid{doi={10.1016/j.spa.2004.11.006}, mr={2128634}}
\end{barticle}
%
\endbibitem

%b19 ###
\bibitem{fil09}
%
\begin{bbook}[vtex]
\bauthor{\bsnm{Filipovi{\'c}},~\bfnm{Damir}\binits{D.}}
(\byear{2009}).
\btitle{Term-Structure Models: A Graduate Course}.
\bpublisher{Springer}, \baddress{Berlin}.
\bid{doi={10.1007/978-3-540-68015-4}, mr={2553163}}
\end{bbook}
%
\endbibitem

%b20 ###
\bibitem{ADPTA}
%
\begin{bincollection}[mr]
\bauthor{\bsnm{Filipovi{\'c}},~\bfnm{Damir}\binits{D.}} \AND
\bauthor{\bsnm{Mayerhofer},~\bfnm{Eberhard}\binits{E.}}
(\byear{2009}).
\btitle{Affine diffusion processes: Theory and applications}.
In \bbooktitle{Advanced Financial Modelling}.
\bseries{Radon Ser. Comput. Appl. Math.}
\bvolume{8}
\bpages{125--164}.
\bpublisher{Walter de Gruyter, Berlin}.
\bid{doi={10.1515/9783110213140.125}, mr={2648460}}
\end{bincollection}
%
\endbibitem

%b21 ###
\bibitem{filipovicoverbeck}
%
\begin{bmisc}[vtex]
\bauthor{\bsnm{Filipovi\'c},~\bfnm{D.}\binits{D.}},
\bauthor{\bsnm{Overbeck},~\bfnm{L.}\binits{L.}} \AND
\bauthor{\bsnm{Schmidt},~\bfnm{T.}\binits{T.}}
(\byear{2009}).
\bhowpublished{Dynamic CDO term structure modeling. Forthcoming in mathematical
finance, Ecole Polytechnique F\'{e}d\'{e}rale de Lausanne}.
\end{bmisc}
%
\endbibitem

%b22 ###
\bibitem{golub}
%
\begin{bbook}[vtex]
\bauthor{\bsnm{Golub},~\bfnm{Gene~H.}\binits{G.~H.}} \AND
\bauthor{\bsnm{Van~Loan},~\bfnm{Charles~F.}\binits{C.~F.}}
(\byear{1996}).
\btitle{Matrix Computations},
\bedition{3rd} ed.
\bpublisher{Johns Hopkins Univ. Press}, \baddress{Baltimore, MD}.
\bid{mr={1417720}}
\end{bbook}
%
\endbibitem

%b23 ###
\bibitem{gourieerouxmonfortsufana}
%
\begin{bmisc}[vtex]
\bauthor{\bsnm{Gourieroux},~\bfnm{C.}\binits{C.}},
\bauthor{\bsnm{Montfort},~\bfnm{A.}\binits{A.}} \AND
\bauthor{\bsnm{Sufana},~\bfnm{R.}\binits{R.}}
(\byear{2007}).
\bhowpublished{International money and stock market contingent claims. Working
paper, CREST, CEPREMAP and Univ. Toronto}.
\end{bmisc}
%
\endbibitem

%b24 ###
\bibitem{gourierouxsufana}
%
\begin{bmisc}[vtex]
\bauthor{\bsnm{Gourieroux},~\bfnm{C.}\binits{C.}} \AND
\bauthor{\bsnm{Sufana},~\bfnm{R.}\binits{R.}}
(\byear{2007}).
\bhowpublished{Wishart quadratic term structure models. Working paper, CREST,
CEPREMAP and Univ. Toronto}.
\end{bmisc}
%
\endbibitem

%b25 ###
\bibitem{gousuf04}
%
\begin{bmisc}[vtex]
\bauthor{\bsnm{Gourieroux},~\bfnm{C.}\binits{C.}} \AND
\bauthor{\bsnm{Sufana},~\bfnm{R.}\binits{R.}}
(\byear{2007}).
\bhowpublished{Derivative pricing with Wishart multivariate stochastic
volatility: Application to credit risk. Working paper, CREST, CEPREMAP and
Univ. Toronto}.
\end{bmisc}
%
\endbibitem

%b26 ###
\bibitem{grasselli}
%
\begin{barticle}[mr]
\bauthor{\bsnm{Grasselli},~\bfnm{Martino}\binits{M.}} \AND
\bauthor{\bsnm{Tebaldi},~\bfnm{Claudio}\binits{C.}}
(\byear{2008}).
\btitle{Solvable affine term structure models}.
\bjournal{Math. Finance}
\bvolume{18}
\bpages{135--153}.
\bid{doi={10.1111/j.1467-9965.2007.00325.x}, mr={2380943}}
\end{barticle}
%
\endbibitem

%b27 ###
\bibitem{heston}
%
\begin{barticle}[auto:SpringerTagBib|2010-03-24|17:41:21]
\bauthor{\bsnm{Heston},~\bfnm{S.}\binits{S.}}
(\byear{1993}).
\btitle{A closed-form solution for options with stochastic volatility with
applications to bond and currency options}.
\bjournal{Rev. of Financial Studies}
\bvolume{6}
\bpages{327--343}.
\end{barticle}
%
\endbibitem

%b28 ###
\bibitem{hiriart}
%
\begin{bbook}[vtex]
\bauthor{\bsnm{Hiriart-Urruty},~\bfnm{Jean-Baptiste}\binits{J.-B.}}
\AND
\bauthor{\bsnm{Lemar{\'e}chal},~\bfnm{Claude}\binits{C.}}
(\byear{1993}).
\btitle{Convex Analysis and Minimization Algorithms. {I}}.
\bseries{Grundlehren der Mathematischen Wissenschaften [Fundamental Principles
of Mathematical Sciences]}
\bvolume{305}.
\bpublisher{Springer}, \baddress{Berlin}.
\bid{mr={1261420}}
\end{bbook}
%
\endbibitem

%b29 ###
\bibitem{horn}
%
\begin{bbook}[mr]
\bauthor{\bsnm{Horn},~\bfnm{Roger~A.}\binits{R.~A.}} \AND
\bauthor{\bsnm{Johnson},~\bfnm{Charles~R.}\binits{C.~R.}}
(\byear{1991}).
\btitle{Topics in Matrix Analysis}.
\bpublisher{Cambridge Univ. Press}, \baddress{Cambridge}.
\bid{mr={1091716}}
\end{bbook}
%
\endbibitem

%b30 ###
\bibitem{jacodj}
%
\begin{barticle}[vtex]
\bauthor{\bsnm{Jacod},~\bfnm{Jean}\binits{J.}}
(\byear{1974/75}).
\btitle{Multivariate point processes: Predictable projection,
{R}adon--{N}ikod\'ym derivatives, representation of martingales}.
\bjournal{Z. Wahrsch. Verw. Gebiete}
\bvolume{31}
\bpages{235--253}.
\bid{mr={0380978}}
\end{barticle}
%
\endbibitem

%b31 ###
\bibitem{jacod}
%
\begin{bbook}[mr]
\bauthor{\bsnm{Jacod},~\bfnm{Jean}\binits{J.}} \AND
\bauthor{\bsnm{Shiryaev},~\bfnm{Albert~N.}\binits{A.~N.}}
(\byear{2003}).
\btitle{Limit Theorems for Stochastic Processes},
\bedition{2nd} ed.
\bseries{Grundlehren der Mathematischen Wissenschaften [Fundamental Principles
of Mathematical Sciences]}
\bvolume{288}.
\bpublisher{Springer}, \baddress{Berlin}.
\bid{mr={1943877}}
\end{bbook}
%
\endbibitem

%b32 ###
\bibitem{karatzas}
%
\begin{bbook}[mr]
\bauthor{\bsnm{Karatzas},~\bfnm{Ioannis}\binits{I.}} \AND
\bauthor{\bsnm{Shreve},~\bfnm{Steven~E.}\binits{S.~E.}}
(\byear{1991}).
\btitle{Brownian Motion and Stochastic Calculus},
\bedition{2nd} ed.
\bseries{Graduate Texts in Mathematics}
\bvolume{113}.
\bpublisher{Springer}, \baddress{New York}.
\bid{mr={1121940}}
\end{bbook}
%
\endbibitem

%b33 ###
\bibitem{kat95}
%
\begin{bbook}[vtex]
\bauthor{\bsnm{Kato},~\bfnm{Tosio}\binits{T.}}
(\byear{1995}).
\btitle{Perturbation Theory for Linear Operators}.
\bpublisher{Springer}, \baddress{Berlin}.
\bid{mr={1335452}}
\end{bbook}
%
\endbibitem

%b34 ###
\bibitem{keller}
%
\begin{bmisc}[vtex]
\bauthor{\bsnm{Keller-Ressel},~\bfnm{M.}\binits{M.}}
(\byear{2009}).
\bhowpublished{Affine processes---theory and applications in mathematical
finance. Ph.D. thesis, Vienna Univ. Technology}.
\end{bmisc}
%
\endbibitem

%b35 ###
\bibitem{kst}
%
\begin{bmisc}[vtex]
\bauthor{\bsnm{Keller-Ressel},~\bfnm{M.}\binits{M.}},
\bauthor{\bsnm{Schachermayer},~\bfnm{W.}\binits{W.}} \AND
\bauthor{\bsnm{Teichmann},~\bfnm{J.}\binits{J.}}
(\byear{2010}).
\bhowpublished{Affine processes are regular.
\textit{Probab. Theory Related Fields}.}
DOI:
\href{http://dx.doi.org/10.1007/s00440-010-0309-4}{10.1007/s00440-010-0309-4}.
\bnote{To appear}.
\end{bmisc}
%
\endbibitem

%b36 ###
\bibitem{lamberton}
%
\begin{bbook}[vtex]
\bauthor{\bsnm{Lamberton},~\bfnm{Damien}\binits{D.}} \AND
\bauthor{\bsnm{Lapeyre},~\bfnm{Bernard}\binits{B.}}
(\byear{2008}).
\btitle{Introduction to Stochastic Calculus Applied to Finance},
\bedition{2nd} ed.
\bpublisher{Chapman and Hall/CRC}, \baddress{Boca Raton, FL}.
\bid{mr={2362458}}
\end{bbook}
%
\endbibitem

%b37 ###
\bibitem{lang}
%
\begin{bbook}[mr]
\bauthor{\bsnm{Lang},~\bfnm{Serge}\binits{S.}}
(\byear{1993}).
\btitle{Real and Functional Analysis},
\bedition{3rd} ed.
\bseries{Graduate Texts in Mathematics}
\bvolume{142}.
\bpublisher{Springer}, \baddress{New York}.
\bid{mr={1216137}}
\end{bbook}
%
\endbibitem

%b38 ###
\bibitem{leippoldtrojani}
%
\begin{bmisc}[vtex]
\bauthor{\bsnm{Leippold},~\bfnm{M.}\binits{M.}} \AND
\bauthor{\bsnm{Trojani},~\bfnm{F.}\binits{F.}}
(\byear{2008}).
\bhowpublished{Asset pricing with matrix affine jump diffusions.
Working paper, University of Zurich---Swiss Banking
Institute (ISB)}.
\end{bmisc}
%
\endbibitem

%b39 ###
\bibitem{maysmi09}
%
\begin{bmisc}[vtex]
\bauthor{\bsnm{Mayerhofer},~\bfnm{E.}\binits{E.}},
\bauthor{\bsnm{Muhle-Karbe},~\bfnm{J.}\binits{J.}} \AND
\bauthor{\bsnm{Smirnov},~\bfnm{A.~G.}\binits{A.~G.}}
(\byear{2011}).
\bhowpublished{A characterization of the martingale property of exponentially
affine processes.
\textit{Stochastic Process. Appl.}
\textbf{121} 568--582}.
\end{bmisc}
%
\endbibitem

%b40 ###
\bibitem{pfaffel}
%
\begin{bmisc}[vtex]
\bauthor{\bsnm{Mayerhofer},~\bfnm{E.}\binits{E.}},
\bauthor{\bsnm{Pfaffel},~\bfnm{O.}\binits{O.}} \AND
\bauthor{\bsnm{Stelzer},~\bfnm{R.}\binits{R.}}
(\byear{2009}).
\bhowpublished{On strong solutions for positive definite
jump--diffusions. VIF
Working Paper No. 30, Vienna Institute of Finance}.
\end{bmisc}
%
\endbibitem

%b41 ###
\bibitem{Narasimhan71}
%
\begin{bbook}[vtex]
\bauthor{\bsnm{Narasimhan},~\bfnm{Raghavan}\binits{R.}}
(\byear{1971}).
\btitle{Several Complex Variables}.
\bpublisher{Univ. Chicago Press},
\baddress{Chicago}.
\bid{mr={0342725}}
\end{bbook}
%
\endbibitem

%b42 ###
\bibitem{reedsimonII}
%
\begin{bbook}[vtex]
\bauthor{\bsnm{Reed},~\bfnm{Michael}\binits{M.}} \AND
\bauthor{\bsnm{Simon},~\bfnm{Barry}\binits{B.}}
(\byear{1975}).
\btitle{Methods of Modern Mathematical Physics. {II}. {F}ourier Analysis,
Self-Adjointness}.
\bpublisher{Academic Press},
\baddress{New York}.
\bid{mr={0493420}}
\end{bbook}
%
\endbibitem

%b43 ###
\bibitem{revuzyor}
%
\begin{bbook}[mr]
\bauthor{\bsnm{Revuz},~\bfnm{Daniel}\binits{D.}} \AND
\bauthor{\bsnm{Yor},~\bfnm{Marc}\binits{M.}}
(\byear{1991}).
\btitle{Continuous Martingales and {B}rownian Motion}.
\bseries{Grundlehren der Mathematischen Wissenschaften [Fundamental Principles
of Mathematical Sciences]}
\bvolume{293}.
\bpublisher{Springer}, \baddress{Berlin}.
\bid{mr={1083357}}
\end{bbook}
%
\endbibitem

%b44 ###
\bibitem{roc97}
%
\begin{bbook}[vtex]
\bauthor{\bsnm{Rockafellar},~\bfnm{R.~Tyrrell}\binits{R.~T.}}
(\byear{1997}).
\btitle{Convex Analysis}.
\bpublisher{Princeton Univ. Press}, \baddress{Princeton, NJ}.
\bid{mr={1451876}}
\end{bbook}
%
\endbibitem

%b45 ###
\bibitem{rogers}
%
\begin{bbook}[vtex]
\bauthor{\bsnm{Rogers},~\bfnm{L.~C.~G.}\binits{L.~C.~G.}} \AND
\bauthor{\bsnm{Williams},~\bfnm{David}\binits{D.}}
(\byear{2000}).
\btitle{Diffusions, {M}arkov Processes, and Martingales. {V}ol. 2}.
\bpublisher{Cambridge Univ. Press}, \baddress{Cambridge}.
\bid{mr={1780932}}
\end{bbook}
%
\endbibitem

%b46 ###
\bibitem{rudinfana}
%
\begin{bbook}[vtex]
\bauthor{\bsnm{Rudin},~\bfnm{Walter}\binits{W.}}
(\byear{1991}).
\btitle{Functional Analysis},
\bedition{2nd} ed.
\bpublisher{McGraw-Hill}, \baddress{New York}.
\bid{mr={1157815}}
\end{bbook}
%
\endbibitem

%b47 ###
\bibitem{sato}
%
\begin{bbook}[vtex]
\bauthor{\bsnm{Sato},~\bfnm{Ken-iti}\binits{K.-i.}}
(\byear{1999}).
\btitle{L\'evy Processes and Infinitely Divisible Distributions}.
\bseries{Cambridge Studies in Advanced Mathematics}
\bvolume{68}.
\bpublisher{Cambridge Univ. Press}, \baddress{Cambridge}.
\bid{mr={1739520}}
\end{bbook}
%
\endbibitem

%b48 ###
\bibitem{semadeni}
%
\begin{bbook}[vtex]
\bauthor{\bsnm{Semadeni},~\bfnm{Zbigniew}\binits{Z.}}
(\byear{1971}).
\btitle{Banach Spaces of Continuous Functions. {V}ol. {I}}.
\bpublisher{PWN---Polish Scientific Publishers}, \baddress{Warsaw}.
\bid{mr={0296671}}
\end{bbook}
%
\endbibitem

%b49 ###
\bibitem{skorohod}
%
\begin{bbook}[vtex]
\bauthor{\bsnm{Skorohod},~\bfnm{A.~V.}\binits{A.~V.}}
(\byear{1991}).
\btitle{Random Processes with Independent Increments}.
\bseries{Mathematics and Its Applications (Soviet Series)}
\bvolume{47}.
\bpublisher{Kluwer Academic}, \baddress{Dordrecht}.
\bid{mr={1155400}}
\end{bbook}
%
\endbibitem

%b50 ###
\bibitem{stokes}
%
\begin{barticle}[mr]
\bauthor{\bsnm{Stokes},~\bfnm{A.~N.}\binits{A.~N.}}
(\byear{1974}).
\btitle{A special property of the matrix {R}iccati equation}.
\bjournal{Bull. Austral. Math. Soc.}
\bvolume{10}
\bpages{245--253}.
\bid{mr={0342748}}
\end{barticle}
%
\endbibitem

%b51 ###
\bibitem{stroock1}
%
\begin{barticle}[mr]
\bauthor{\bsnm{Stroock},~\bfnm{Daniel~W.}\binits{D.~W.}}
(\byear{1975}).
\btitle{Diffusion processes associated with {L}\'evy generators}.
\bjournal{Z. Wahrsch. Verw. Gebiete}
\bvolume{32}
\bpages{209--244}.
\bid{mr={0433614}}
\end{barticle}
%
\endbibitem

%b52 ###
\bibitem{Volkmann1973}
%
\begin{barticle}[mr]
\bauthor{\bsnm{Volkmann},~\bfnm{Peter}\binits{P.}}
(\byear{1973}).
\btitle{\"{U}ber die {I}nvarianz konvexer {M}engen und
{D}ifferentialungleichungen in einem normierten {R}aume}.
\bjournal{Math. Ann.}
\bvolume{203}
\bpages{201--210}.
\bid{mr={0322305}}
\end{barticle}
%
\endbibitem

\end{thebibliography}
\end{document}